\documentclass[final,11pt,3p,doublespacing]{elsarticle/elsarticle}
\journal{journal for consideration.}
\usepackage{times}
\usepackage{booktabs}
\emergencystretch=\maxdimen
\usepackage{titlesec}  
\usepackage{sectsty}
\sectionfont{\color{blue}\large}
\subsectionfont{\color{blue}}
\subsubsectionfont{\color{blue}\itshape}
\usepackage{algorithm}
\usepackage{algpseudocode}
\usepackage{amsmath, amssymb}
\biboptions{sort&compress}
\usepackage[colorlinks]{hyperref}
\AtBeginDocument{%
	\hypersetup{
		plainpages = false, 
		bookmarksopen = true,
		bookmarksnumbered = true,
		breaklinks = true,
		linktocpage,
		colorlinks = true,
		linkcolor = purple,
		urlcolor  = black,
		citecolor = Magenta,
}}
\usepackage{multicol}
\usepackage{multirow}
\usepackage{imakeidx}
\makeindex
\usepackage{nomencl}
\makenomenclature
\usepackage[xindy]{glossaries}
\makeglossaries
\usepackage{graphics}
\usepackage{graphicx}
\usepackage{epsf,psfrag}
\usepackage{epstopdf}
\usepackage[para]{threeparttable}
\usepackage{subfigure}
\usepackage{epsfig}
\usepackage{pdflscape}
\usepackage{rotating}
\usepackage{lscape}
\usepackage{float}
\usepackage{stfloats}
\usepackage{textgreek}
\usepackage{longtable}
\usepackage{lscape}
\usepackage{capt-of}
\usepackage{tablefootnote}
\usepackage{footnote}
\usepackage{scrextend}
\usepackage{caption}
\usepackage[autopunct=true]{csquotes}
\usepackage{rotating}
\usepackage[normalem]{ulem}
\usepackage{resizegather}
\usepackage{siunitx}
\usepackage{cases}
\usepackage[nodisplayskipstretch]{setspace}
\usepackage{color,soul}
\usepackage[usenames,dvipsnames]{xcolor}
\usepackage{amsfonts,amsmath,amssymb,gensymb,stmaryrd}
\usepackage{latexsym,array}
\usepackage{textcomp}

\newcommand{\RomanNumeralCaps}[1]
\linenumbers
\makesavenoteenv{tabular}
\makesavenoteenv{table}
\usepackage[left,displaymath, mathlines]{lineno}
\DeclareMathAlphabet{\mathpzc}{OT1}{pzc}{m}{it}
\def\fig{Fig.~}
\def\figs{Figs.~}
\def\eqn{Eq.~}
\def\eqns{Eqs.~}
\def\tab{Table~}
\def\tabs{Tables~}
\def\sect{Section~}

\def\tsc#1{\csdef{#1}{\textsc{\lowercase{#1}}\xspace}}
\tsc{WGM}
\tsc{QE}
\tsc{EP}
\tsc{PMS}
\tsc{BEC}
\tsc{DE}
\usepackage{easyReview} \setreviewsoff
\newcommand{\removeEq}[1]{\ifistoreview{\@\expandafter\removeColor{#1}\hspace{-0.6em}} \else {}\fi} 

\graphicspath{{../}{Figures/}{Figures/Validation/}}
\DeclareGraphicsExtensions{.png, .PNG,.pdf,.PDF,.jpg,.JPG,.eps,.EPS}
\begin{document}
\setcounter{page}{1}
\begin{frontmatter} 
%
%
%
%
\title{Efficient Cell-Centered Nodal Integral Method for Multi-Dimensional Burgers' Equations}
\author[labela]{Nadeem Ahmed}
\emailauthor{26ahmedansari@gmail.com}{NA}
\author[labela]{Ram Prakash {Bharti}\corref{coradd}}
\emailauthor{rpbharti@iitr.ac.in}{RPB}
\author[labelb]{Suneet Singh}
\emailauthor{suneet.singh@iitb.ac.in}{SS}
\address[labela]{Complex Fluid Dynamics and Microfluidics (CFDM) Lab, Department of Chemical Engineering, Indian Institute of Technology Roorkee, Roorkee - 247667, Uttarakhand, India}
\address[labelb]{Fluid Flow and System Simulation Lab, Department of Energy Science and Engineering, Indian Institute of Technology Bombay, Mumbai 400076, India}
%
%
\cortext[coradd]{\textit{Corresponding author. }}
%
\begin{abstract}
\fontsize{11}{16pt}\selectfont
%
An efficient coarse-mesh nodal integral method (NIM), based on cell-centered variables and termed the cell-centered NIM (CCNIM), is developed and applied to solve multi-dimensional, time-dependent, nonlinear Burgers’ equations, extending the applicability of CCNIM to nonlinear problems. \add{To overcome the existing limitation of CCNIM to linear problems, the convective velocity in the nonlinear convection term is approximated using two different approaches, both demonstrating accuracy comparable to or better than traditional NIM for nonlinear Burgers’ problems.} Unlike traditional NIM, which utilizes surface-averaged variables as discrete unknowns, this innovative approach formulates the final expression of the numerical scheme using discrete unknowns represented by cell-centered (or node-averaged) variables. Using these cell centroids, the proposed CCNIM approach presents several advantages compared to traditional NIM. These include a simplified implementation process in terms of local coordinate systems, enhanced flexibility regarding the higher order of accuracy in time, straightforward formulation for higher-degree temporal derivatives, and offering a viable option for coupling with other physics. The multi-dimensional time-dependent Burgers' problems (propagating shock, propagation, and diffusion of an initial sinusoidal wave, shock-like formation) with known analytical solutions are solved in order to validate the developed scheme. Furthermore, a detailed comparison between the proposed CCNIM approach and other traditional NIM schemes is conducted to demonstrate its effectiveness. \add{The proposed approach has shown quadratic convergence in both space and time, i.e., O[$(\Delta x)^2, (\Delta t)^2$], for the considered test problems. }The simplicity and robustness of the approach provide a strong foundation for its seamless extension to more complex fluid flow problems. 
\end{abstract}
\end{frontmatter}
\section{Introduction}
\label{sec:1}
Over the past few decades, numerous nodal schemes have been developed to solve complex partial differential equations (PDEs) in various scientific and engineering applications \citep{01Delp_1964,02Steinke_1973,03Burns_1975,04Lawrence_1980,05Dorning_1979,06Hennart_1986,07Lawrence_1986}. These methods are specifically designed to address the computational challenges associated with large-scale computations, such as extensive memory requirements and high computational expense \citep{05Dorning_1979,07Lawrence_1986}. One of the main advantages of nodal schemes is their ability to produce numerical solutions with accuracy comparable to that of more conventional approaches but in significantly less computational time \citep{05Dorning_1979}. This efficiency is primarily achieved by using a coarser mesh size, which reduces the number of computational elements, and consequently, the overall computational workload. For a given mesh size, a second-order coarse-mesh scheme typically results in smaller numerical errors compared to a second-order finite-difference scheme \citep{04Lawrence_1980,11Azmy_1983,azmy1985nodal}. This is because the nodal schemes include the analytical pre-processing (or local analytical solution) at each node in the development procedure \citep{11Azmy_1983,esser1993upwind,22Rizwan_uddin_1997}. This makes them particularly useful for large-scale problems where computational resources are limited or where quick turnaround times are critical.

Various approaches to nodal methods have been explored in the literature to solve the neutron diffusion equation \citep{01Delp_1964,02Steinke_1973,03Burns_1975,09Shober_1977,05Dorning_1979,04Lawrence_1980,11Azmy_1983,06Hennart_1986,07Lawrence_1986,13Shober1978,08Raj_2022,10Ferrer_2009,wang2010improved,12Raj_2017,guessous2016three}. Some of these approaches formulate the final discretized scheme in terms of surface-averaged variables \citep{11Azmy_1983,azmy1985nodal,12Raj_2017}, while others are based on node-averaged variables \citep{08Raj_2022,09Shober_1977,13Shober1978}. In the nuclear industry, both methods have proven to be highly accurate and efficient, sharing a common primary step known as the transverse integration procedure (TIP). This procedure simplifies the PDEs at each node by transforming them into a system of ordinary differential equations (ODEs). The analytical solutions of these transversely averaged ODEs are then computed at each node. Subsequently, these analytical solutions are subjected to physically relevant constraint conditions, which results in the final expression of the numerical scheme. The final numerical scheme for the surface-averaged (and node-averaged) based nodal method includes surface-averaged  (and node-averaged) variables as discrete unknowns. Both surface-averaged and node-averaged schemes have been successfully applied in the nuclear industry and are effective at solving PDEs efficiently. 
Surface-averaged nodal methods, commonly referred to as the ``nodal integral method'' (NIM) or ``analytical nodal methods'' (ANM), have been successfully applied to a wide range of fluid flow problems. These include the convection-diffusion equation (CDE) \citep{22Rizwan_uddin_1997,21Michael_2001,nezami2009nodal,neeraj2013convection,zhou2016general,jarrah2021nodal,jarrah2022nodal,Neeraj2024coupled}, Burgers' equations \citep{23Rizwan_uddin_1997,25Wescott_2001,kumar2016numerical,kumar2019physics,Niteen2020predictor,Niteen2021novel,namala2019hybrid,ahmed2022modified,ahmed2022picard}, and Navier-Stokes equations (NSE) \citep{11Azmy_1983,16Wang_2003,wang2005modified,neeraj2013navier,zhou2016general,jarrah2024nodal,ahmed2024physics,14Ahmed_2024,singh2009k,kumar2012pressure}, among other scientific and engineering applications \citep{Niteen2022nodal,gander2022new}. However, nodal schemes based on node-averaged variables have been less popular in the fluid flow community, primarily due to the complexities involved in their development process \citep{09Shober_1977,13Shober1978}. 

Recently, advancements in nodal schemes utilizing node-averaged variables for fluid flow problems have emerged, collectively termed as the “cell-centered nodal integral method” (CCNIM) \citep{19Ahmed_2021,18Ahmed_2023,17Ahmed_2024,Ahmed_Diffusion_2024}. This terminology aligns with the conventions in the fluid flow community, where “node-averaged” is often referred to as “cell-centered” (CC) \citep{18Ahmed_2023}.
For instance, \citet{09Shober_1977} is credited with being the first to develop a variant of the nodal method that focuses on cell-centered variables, specifically designed to solve the neutron diffusion equation. Building on this foundation, Ahmed et al. \citep{19Ahmed_2021,18Ahmed_2023,17Ahmed_2024,Ahmed_Diffusion_2024} were the first to adapt the cell-centered nodal methodology to fluid flow problems, introducing a simplified variant specifically designed to solve diffusion \citep{19Ahmed_2021,Ahmed_Diffusion_2024} and convection-diffusion problems \citep{18Ahmed_2023,17Ahmed_2024}. They  introduced the term cell-centered NIM (CCNIM) and provided a detailed discussion of this methodology for addressing fluid flow problems \citep{19Ahmed_2021,18Ahmed_2023,17Ahmed_2024,Ahmed_Diffusion_2024}. The CCNIM builds upon the same foundational principle as the TIP, a common step in various nodal approaches. 
 However, CCNIM \citep{18Ahmed_2023,17Ahmed_2024} differs significantly from the standard procedure of traditional NIM \citep{22Rizwan_uddin_1997,23Rizwan_uddin_1997}, which typically employs surface-averaged variables at the node boundaries as local boundary conditions for solving the ODEs at each node. Instead, CCNIM utilizes `interface-averaged' flux and variables from neighboring nodes as local boundary conditions for these ODEs, rather than relying on just interface-averaged variables. This modification leads to a formulation in which the discrete unknowns at each node are cell-centered values, in contrast to traditional NIMs, which use surface-averaged variables per node \citep{18Ahmed_2023,17Ahmed_2024}. The temporal derivative is treated explicitly in the earlier development of CCNIM \citep{18Ahmed_2023}, which offers certain benefits over conventional NIM. These include a straightforward formulation for higher-degree temporal derivatives, improved flexibility with respect to the higher level of accuracy in time, a simplified implementation procedure in terms of local coordinate systems, and a promising alternative for coupling with other physics. A detailed explanation of these advantages is provided in the literature \citep{18Ahmed_2023,17Ahmed_2024}. Due to the explicit treatment of the temporal derivative, the final system produced by CCNIM links algebraic equations with ODEs, creating a complicated system known as a system of differential-algebraic equations (DAEs) \cite{18Ahmed_2023}. DAEs present considerable challenges in terms of complexity, largely stemming from the intricate coupling between their algebraic and differential components, along with other contributing factors. Furthermore, CCNIM has been hindered by challenges like managing Neumann boundary conditions and not being applicable to one-dimensional problems \cite{17Ahmed_2024}. 
 
 Consequently, a modified CCNIM (known as MCCNIM) is developed in order to address the issues of CCNIM \citep{17Ahmed_2024}. MCCNIM discretizes both space and time in the nodal framework, ensuring second-order accuracy in both dimensions. This approach eliminates complex DAE systems, significantly simplifying the computational process by generating an algebraic set of equations for the discrete unknowns at each node. Furthermore, the use of a straightforward flux definition facilitates the seamless incorporation of Neumann boundary conditions. MCCNIM is also adaptable to one-dimensional problems. Notably, by treating the temporal derivative explicitly, the same MCCNIM framework, along with its flux definition, can be readily converted into a CCNIM method. Until now, both schemes, CCNIM and MCCNIM, have been exclusively developed for solving linear fluid flow problems, such as heat conduction and linear convection-diffusion equations. 

In this work, the modified CCNIM (MCCNIM) has been extended to solve nonlinear PDEs, with a particular focus on Burgers' equations. A novel approach has been introduced for handling the convective velocity in the nonlinear convection term, ensuring that MCCNIM retains its superior accuracy compared to traditional NIMs when applied to nonlinear problems. It is worth noting that we chose to develop MCCNIM because it can be applied to both one- and two-dimensional problems, whereas CCNIM is limited to two-dimensional problems. However, if necessary, the CCNIM approach could also be developed by taking the temporal derivative outside and following a development procedure similar to that outlined in the present work for MCCNIM. The proposed scheme has been validated by solving various one- and two-dimensional Burgers' problems, the analytical solutions of which are available. Additionally, a comprehensive comparison of the proposed scheme with previously published results is presented to demonstrate its effectiveness. The rate of convergence is rigorously tested for the chosen problems, confirming that the scheme achieves second-order accuracy in both space and time.

The paper is organized as follows: \sect\ref{sec:2} develops the MCCNIM scheme for one- and two-dimensional time-dependent Burgers' equations. Numerical results are presented in \sect\ref{sec:3} to demonstrate the performance of the scheme, and, finally, conclusions and discussions are provided in \sect\ref{sec:4}.
%
\section{MCCNIM Formulation}
\label{sec:2}
The essential steps of the MCCNIM are briefly outlined here, with a detailed discussion of their application to one- and two-dimensional Burgers' equations provided in the following sections. The issues related to handling nonlinear convective velocity in both one- and two-dimensional cases are addressed in the concluding sections. 
In general, there are six significant steps in the MCCNIM implementation procedure.
\begin{enumerate}
	\item 	Finite-sized rectangular elements, known as nodes or cells, are used to divide the space-time domain. The width of each cell is equal to $2a_i$, denoted as $\Delta x$. Similarly, the height of each cell is equal to $2\tau_{\ell}$, denoted as $\Delta t$. This may be seen in panel (b) of \fig\ref{fig:01}. After discretization, the ``transverse integration process'' (TIP) is applied to each cell, transforming the partial differential equations (PDEs) into a group of ordinary differential equations (ODEs). 
	\item The set of ODEs is separated into homogeneous and inhomogeneous parts, where the former is easily integrable, whereas the latter poses integration challenges. The homogeneous part is placed on the left side of the differential equation, and the inhomogeneous part on the right side is denoted as the pseudo-source term. These pseudo-source terms are expanded to a specific order using Legendre polynomials. Subsequently, within each cell, the homogeneous components are solved analytically with the interface-averaged flux and variables acting as local boundary conditions.
	\item The analytical solutions of these ODEs are then integrated (averaged) along the remaining direction to derive the expression of the cell-centered (node-averaged) variables.
	\item Ensuring that the continuity criteria are satisfied at the shared boundary of neighboring nodes yields equations for surface-averaged variables.
	\item Following that, physically relevant constraint conditions are established in order to close the system of equations. The initial two constraint conditions are obtained through systematic manipulation of the definitions of pseudo-source terms from Step 2 and the utilization of surface-averaged variable expressions from Step 4. This process establishes linkages between surface-averaged variables and pseudo-source terms. A supplementary constraint condition is derived by calculating the average of the original PDEs over each dimension.
	\item At last, the final set of algebraic equations in terms of cell-centered (node-averaged) variables is derived by incorporating the expressions of the surface-averaged variables obtained in step 4 into the constraint conditions.
\end{enumerate}
The steps discussed for MCCNIM are further elaborated by applying them to the one- and two-dimensional time-dependent Burgers’ equations in the following sections. 
The one-dimensional transient Burgers' equation has been used to provide a detailed formulation of our proposed approach for addressing a nonlinear problem. Following a similar procedure, a brief derivation of the two-dimensional Burgers’ equation is also presented in the subsequent section.
\subsection{One-dimensional Burgers’ equation}
\label{sec:2.1}
\begin{figure}[!b]%
	\centering
	\subfigure[~]{\includegraphics[width=0.48\linewidth]{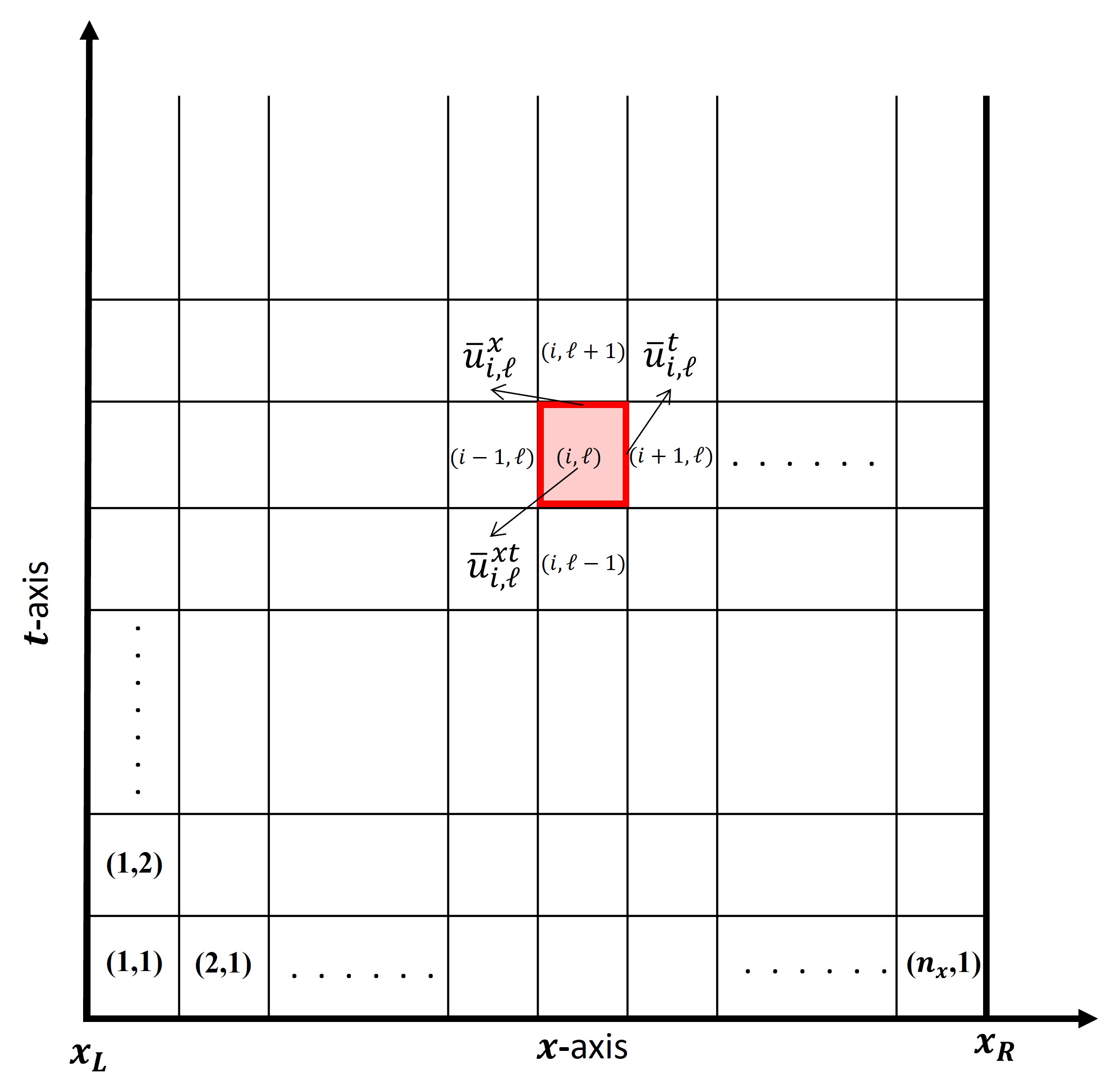}\label{fig1a}}
	\subfigure[~] {\includegraphics[width=0.48\linewidth]{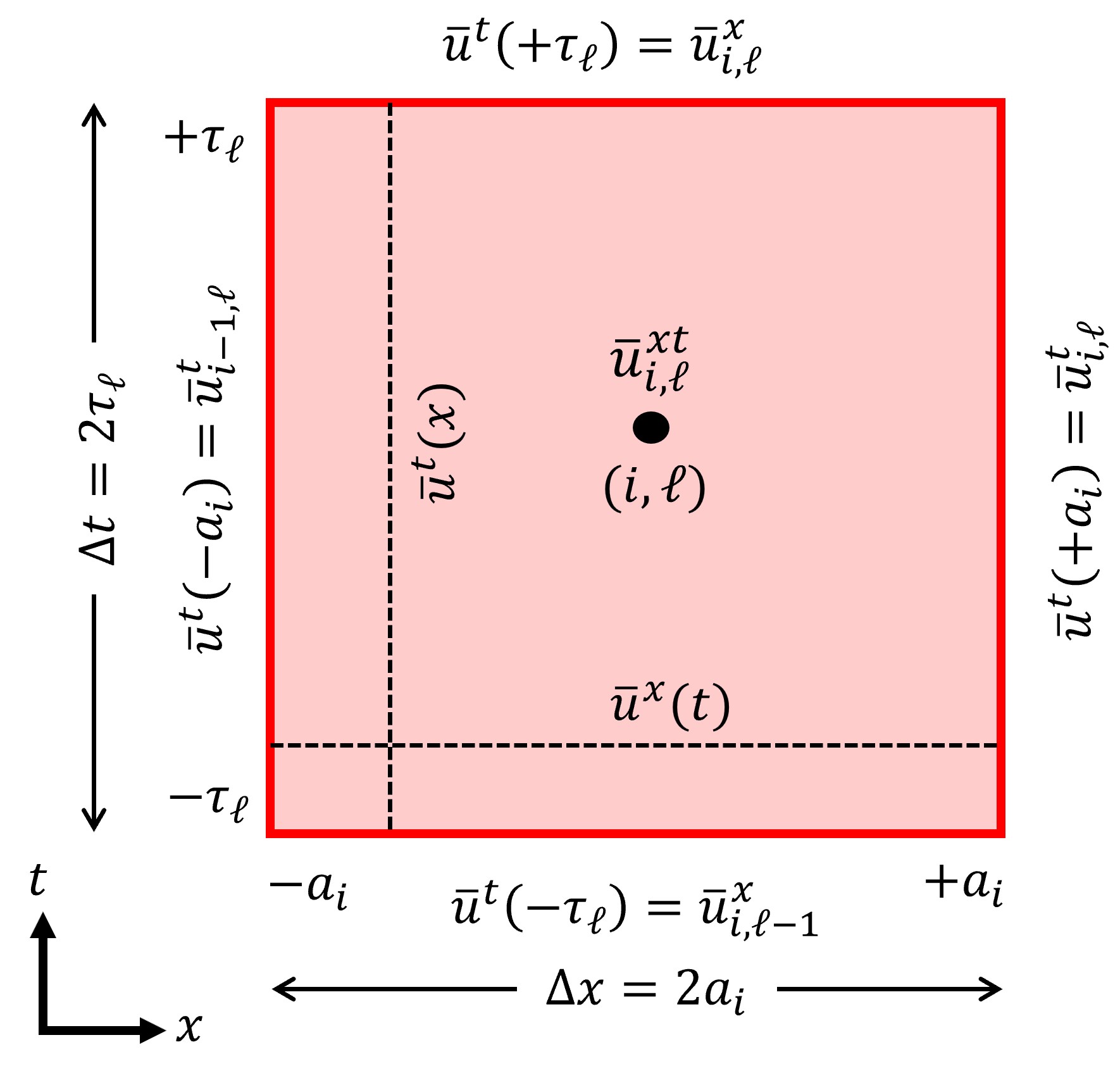}\label{fig1b}}%
	\caption{(a) A schematic representation of the global ($x,t$) space, and its division into the rectangular nodes or elements for one-dimensional transient problems in MCCNIM. (b) Local coordinate system and transverse-averaged quantities within the node ($i,\ell$).}%
	\label{fig:01}%
\end{figure}

The one-dimensional time-dependent Burgers’ equation is given as
\begin{gather}
\frac{\partial u(x,t)}{\partial t}+u\left(x,t\right)\frac{\partial u\left(x,t\right)}{\partial x}=\frac{1}{Re}\frac{\partial^2u\left(x,t\right)}{\partial x^2} 
\label{eq:001}
\end{gather}
where $u(x,t)$ is the velocity, and $Re$ is the Reynolds number. As shown in panel (a) of \fig\ref{fig:01}, the space-time domain for the one-dimensional problem described by \eqn\eqref{eq:001} is discretized into rectangular nodes indexed by $i$ and ${\ell}$, where $i$ represents the spatial index and ${\ell}$ the temporal index. Each node ($i,\ell$) is defined by dimensions ($\Delta x\times\Delta t$), with the width of each cell equal to $2a_i$ (i.e., $\Delta x=2a_i$) and the height equal to $2\tau_{\ell}$ (i.e., $\Delta t=2\tau_{\ell}$). The origin of the node is positioned at the node center, as illustrated in panel (b) of \fig\ref{fig:01}.

Prior to the development of the numerical scheme, the one-dimensional Burgers’ equation  (\eqn\ref{eq:001}) is reformulated in terms of the local coordinate system, as follows.
\begin{gather}
\frac{\partial u(x,t)}{\partial t}=\frac{1}{Re}\frac{\partial^2u\left(x,t\right)}{\partial x^2} - {\bar{u}}^0\frac{\partial u\left(x,t\right)}{\partial x}
\label{eq:002}
\end{gather}
where ${\bar{u}}^0$ is approximated (i.e., the node-averaged) convective velocity at the present time step. The definition of ${\bar{u}}^0$ is discussed in \sect\ref{sec:2.1.7}.
\subsubsection{Transverse integration process}
\label{sec:2.1.1}
This step, known as the transverse integration procedure (TIP), is common to all nodal schemes \citep{22Rizwan_uddin_1997} in which, after discretizing the space-time domain, the PDE is integrated transversely (averaged) in each direction. As we focus on developing the scheme for the one-dimensional Burgers' equation (\eqn\ref{eq:001}), which involves both $x-$ and $t-$directions, the integration is performed sequentially by using the spatial, temporal and double integration operators defined as follows.
\begin{gather*}
	I_{x}(g)= \frac{1}{2g}\int_{-g}^{+g} F \mathrm{d}x; \quad 
	I_{t}(h) = \frac{1}{2h}\int_{-h}^{+h} F \mathrm{d}t;  \quad	
	I_{xt}(g,h) = \frac{1}{4gh}\int_{-g}^{+g}\int_{-h}^{+h} F \mathrm{d}t\mathrm{d}x
\end{gather*}  
where, $F$ is the integrand, and $g~(=a_i \text{ or } a_{i+1})$ and $h~ (=\tau_{\ell})$ refers to the limits of integration.

First, the integration is performed in the $x$-direction using the spatial integration operator $I_x(a_i)$ on \eqn\eqref{eq:001}, followed by in the $t$-direction using the temporal integration operator $I_t(\tau_{\ell})$ on \eqn\eqref{eq:001} over each node ($i,\ell$) of the discretized domain. 

Employing the spatial integration operator $I_x(a_i)$ on \eqn\eqref{eq:002} yields
\begin{gather}
\frac{1}{2a_i}\int_{-a_i}^{+a_i}\left(\frac{\partial u(x,t)}{\partial t}\right)\mathrm{d}x =\frac{1}{2a_i}\int_{-a_i}^{+a_i}\left(\frac{1}{Re}\frac{\partial^2u\left(x,t\right)}{\partial x^2}-{\bar{u}}^0\frac{\partial u\left(x,t\right)}{\partial x}\right)\mathrm{d}x 
\label{eq:003}
\end{gather}
Note that \eqn\eqref{eq:003} is arranged so that the terms straightforward integrable are placed on the left-hand side (LHS), while the remaining terms are shifted to the right-hand side (RHS) to serve as pseudo-source terms. Thus, from \eqn\eqref{eq:003}, the first ODE which is averaged in space and dependent on time for node ($i,\ell$) is then written as
\begin{gather}
\frac{\mathrm{d}{\bar{u}}^x(t)}{\mathrm{d}t}={\bar{S}}_1^x(t) 
\label{eq:004}
\end{gather}
Likewise, employing the temporal integration operator $I_t(\tau_{\ell})$ on \eqn\eqref{eq:002}, and adhering to the analogous steps in achieving \eqn\eqref{eq:004}, an additional time-averaged and space-dependent ODE can be obtained as
\begin{gather}
\frac{1}{Re}\frac{\mathrm{d}^2{\bar{u}}^t(x)}{\mathrm{d}x^2}-{\bar{u}}^0\frac{\mathrm{d}{\bar{u}}^t\left(x\right)}{\mathrm{d}x}={\bar{S}}_2^t\left(x\right)  
\label{eq:005}
\end{gather}
where\add{, ${\bar{S}}_1^x(t)$ and ${\bar{S}}_2^t(x)$ are the space-averaged (time-dependent) and time-averaged (spatial-dependent) pseudo-source terms, respectively. The} transverse-averaged physical quantities (${\bar{u}}^t$ and ${\bar{u}}^x$) and pseudo-source terms (${\bar{S}}_1^x$ and ${\bar{S}}_2^t$) \add{appearing in \eqns\eqref{eq:004} and \eqref{eq:005}} are defined as follows.
\begin{gather}
	{\bar{u}}^t\left(x\right)=\frac{1}{2\tau_{\ell}}\int_{-\tau_{\ell}}^{+\tau_{\ell}}{u(x,t)\mathrm{d}t} 
\label{eq:006}
\end{gather}
\begin{gather}
	{\bar{u}}^x\left(t\right)=\frac{1}{2a_i}\int_{-a_i}^{+a_i}{u(x,t)\mathrm{d}x} 
\label{eq:007}
\end{gather}
\begin{gather}
	{\bar{S}}_1^x\left(t\right)=\frac{1}{2a_i}\int_{-a_i}^{+a_i}\left(\frac{1}{Re}\frac{\partial^2u\left(x,t\right)}{\partial x^2}-{\bar{u}}^0\frac{\partial u\left(x,t\right)}{\partial x}\right)\mathrm{d}x 
\label{eq:008}
\end{gather}
\begin{gather}
{\bar{S}}_2^t\left(x\right)=\frac{1}{2\tau_{\ell}}\int_{-\tau_{\ell}}^{+\tau_{\ell}}{\frac{\partial u(x,t)}{\partial t}\mathrm{d}t} 
\label{eq:009}
\end{gather}
The product of the averages has been used to estimate the average of the product, representing a second-order approximation \citep{21Michael_2001,22Rizwan_uddin_1997,azmy1985nodal}. The aforementioned progression holds true for all cells. 

One crucial step in solving locally transverse-integrated ODEs (\eqns\ref{eq:004} and \ref{eq:005}) is the use of Legendre polynomials to expand the pseudo-source terms (${\bar{S}}_1^x$ and ${\bar{S}}_2^t$). The expansion of these pseudo-source terms is truncated to the zeroth order during this operation, resulting in a constant. 
Prior studies \citep{11Azmy_1983,20Elnawawy_1990,21Michael_2001} have established that by truncating the expansion at the zeroth order level, one can obtain a scheme that is accurate to the second order. Nevertheless, it is important to acknowledge, that higher-order schemes are possible as a result of choosing to truncate the expansion at higher orders \citep{21Michael_2001,guessous2002higher,Guessous2003high,guessous2006order,guessous2016three}. 

\add{The Legendre polynomial expansions of the pseudo-source terms, i.e., ${\bar{S}}_1^x$ in time and ${\bar{S}}_2^t$ in space, have been truncated to zeroth order, as discussed above, resulting in the cell-centered quantities denoted as ${\bar{S}}_{1i,\ell}^{xt}$ and ${\bar{S}}_{2i,\ell}^{xt}$, respectively.  Both terms are assigned the same superscript `${xt}$' as they have now been integrated/expanded in both $x$- and $t$-directions at each node. However, despite sharing the same superscript, the two terms are fundamentally different (\eqns\ref{eq:008} and \ref{eq:009}) and distinguished by the subscripts `1' and `2' along with the nodal indices $(i,\ell)$, aligning with the main variable $\bar{u}_{i,\ell}^{xt}$ in the final set of algebraic equations.}
Following truncation, \eqns\eqref{eq:004} and \eqref{eq:005} simplify into two ODEs as follows.
\begin{gather}
\frac{\mathrm{d}{\bar{u}}^x(t)}{\mathrm{d}t}={\bar{S}}_{1i,\ell}^{xt} 
\label{eq:010}
\end{gather}
\begin{gather}
	\frac{1}{Re}\frac{\mathrm{d}^2{\bar{u}}^t(x)}{\mathrm{d}x^2}-{\bar{u}}_{i,\ell}^0\frac{\mathrm{d}{\bar{u}}^t(x)}{\mathrm{d}x}={\bar{S}}_{2i,\ell}^{xt} 
\label{eq:011}
\end{gather}
\subsubsection{Analytical solution of ODEs}
\label{sec:2.1.2}
The analytical solutions to the ODEs derived in the previous subsection is then implemented for individual nodes. While the first ODE (\eqn\ref{eq:010}) is solved analytically in a manner similar to the traditional NIM, the second ODE (\eqn\ref{eq:011}) is solved using a different approach. This approach involves applying modified cell boundary conditions, in contrast to conventional NIM \citep{22Rizwan_uddin_1997,23Rizwan_uddin_1997}. Specifically, this method uses the average velocity and average flux of the interface as local boundary conditions, as illustrated in \fig\ref{fig:02}.

By solving \eqn\eqref{eq:010} (first ODE) analytically and applying the nodal boundary condition, ${\bar{u}}^x(t)={\bar{u}}^x\left({+\tau}_{\ell}\right)$ at the surface $t={+\tau}_{\ell}$, as depicted in the panel (a) of \fig\ref{fig:02}, we obtain: 
\begin{figure}[!b]%
	\centering
	\subfigure[~]{\includegraphics[width=0.48\linewidth]{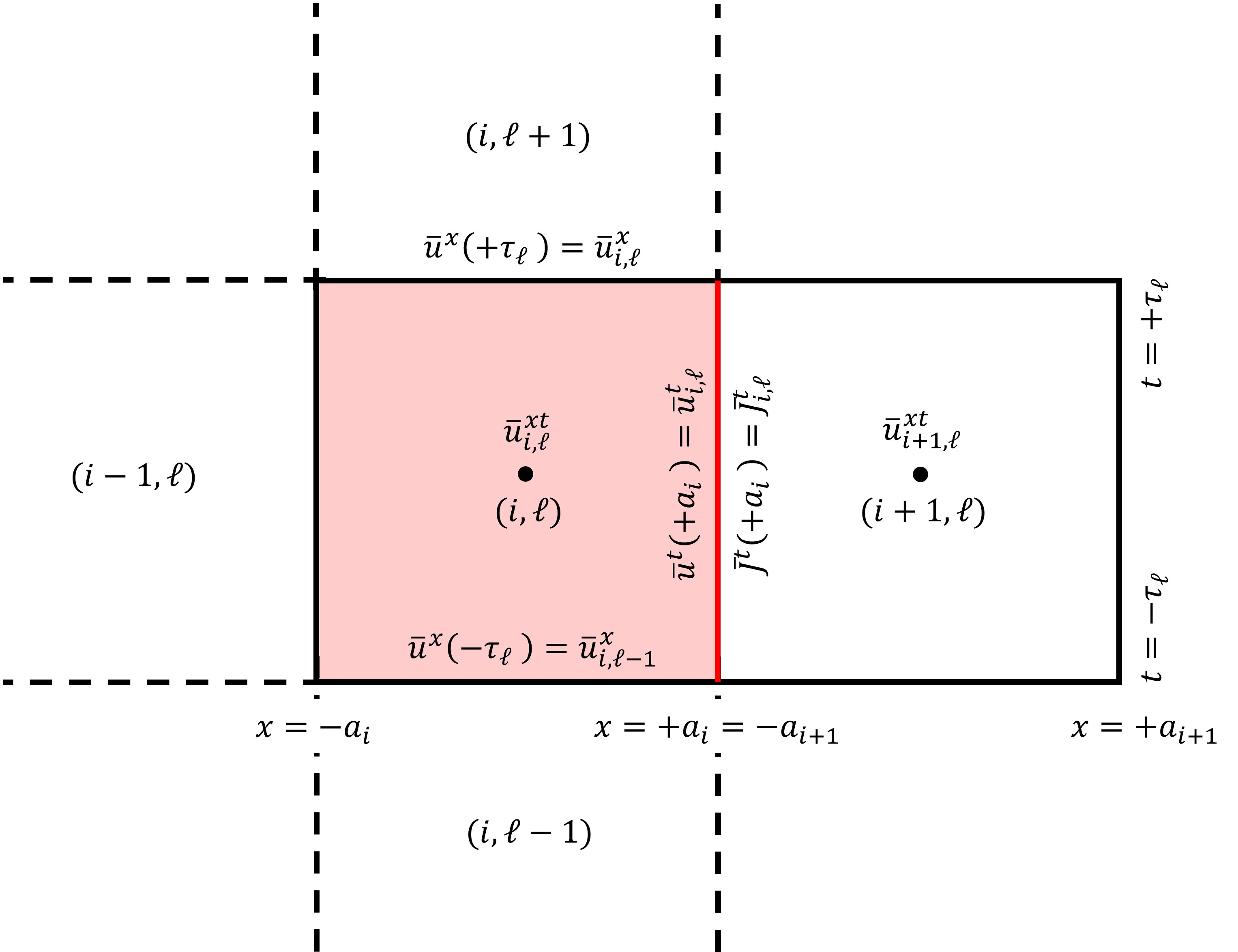}\label{fig2a}}
	\subfigure[~] {\includegraphics[width=0.48\linewidth]{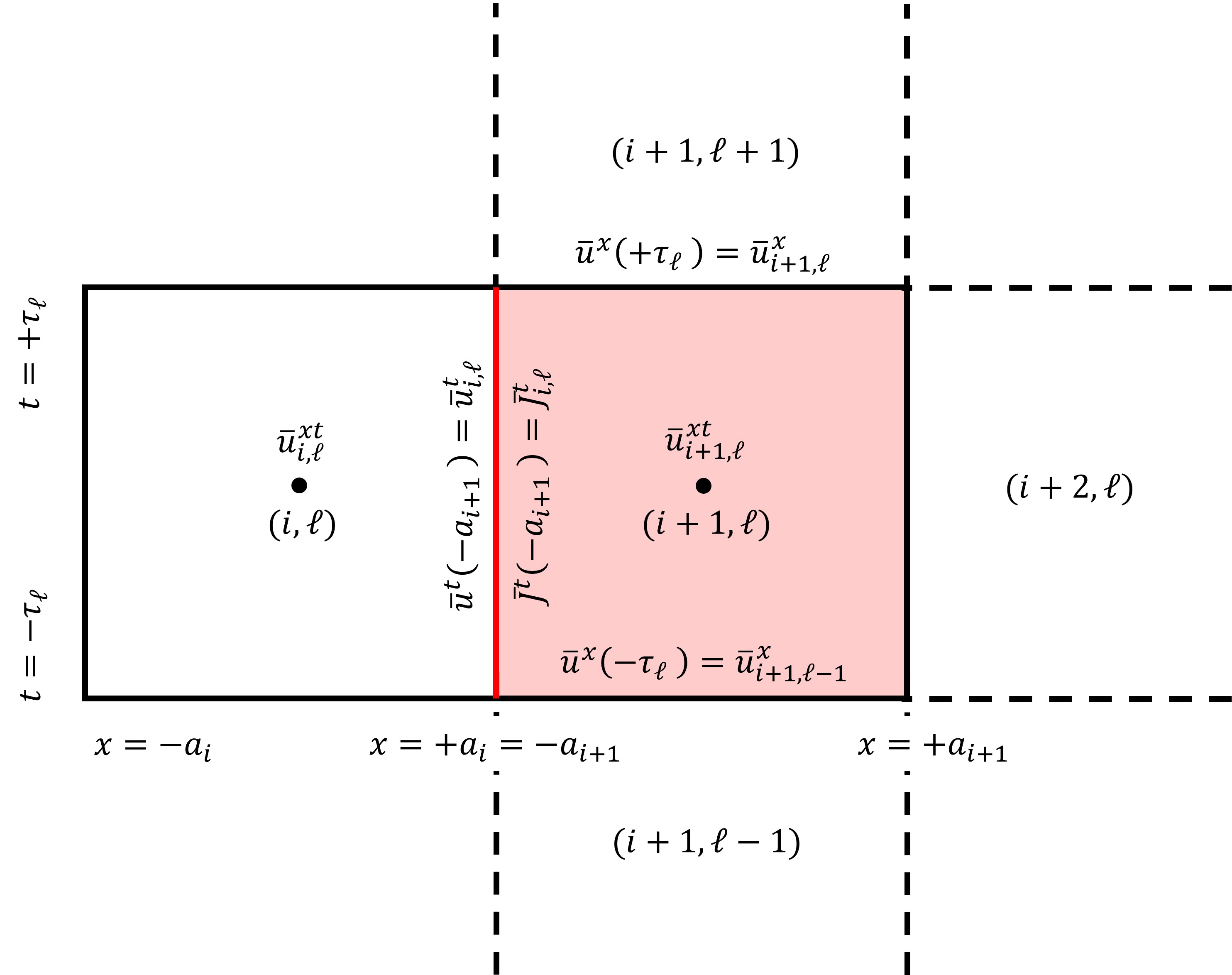}\label{fig2b}}%
	\caption{Local boundary conditions for MCCNIM. (a) For node ($i,\ell$). (b) For node ($i+1,\ell$)}%
	\label{fig:02}%
\end{figure}
\begin{gather}
	{\bar{u}}^x\left(t\right)={\bar{S}}_{1i,\ell}^{xt}\left(t-\tau_{\ell}\right)+{\bar{u}}^x\left({+\tau}_{\ell}\right)
\label{eq:012}
\end{gather}
\eqn\eqref{eq:012} can be rewritten in ($i,\ell$) indexing using ${\bar{u}}^x\left({+\tau}_{\ell}\right)={\bar{u}}_{i,\ell}^x$ as,
\begin{gather}
	{\bar{u}}^x\left(t\right)={\bar{S}}_{1i,\ell}^{xt}\left(t-\tau_{\ell}\right)+{\bar{u}}_{i,\ell}^x
\label{eq:013}
\end{gather}
For the solution of second ODE (\eqn\ref{eq:011}), current approach employs the identical analytical solution procedure as the prior CCNIM, however, with a simpler definition of interface flux for the boundary conditions. The convective term incorporated into the definition of flux in the previous CCNIM \citep{18Ahmed_2023} poses challenge when attempting to apply Neumann or insulated boundary conditions with the nodal spirit.  Consequently, the insulated boundary conditions were applied using a second-order finite-difference approximation, which compromises the fundamental nature of the nodal scheme. To address this issue and enhance the adaptability and adherence to physical principles, a straightforward definition of flux was introduced in previous work \citep{17Ahmed_2024}, as illustrated in \fig\ref{fig:02}. This revised formulation effectively distinguished between the convective and diffusive aspects of the flux and facilitated the independent application of distinct boundary conditions to each componentry, thereby enhancing the practicality and accuracy of the method. The modified flux at the common edge of the neighbouring nodes is defined as
\begin{gather}
{\bar{J}}^t\left(x\right)=-\frac{1}{2\tau_{\ell}Re}\int_{-\tau_{\ell}}^{+\tau_{\ell}}{\frac{\partial u(x,t)}{\partial x}dt}=-\frac{1}{Re}\frac{\mathrm{d}{\bar{u}}^t(x)}{\mathrm{d}x} 
\label{eq:014}
\end{gather}
By using the surface-averaged flux, defined in \eqn\eqref{eq:014}, and the surface-averaged velocity as local boundary conditions at the interface between nodes ($i,\ell$) and ($i+1,\ell$), as shown in panel (a) of \fig\ref{fig:02}, we obtain the solution to the second ODE (\eqn\ref{eq:011}). 

For node ($i,\ell$), the analytical solution of \eqn\eqref{eq:011} is obtained by applying the average velocity ${\bar{u}}^t\left(x\right)={\bar{u}}^t\left(+a_i\right)$ and averaged flux ${\bar{J}}^t\left(x\right)={\bar{J}}^t\left(+a_i\right)$ at the interface $x={+a}_i$, as depicted in panel (a) of \fig\ref{fig:02}, leading to:
\begin{gather}
	\begin{split}
\left[{\bar{u}}^t\left(x\right)\right]_{i,\ell}={\bar{u}}^t\left({+a}_i\right)-\frac{-1+e^{Re{\bar{u}}_{i,\ell}^0\left(x-a_i\right)}}{\ {\bar{u}}_{i,\ell}^0}{\bar{J}}^t\left(+a_i\right) 
+\frac{\left(-1+e^{Re{\bar{u}}_{i,\ell}^0\left(x-a_i\right)}-Re{\bar{u}}_{i,\ell}^0(x-a_i)\right)}{Re\left({\bar{u}}_{i,\ell}^0\right)^2}{\bar{S}}_{2i,\ell}^{xt} 
	\end{split}
\label{eq:015}
\end{gather}
Similarly, the solution at node ($i+1,\ell$) can be derived  at $x={-a}_{i+1}$ using boundary conditions, ${\bar{u}}^t\left(x\right)={\bar{u}}^t\left(-a_{i+1}\right)$ and ${\bar{J}}^t\left(x\right)={\bar{J}}^t\left(-a_{i+1}\right)$, as depicted in panel (b) of \fig\ref{fig:02}, given as follows:
\begin{gather}
	\begin{split}
\left[{\bar{u}}^t\left(x\right)\right]_{i+1,\ell}={\bar{u}}^t\left({-a}_{i+1}\right)-\frac{-1+e^{Re{\bar{u}}_{i+1,\ell}^0\left(x+a_{i+1}\right)}}{\ {\bar{u}}_{i+1,\ell}^0}{\bar{J}}^t\left(-a_{i+1}\right)+\frac{\left(-1+e^{Re{\bar{u}}_{i+1,\ell}^0\left(x+a_{i+1}\right)}-Re{\bar{u}}_{i+1,\ell}^0(x+a_{i+1})\right)}{Re\left({\bar{u}}_{i+1,\ell}^0\right)^2}{\bar{S}}_{2i+1,\ell}^{xt} 
	\end{split}
\label{eq:016}
\end{gather}
Up to this point, the procedure has strictly adhered to conventional NIM practices. However, a key distinction arises in the boundary conditions for solving \eqn\eqref{eq:011}. In conventional NIM, surface-averaged variables at both cell edges are used. In contrast, MCCNIM employs both surface-averaged variables and fluxes at the interface between adjacent nodes, as shown in \fig\ref{fig:02}. A detailed comparison between these approaches is provided in the  CCNIM literature \citep{17Ahmed_2024,18Ahmed_2023,19Ahmed_2021}.
\subsubsection{Evaluation of the cell-centered (CC) variables}
\label{sec:2.1.3}
The solutions derived in the previous subsection for ${\bar{u}}^x\left(t\right)$ given by \eqn\eqref{eq:013} and ${\bar{u}}^t\left(x\right)$ given by \eqns\eqref{eq:015} and \eqref{eq:016} are now subjected to a second averaging (integration) process to determine the cell-centered variables. 
Initially, by applying the temporal integration operator $I_t(\tau_{\ell})$ on \eqn\eqref{eq:013}, we obtain:
\begin{gather}
	{\bar{u}}_{i,\ell}^{xt}=-\tau_{\ell}{\bar{S}}_{1i,\ell}^{xt}+{\bar{u}}_{i,\ell}^x 
\label{eq:017}
\end{gather}
where ${\bar{u}}_{i,\ell}^{xt}$ represents the cell-centered (node-averaged) velocity. Note that the pseudo-source term ${\bar{S}}_{1i,\ell}^{xt}$ is inherently a cell-centered value following truncation, as discussed above. 
Now, by applying the spatial integration operator $I_x(a_i)$ on \eqn\eqref{eq:015} and  $I_x(a_{i+1})$ on \eqn\eqref{eq:016}, we obtain:
\begin{gather}
	\begin{split}
	{\bar{u}}_{i,\ell}^{xt}&={\bar{u}}^t\left({+a}_i\right)-\frac{\left(1-e^{-{Reu}_{i,\ell}}-2a_iRe{\bar{u}}_{i,\ell}^0\right)}{2a_i{Re}\left({\bar{u}}_{i,\ell}^0\right)^2}{\bar{J}}^t\left(+a_i\right)\\&+\frac{\left(1-e^{-{Reu}_{i,\ell}}+2a_iRe{\bar{u}}_{i,\ell}^0\left(-1+a_iRe{\bar{u}}_{i,\ell}^0\right)\right)}{2a_i{Re}^2\left({\bar{u}}_{i,\ell}^0\right)^3}{\bar{S}}_{2i,\ell}^{xt} 
	\end{split}
\label{eq:018}
\end{gather}
\begin{gather}
	\begin{split}
		{\bar{u}}_{i+1,\ell}^{xt}&={\bar{u}}^t\left({-a}_{i+1}\right)+\frac{\left(1-e^{{Reu}_{i+1,\ell}}-2a_{i+1}Re{\bar{u}}_{i+1,\ell}^0\right)}{2a_{i+1}{Re}\left({\bar{u}}_{i+1,\ell}^{0}\right)^2}{\bar{J}}^t\left(-a_{i+1}\right)\\ &+\frac{\left(-1+e^{{Reu}_{i+1,\ell}}-2a_{i+1}Re{\bar{u}}_{i+1,\ell}^0\left(1+a_{i+1}Re{\bar{u}}_{i+1,\ell}^0\right)\right)}{2a_{i+1}{Re}^2\left({\bar{u}}_{i+1,\ell}^0\right)^3}{\bar{S}}_{2i+1,\ell}^{xt} 
	\end{split}
\label{eq:019}
\end{gather}
In \eqns\eqref{eq:018} and \eqref{eq:019}, each power of the exponential term is written in terms of the local Reynolds number. For instance, at node ($i,\ell$) the local Reynolds number is defined as ${Reu}_{i,\ell}=2a_iRe{\bar{u}}_{i,\ell}^{0}$. Furthermore, \eqns\eqref{eq:018} and \eqref{eq:019} can be rearranged to yield,
\begin{gather}
	{\bar{J}}^t\left(+a_i\right)=A_{31}\left({\bar{u}}_{i,\ell}^{xt}-{\bar{u}}^t\left({+a}_i\right)\right)+A_{32}{\bar{S}}_{2i,\ell}^{xt}  
\label{eq:020}
\end{gather}
\begin{gather}
{\bar{J}}^t\left(-a_{i+1}\right)=A_{51,i+1}\left({\bar{u}}_{i+1,\ell}^{xt}-{\bar{u}}^t\left({-a}_{i+1}\right)\right)+A_{52,i+1}{\bar{S}}_{2i+1,\ell}^{xt} 
\label{eq:021}
\end{gather}
where $A_{31}$, $A_{32}$, $A_{51}$ and $A_{52}$ are the coefficients that are dependent on the various parameters ($a_i$, ${Reu}_{i,\ell}$, ${\bar{u}}_{i,\ell}^{0}$). The explicit definitions of these coefficients are presented in \ref{app:A}. 
\subsubsection{Continuity conditions}
\label{sec:2.1.4}
Since the surface-averaged velocity is solely needed as a boundary condition for the solution of the first ODE (\eqn\ref{eq:010}), it is evident that no continuity condition is required in the temporal direction. As a result, \eqn\eqref{eq:017} can be reformulated with all cell-centered values on RHS, which is presented as follows.
\begin{gather}
	{\bar{u}}_{i,\ell}^x={\bar{u}}_{i,\ell}^{xt}+\tau_{\ell}{\bar{S}}_{1i,\ell}^{xt} 
\label{eq:022}
\end{gather}
\eqn\eqref{eq:022} is applicable to node ($i,\ell$), a corresponding equation for node ($i,\ell-1$) can be obtained through a simple index shift, given as
\begin{gather}
	{\bar{u}}_{i,\ell-1}^x={\bar{u}}_{i,\ell-1}^{xt}+\tau_{\ell-1}{\bar{S}}_{1i,\ell-1}^{xt} 
\label{eq:023}
\end{gather}
Spatially, two continuity conditions must be imposed: (a) continuity of the surface-averaged velocity at the common edge between adjacent nodes, and (b) continuity of the flux at the shared edge between the same neighbouring nodes. The continuity of the surface-averaged velocity between nodes ($i,\ell$) and ($i+1,\ell$) can be written as
\begin{gather}
	{\bar{u}}^t\left({+a}_i\right)={\bar{u}}^t\left({-a}_{i+1}\right)={\bar{u}}_{i,\ell}^t 
\label{eq:024}
\end{gather}
Similarly, the continuity of the surface averaged flux between nodes ($i,\ell$) and ($i+1,\ell$) can be expressed as
\begin{gather}
{\bar{J}}^t\left(+a_i\right)={\bar{J}}^t\left(-a_{i+1}\right)={\bar{J}}_{i,\ell}^t 
\label{eq:025}
\end{gather}
On applying the continuity conditions, utilizing \eqns\eqref{eq:020} and \eqref{eq:021}, and solving for the surface-averaged velocity, ${\bar{u}}_{i,\ell}^t$, we get,
\begin{gather}
	{\bar{u}}_{i,\ell}^t=\frac{A_{32}{\bar{S}}_{2i,\ell}^{xt}-\ A_{52,i+1}{\bar{S}}_{2i+1,\ell}^{xt}+\ A_{31}{\bar{u}}_{i,\ell}^{xt}-\ A_{51,i+1}{\bar{u}}_{i+1,\ell}^{xt}}{A_{31}-A_{51,i+1}} 
\label{eq:026}
\end{gather}
Now, substituting the expression for ${\bar{u}}_{i,\ell}^t$ from \eqn\eqref{eq:026} into either flux equation (\eqns\ref{eq:020} and \ref{eq:021}) yields the same expression for ${\bar{J}}_{i,\ell}^t$, conforming that the continuity conditions have been met. The resulting expression of ${\bar{J}}_{i,\ell}^t$ is as follows:
\begin{gather}
	{\bar{J}}_{i,\ell}^t=\frac{-A_{32}A_{51,i+1}{\bar{S}}_{2i,\ell}^{xt}+\ {A_{31}A}_{52,i+1}{\bar{S}}_{2i+1,\ell}^{xt}+\ A_{31}A_{51,i+1}({\bar{u}}_{i+1,\ell}^{xt}-{\bar{u}}_{i,\ell}^{xt})}{A_{31}-A_{51,i+1}}  
\label{eq:027}
\end{gather}
Similarly, by changing the index from ($i,\ell$) to ($i-1,\ell$), we can readily write analogous expressions for ${\bar{u}}_{i-1,\ell}^t$ and ${\bar{J}}_{i-1,\ell}^t$ given by,
\begin{gather}
{\bar{u}}_{i-1,\ell}^t=\frac{A_{32,i-1}{\bar{S}}_{2i-1,\ell}^{xt}-\ A_{52}{\bar{S}}_{2i,\ell}^{xt}+\ A_{31,i-1}{\bar{u}}_{i-1,\ell}^{xt}-\ A_{51}{\bar{u}}_{i,\ell}^{xt}}{A_{31,i-1}-A_{51}} 
\label{eq:028}
\end{gather}
\begin{gather}
	{\bar{J}}_{i-1,\ell}^t=\frac{-A_{32,i-1}A_{51}{\bar{S}}_{2i-1,\ell}^{xt}+\ {A_{31,i-1}A}_{52}{\bar{S}}_{2i,\ell}^{xt}+\ A_{31,i-1}A_{51}({\bar{u}}_{i,\ell}^{xt}-{\bar{u}}_{i-1,\ell}^{xt})}{A_{31,i-1}-A_{51}} 
\label{eq:029}
\end{gather}
It is  crucial to note that, for each node ($i,\ell$), we have three equations (\eqns\ref{eq:022}, \ref{eq:026} and \ref{eq:027}) but six unknowns (${\bar{u}}_{i,\ell}^x$, ${\bar{u}}_{i,\ell}^t$, ${\bar{J}}_{i,\ell}^t$, ${\bar{u}}_{i,\ell}^{xt}$, ${\bar{S}}_{1i,\ell}^{xt}$ and ${\bar{S}}_{2i,\ell}^{xt}$). Given that the scheme must rely solely on cell-centered values, it becomes necessary to eliminate the surface-averaged terms from the final set of algebraic equations. To achieve this, three additional constraint conditions are required.
\subsubsection{Constraint equations}
\label{sec:2.1.5}
In order to successfully close the system of equations, we define three physically realistic constraint conditions in this section. Averaging the original PDE across all dimensions while taking the concept of pseudo-source terms into account yields the first constraint condition. Using the definitions of surface-averaged variables and carefully adjusting the pseudo-source terms, we systematically establish two more constraint conditions. This procedure creates a crucial connection between the surface-averaged terms and the pseudo-sources.

The first constraint equation is obtained by integrating the original Burgers' equation (\eqn\ref{eq:001}) over the node ($i,\ell$) using the double integration operator, $I_{xt}(a_i,\tau_{\ell})$. 
By incorporating the definitions of ${\bar{S}}_{1i,\ell}^{xt}$ from \eqn\eqref{eq:010} and ${\bar{S}}_{2i,\ell}^{xt}$ from \eqn\eqref{eq:011}, a relationship between ${\bar{S}}_{1i,\ell}^{xt}$ and ${\bar{S}}_{2i,\ell}^{xt}$ is established, expressed as follows.
\begin{gather}
	{\bar{S}}_{1i,\ell}^{xt}={\bar{S}}_{2i,\ell}^{xt}
\label{eq:030}
\end{gather}
The two additional constraint conditions are derived from the definition of the pseudo-source terms. The second constraint equation is obtained by averaging \eqn\eqref{eq:008} employing the integration operator $I_{t}(\tau_{\ell})$, which yields,
\begin{gather}
{\bar{S}}_{1i,\ell}^{xt} = \frac{1}{2a_i}\int_{-a_i}^{+a_i}\left(\frac{1}{Re}\frac{\mathrm{d}^2{\bar{u}}^t(x)}{\mathrm{d}x^2}-{\bar{u}}_{i,\ell}^0\frac{\mathrm{d}{\bar{u}}^t\left(x\right)}{\mathrm{d}x}\right)\mathrm{d}x  
\label{eq:031}
\end{gather}
\eqn\eqref{eq:031} can be rewritten as follows.
\begin{gather}
{\bar{S}}_{1i,\ell}^{xt} = \frac{1}{2a_i}\int_{-a_i}^{+a_i}\left(-\frac{\mathrm{d}{\bar{J}}^t\left(x\right)}{\mathrm{d}x}-{\bar{u}}_{i,\ell}^0\frac{\mathrm{d}{\bar{u}}^t\left(x\right)}{\mathrm{d}x}\right)\mathrm{d}x , \qquad\text{where}\qquad  {\bar{J}}^t(x)=-\frac{1}{Re}\frac{\mathrm{d}{\bar{u}}^t(x)}{\mathrm{d}x}
\label{eq:032}
\end{gather}
where ${\bar{J}}^t(x)$ is the flux. 
It is then reformulated in terms of the surface-averaged velocity and flux, leading to the second constraint condition, expressed as follows.
\begin{gather}
{\bar{S}}_{1i,\ell}^{xt}=-\frac{{\bar{J}}_{i,\ell}^t-{\bar{J}}_{i-1,\ell}^t}{2a_i}-{\bar{u}}_{i,\ell}^0\frac{{\bar{u}}_{i,\ell}^t-{\bar{u}}_{i-1,\ell}^t}{2a_i} 
\label{eq:033}
\end{gather}
Likewise, applying the definition of the second pseudo-source term from \eqn\eqref{eq:009}, the third constraint equation can be derived as follows.
\begin{gather}
{\bar{S}}_{2i,\ell}^{xt}=\frac{{\bar{u}}_{i,\ell}^x-{\bar{u}}_{i,\ell-1}^x}{2\tau_{\ell}} 
\label{eq:034}
\end{gather}
The above six equations (\eqns\ref{eq:022}, \ref{eq:026}, \ref{eq:027}, \ref{eq:030}, \ref{eq:033}, \ref{eq:034}) govern six unknowns (${\bar{u}}_{i,\ell}^x$, ${\bar{u}}_{i,\ell}^t$, ${\bar{J}}_{i,\ell}^t$, ${\bar{u}}_{i,\ell}^{xt}$, ${\bar{S}}_{1i,\ell}^{xt}$, ${\bar{S}}_{2i,\ell}^{xt}$). By eliminating the surface-averaged terms, the final set of algebraic equation can now be readily derived.
\subsubsection{Set of discrete equations}
\label{sec:2.1.6}
The first algebraic equation is derived by substituting the time-averaged velocity expressions (\eqns\ref{eq:026} and \ref{eq:028}) and the flux expressions (\eqns\ref{eq:027} and \ref{eq:029}) into the second constraint equation (\eqn\ref{eq:033}) as follows.
\begin{gather}
	{\bar{S}}_{1i,\ell}^{xt}=F_{31}{\bar{S}}_{2i,\ell}^{xt}+F_{32}{\bar{S}}_{2i-1,\ell}^{xt}+F_{33}{\bar{S}}_{2i+1,\ell}^{xt}+F_{34}{\bar{u}}_{i,\ell}^{xt}+F_{35}{\bar{u}}_{i-1,\ell}^{xt}+F_{36}{\bar{u}}_{i+1,\ell}^{xt} 
\label{eq:035}
\end{gather}
Further, substituting the space-averaged velocity expressions (\eqns\ref{eq:022} and \ref{eq:023}) into the third constraint equation (\eqn\ref{eq:034}) yields the second algebraic equation as follows.
\begin{gather}
	{\bar{S}}_{2i,\ell}^{xt}=\frac{{\bar{S}}_{1i,\ell}^{xt}-{\bar{S}}_{1i,\ell-1}^{xt}}{2}+\frac{{\bar{u}}_{i,\ell}^{xt}-{\bar{u}}_{i,\ell-1}^{xt}}{2\tau_{\ell}} 
\label{eq:036}
\end{gather}
By substituting the expressions for ${\bar{S}}_{1i,\ell}^{xt}$ and ${\bar{S}}_{2i,\ell}^{xt}$  (\eqns\ref{eq:035} and \ref{eq:036}) into the first constraint equation (\eqn\ref{eq:030}), the third algebraic equation is obtained as follows.
\begin{gather}
	{\bar{u}}_{i,\ell}^{xt}=F_{51}{\bar{S}}_{2i,\ell}^{xt}+F_{52}{\bar{S}}_{2i-1,\ell}^{xt}+F_{53}{\bar{S}}_{2i+1,\ell}^{xt}+F_{54}{\bar{S}}_{1i,\ell}^{xt}+F_{55}{\bar{S}}_{1i,\ell-1}^{xt}+F_{56}{\bar{u}}_{i,\ell-1}^{xt}+F_{57}{\bar{u}}_{i-1,\ell}^{xt}+F_{58}{\bar{u}}_{i+1,\ell}^{xt} 
\label{eq:037}
\end{gather}
The coefficients $F$’s are detailed in \ref{app:A}. The discrete unknowns for the cell ($i,\ell$) include the variables averaged at the cell centers (${\bar{u}}_{i,\ell}^{xt}$, ${\bar{S}}_{1i,\ell}^{xt}$, ${\bar{S}}_{2i,\ell}^{xt}$).
\subsubsection{Approximation of the non-linear convective velocity}
\label{sec:2.1.7}
In the context of Burgers' equation (\eqn\ref{eq:001}), the velocity $u(x,t)$ in the convection term $u(x,t)\frac{\partial u(x,t)}{\partial x}$ is an unknown function, which introduces non-linearity into the equation. As a result, during the implementation of the key steps in the development process, assumptions must be made regarding the approximation of the convective velocity $u(x,t)$ in the convection term. In MNIM \citep{23Rizwan_uddin_1997}, the convective velocity at a node ($i,\ell$), denoted as ${\bar{u}}_{i,\ell}^0$, is approximated by averaging the surface-averaged velocities at the current time step, as follows.
\begin{gather}
	{\bar{u}}_{i,\ell}^0=\frac{{\bar{u}}_{i,\ell}^t+{\bar{u}}_{i-1,\ell}^t}{2} 
\label{eq:038}
\end{gather} 
In the context of MCCNIM, we can adopt the same definition provided in MNIM, as outlined in \eqn\eqref{eq:038}.  However, it is important to highlight that MCCNIM employs cell-centered values in the final discretized equation, unlike MNIM, which uses surface-averaged velocities (${\bar{u}}_{i,\ell}^t$ and ${\bar{u}}_{i-1,\ell}^t$) to approximate the convective velocity. To ensure consistency, applying this approximation in MCCNIM requires reformulating convective velocity in terms of cell-centered values. Accordingly, the convective velocity in MCCNIM can be similarly defined using cell-centered velocities as follows.
\begin{gather}
	{\bar{u}}_{i,\ell}^0=\frac{{\bar{u}}_{i-1,\ell}^{xt}+{\bar{u}}_{i,\ell}^{xt}+{\bar{u}}_{i+1,\ell}^{xt}}{3} 
\label{eq:039}
\end{gather}
The approximation given by \eqn\eqref{eq:039} in MCCNIM can yield reasonable results comparable to those of the traditional MNIM. However, MCCNIM is recognized for its superior accuracy compared to the traditional NIM, as discussed in the literature \citep{17Ahmed_2024,18Ahmed_2023}. 
To fully realize the potential of MCCNIM for non-linear problems, a more precise approximation of the node-averaged velocity is required. This need arises from the fact that MCCNIM incorporates both cell-centered velocities and cell-centered pseudo-source terms in its final discretized equations. Therefore, approximating the convective velocity using only cell-centered velocities is inconsistent with the MCCNIM framework. To resolve this inconsistency, an expression for the convective velocity must be derived, accounting for both cell-centered values (velocities and pseudo-source terms). Hence, \eqns\eqref{eq:026} and \eqref{eq:028} are utilized to establish the connection between the surface-averaged and cell-centered values, and to express \eqn\eqref{eq:038} in terms of cell-centered values as follows.
\begin{gather}
	{\bar{u}}_{i,\ell}^0=F_{71}{\bar{S}}_{2i-1,\ell}^{xt}+F_{72}{\bar{S}}_{2i,\ell}^{xt}+F_{73}{\bar{S}}_{2i+1,\ell}^{xt}+F_{74}{\bar{u}}_{i-1,\ell}^{xt}+F_{75}{\bar{u}}_{i,\ell}^{xt}+F_{76}{\bar{u}}_{i+1,\ell}^{xt} 
\label{eq:040}
\end{gather}
where the coefficients $F$'s are given in \ref{app:A}. 
The results section includes a comparison of two definitions of node-averaged velocities (\eqns\ref{eq:039} and \ref{eq:040}), highlighting the necessity of deriving the expression for node-averaged velocity (\eqn\ref{eq:040}). This derivation is essential to ensure the superior accuracy of the developed MCCNIM scheme for addressing non-linear problems.
\subsubsection{Boundary conditions}
\label{sec:2.1.8}
\begin{figure}[!b]%
	\centering
	\includegraphics[width=0.8\linewidth]{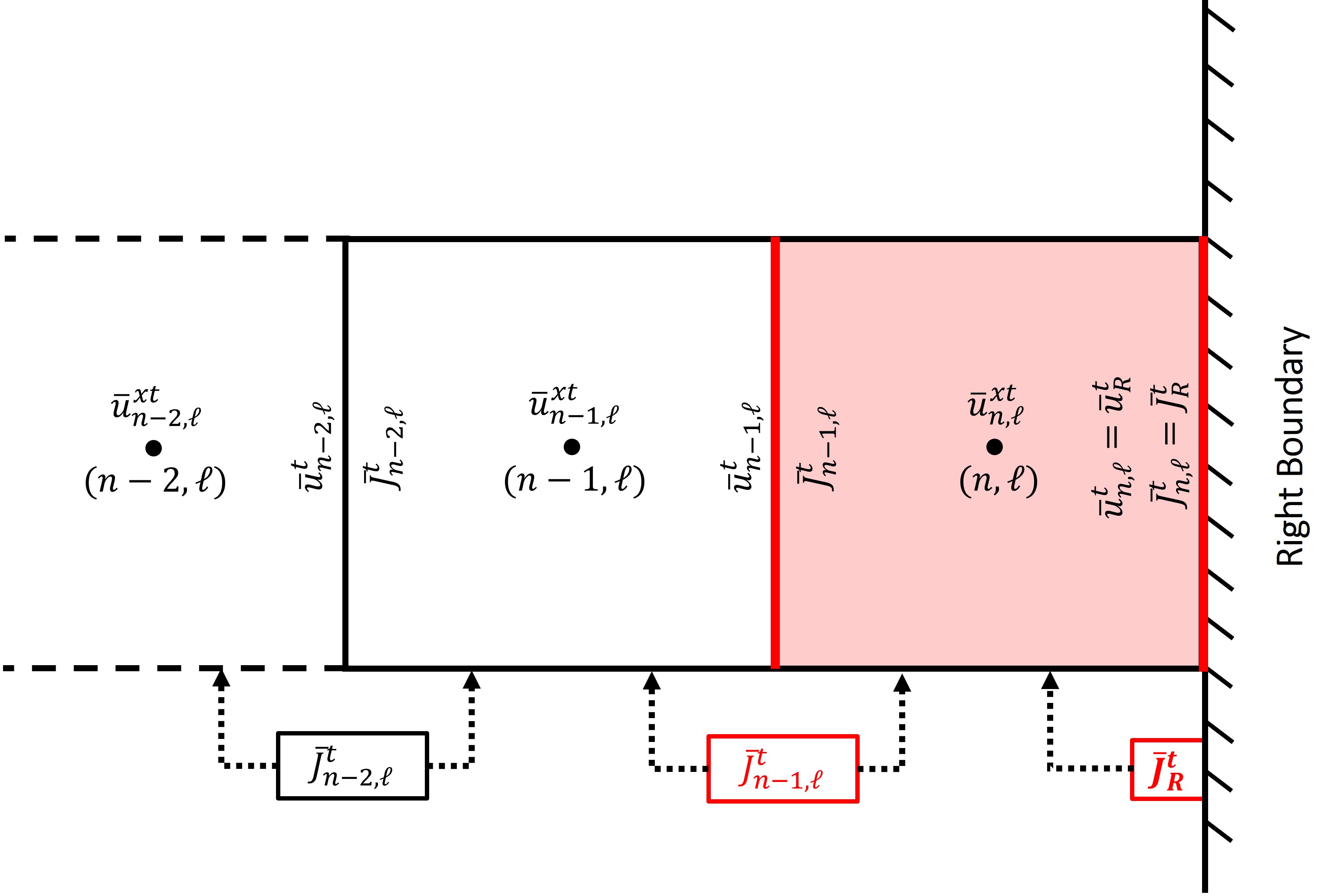}
	\caption{Stencil for right boundary node.}%
	\label{fig:03}%
\end{figure}
In the above-descrived scheme, the discrete unknowns are represented by dependent variables averaged at the center of each cell. This contrasts with the conventional NIM approach, where unknowns are defined on the cell surfaces \citep{23Rizwan_uddin_1997}. Consequently, since boundary conditions are typically specified at the edges or surfaces of the domain, the implementation of Dirichlet and Neumann boundary conditions in CCNIM schemes exhibits subtle differences compared to the standard NIM schemes. The expressions for the surface-averaged velocity ${\bar{u}}_{i,\ell}^t$ (\eqn\ref{eq:026}) and the surface-averaged flux ${\bar{J}}_{i,\ell}^t$  (\eqn\ref{eq:027}) for interior nodes reveal a dependence on the cell-centered velocities (${\bar{u}}_{i,\ell}^{xt}$ and ${\bar{u}}_{i+1,\ell}^{xt}$) of two adjacent nodes ($i,\ell$ and $i,\ell+1$), as detailed in \sect\ref{sec:2.1.4}. However, this relationship does not apply at boundary nodes because, at the boundaries, there is only one neighbouring node, as depicted in \fig\ref{fig:03} for the right boundary.

As illustrated in \fig\ref{fig:03}, the flux ${\bar{J}}_{n-1,\ell}^t$ at the right boundary node ($n,\ell$) is influenced by the cell-centered velocities, ${\bar{u}}_{n-1,\ell}^{xt}$ and ${\bar{u}}_{n,\ell}^{xt}$, from the adjacent nodes, ($n-1,\ell$) and ($n,\ell$). However, at the right boundary surface, where the flux is given by ${\bar{J}}_{n,\ell}^t={\bar{J}}_R^t$, only the node ($n,\ell$) remains, meaning it solely depends on the cell-centered velocity ${\bar{u}}_{n,\ell}^{xt}$. A detailed derivation for the right boundary is provided, and the equation for the left boundary node can be derived in an analogous manner.
\subsubsection*{(a) Dirichlet boundary conditions}
\label{sec:2.1.8a}
For the right boundary, referring to \eqn\eqref{eq:020}, ${\bar{J}}_{i,\ell}^t$ can be expressed as follows.
\begin{gather}
	{\bar{J}}_{i,\ell}^t=A_{31}\left({\bar{u}}_{i,\ell}^{xt}-{\bar{u}}_{i,\ell}^t\right)+A_{32}{\bar{S}}_{2i,\ell}^{xt} 
\label{eq:041}
\end{gather}
To simplify the explanation, we define the flux and velocity at the right boundary as ${\bar{J}}_R^t$ and ${\bar{u}}_R^t$, respectively. Using these definitions, \eqn\eqref{eq:041} can be reformulated as follows.
\begin{gather}
	{\bar{J}}_R^t=A_{31}\left({\bar{u}}_{i,\ell}^{xt}-{\bar{u}}_R^t\right)+A_{32}{\bar{S}}_{2i,\ell}^{xt} 
\label{eq:042}
\end{gather}
In a similar manner, by applying \eqn\eqref{eq:021}, the flux equation for the left boundary can be expressed as follows.
\begin{gather}
	{\bar{J}}_L^t=A_{51}\left({\bar{u}}_{i,\ell}^{xt}-{\bar{u}}_L^t\right)+A_{52}{\bar{S}}_{2i,\ell}^{xt} 
\label{eq:043}
\end{gather}
Here, the subscripts $L$ and $R$ denote the left and right boundary node, respectively. The process for deriving the final set of discrete equations for the right boundary node remains unchanged, except for the second constraint equation (\eqn\ref{eq:033}), which can be written for the right boundary as follows.
\begin{gather}
	{\bar{S}}_{1i,\ell}^{xt}=-\frac{{\bar{J}}_R^t-{\bar{J}}_{i-1,\ell}^t}{2a_i}-{\bar{u}}_{i,\ell}^0\frac{{\bar{u}}_R^t-{\bar{u}}_{i-1,\ell}^t}{2a_i} 
\label{eq:044}
\end{gather}
The first algebraic equation can now be derived by substituting the expression for ${\bar{u}}_{i-1,\ell}^t$  (\eqn\ref{eq:028}), ${\bar{J}}_{i-1,\ell}^t$ (\eqn\ref{eq:029}), and ${\bar{J}}_R^t$ (\eqn\ref{eq:042}) into the first constraint equation for the right boundary (\eqn\ref{eq:044}), resulting as follows. 
\begin{gather}
	{\bar{S}}_{1i,\ell}^{xt}=F_{31}^R{\bar{S}}_{2i,\ell}^{xt}+F_{32}^R{\bar{S}}_{2i-1,\ell}^{xt}+F_{34}^R{\bar{u}}_{i,\ell}^{xt}+F_{35}^R{\bar{u}}_{i-1,\ell}^{xt}+F_{36}^R{\bar{u}}_R^t 
\label{eq:045}
\end{gather}
It is important to note that the second algebraic equation will be the same for both interior and boundary nodes. By substituting the expressions for ${\bar{S}}_{2i,\ell}^{xt}$ (\eqn\ref{eq:035}) and ${\bar{S}}_{1i,\ell}^{xt}$  (\eqn\ref{eq:045}) into the third constraint equation (\eqn\ref{eq:030}), the third algebraic equation for the right boundary node is obtained as follows.
\begin{gather}
	{\bar{u}}_{i,\ell}^{xt}=F_{51}^R{\bar{S}}_{2i,\ell}^{xt}+F_{52}^R{\bar{S}}_{2i-1,\ell}^{xt}+F_{54}^R{\bar{S}}_{1i,\ell}^{xt}+F_{55}^R{\bar{S}}_{1i,\ell-1}^{xt}+F_{56}^R{\bar{u}}_{i,\ell-1}^{xt}+F_{57}^R{\bar{u}}_{i-1,\ell}^{xt}+F_{58}^R{\bar{u}}_R^t 
\label{eq:046}
\end{gather}
The superscript $R$ in all the coefficients signifies that they correspond to the right boundary node. The definitions of $F^R$ coefficients are provided in \ref{app:A}. The algebraic equations for the left boundary node can be derived using a similar approach.
\subsubsection*{(b) Neumann boundary conditions}
\label{sec:2.1.8b}
Deriving Neumann boundary conditions in MCCNIM is relatively straightforward. It involves setting the flux, either ${\bar{J}}_L^t$ or ${\bar{J}}_R^t$, to zero or applying the appropriate boundary condition as required. For instance, if the left boundary is insulated and the right boundary is governed by Dirichlet conditions, ${\bar{J}}_L^t$ is set to zero, while ${\bar{J}}_R^t$ is calculated using \eqn\eqref{eq:042}. Conversely, if the right boundary is insulated and the left boundary follows Dirichlet conditions, ${\bar{J}}_R^t$ is set to zero, and ${\bar{J}}_L^t$ is determined from \eqn\eqref{eq:043}. In cases where both boundaries are insulated, the flux at both surfaces is set to zero. Importantly, the method for obtaining the final set of algebraic equations remains identical to the approach outlined for Dirichlet boundary conditions.
\subsection{Two-dimensional coupled non-linear Burgers’ equations}
\label{sec:2.2}
The two-dimensional time-dependent coupled Burgers’ equation is written as follows.
\begin{gather}
	\frac{\partial u(x,y,t)}{\partial t}+u\left(x,y,t\right)\frac{\partial u\left(x,y,t\right)}{\partial x}+v\left(x,y,t\right)\frac{\partial u\left(x,y,t\right)}{\partial y}=\frac{1}{Re}\left(\frac{\partial^2u\left(x,y,t\right)}{\partial x^2}+\frac{\partial^2u\left(x,y,t\right)}{\partial y^2}\right)+f_x(x,y,t) 
\label{eq:047}
\end{gather}
\begin{gather}
	\frac{\partial v(x,y,t)}{\partial t}+u\left(x,y,t\right)\frac{\partial v\left(x,y,t\right)}{\partial x}+v\left(x,y,t\right)\frac{\partial v\left(x,y,t\right)}{\partial y}=\frac{1}{Re}\left(\frac{\partial^2v\left(x,y,t\right)}{\partial x^2}+\frac{\partial^2v\left(x,y,t\right)}{\partial y^2}\right)+f_y(x,y,t) 
\label{eq:048}
\end{gather}
The development procedure for two-dimensional Burgers’ equation will be the same as that for one-dimensional Burgers’ equation derived in \sect\ref{sec:2.1}. Some essential steps for a two-dimensional case are presented in this section for completeness.

The spatial domain defined by the independent variables ($x,y$) is divided into nodes of size $\Delta x\times\Delta y$, indexed by $i$ in the $x$-direction and $j$ in the $y$-direction. The solution is computed at discrete time intervals of $\Delta t$, indexed by $\ell$. For transverse integration, each space-time node is assigned dimensions, $\Delta x_i=2a_i$, $\Delta y_j=2b_j$, and $\Delta t_{\ell}=2\tau_{\ell}$, centered at the node. The local coordinate system for the node ($i,j,\ell$) is defined by (${-a}_i\le x_i\le{+a}_i$, ${-b}_j\le y_j\le{+b}_j$, ${-\tau}_{\ell}\le t_{\ell}\le{+\tau}_{\ell}$). 

\begin{figure}[!b]%
	\centering
	\includegraphics[width=0.75\linewidth]{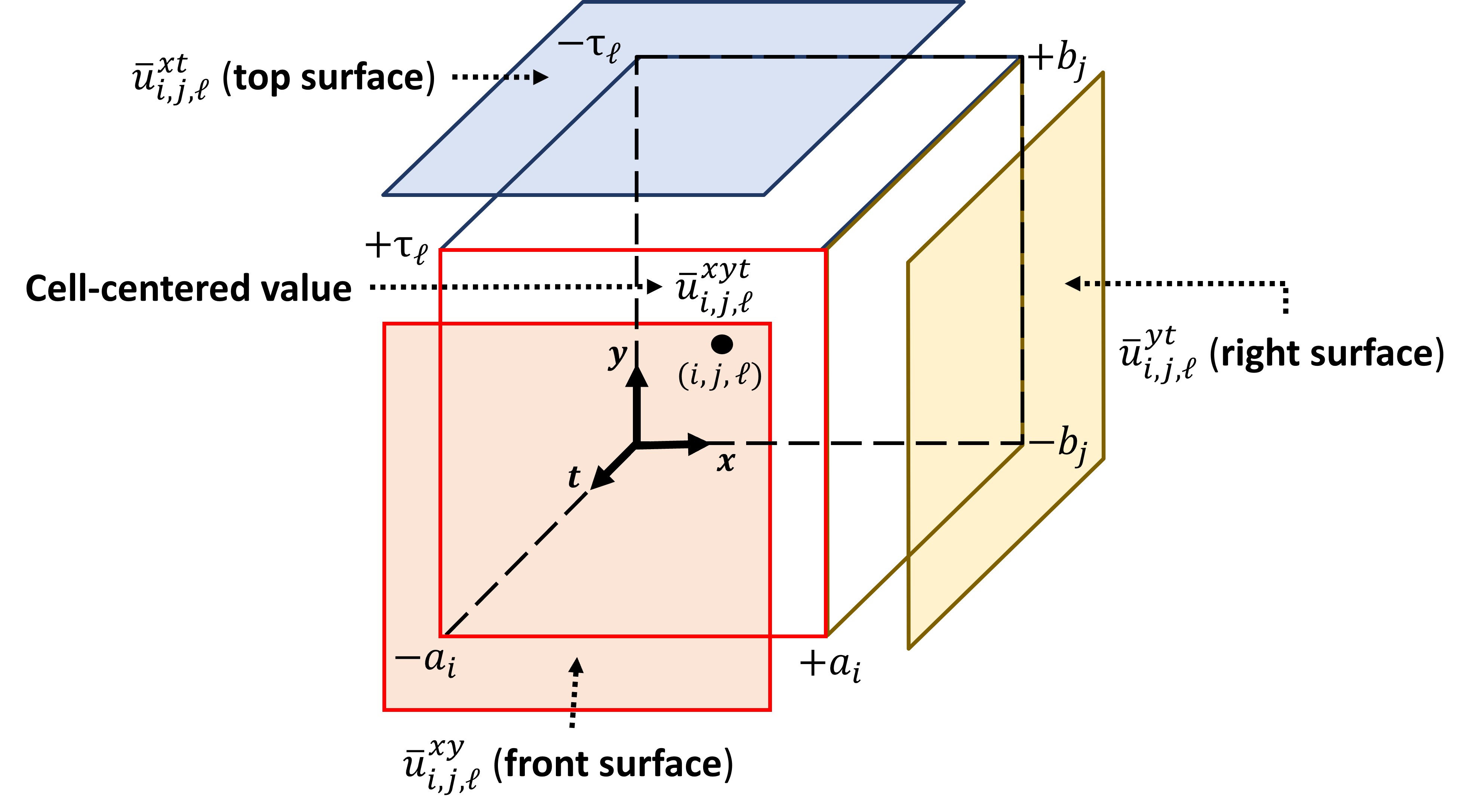}
	\caption{Local coordinate system and transverse-averaged quantities of node ($i,j,\ell$) for two-dimensional transient problems in MCCNIM.}%
	\label{fig:04}%
\end{figure}

\fig\ref{fig:04} illustrates the transverse-averaged quantities on various surfaces, highlighted in different colours, along with the local coordinate system.
Prior to the development of the numerical scheme, the Burgers’ equations (\eqns\ref{eq:047} and \ref{eq:048}) are reformulated in terms of the local coordinate system as follows.
\begin{gather}
	\frac{\partial u(x,y,t)}{\partial t}+{\bar{u}}^0\frac{\partial u\left(x,y,t\right)}{\partial x}+{\bar{v}}^0\frac{\partial u\left(x,y,t\right)}{\partial y}=\frac{1}{Re}\left(\frac{\partial^2u\left(x,y,t\right)}{\partial x^2}+\frac{\partial^2u\left(x,y,t\right)}{\partial y^2}\right)+f_x(x,y,t) 
\label{eq:049}
\end{gather}
\begin{gather}
\frac{\partial v(x,y,t)}{\partial t}+{\bar{u}}^0\frac{\partial v\left(x,y,t\right)}{\partial x}+{\bar{v}}^0\frac{\partial v\left(x,y,t\right)}{\partial y}=\frac{1}{Re}\left(\frac{\partial^2v\left(x,y,t\right)}{\partial x^2}+\frac{\partial^2v\left(x,y,t\right)}{\partial y^2}\right)+f_y(x,y,t) 
\label{eq:050}
\end{gather}
Here, ${\bar{u}}^0$ and ${\bar{v}}^0$ represents the approximation of the convective (node-averaged) velocities in the $x-$ and $y-$directions, respectively. The definitions of these velocities for the two-dimensional case are provided later in the derivation.

By performing transverse integration in each direction and truncating the pseudo-source terms to zeroth order, analogous to the one-dimensional case discussed in \sect\ref{sec:2.1.1}, the resulting set of ODEs for the node ($i,j,\ell$) in the context of the two-dimensional Burgers’ equations (\eqns\ref{eq:047} and \ref{eq:048}) is presented as follows.
\begin{gather}
	\frac{1}{Re}\frac{\mathrm{d}^2{\bar{u}}^{yt}(x)}{\mathrm{d}x^2}-{\bar{u}}_{i,j,\ell}^0\frac{\mathrm{d}{\bar{u}}^{yt}(x)}{\mathrm{d}x}={\bar{S}}_{3i,j,\ell}^{xyt} 
\label{eq:051}
\end{gather}
\begin{gather}
	\frac{1}{Re}\frac{\mathrm{d}^2{\bar{u}}^{xt}\left(y\right)}{\mathrm{d}y^2}-{\bar{v}}_{i,j,\ell}^0\frac{\mathrm{d}{\bar{u}}^{xt}(y)}{\mathrm{d}y}={\bar{S}}_{2i,j,\ell}^{xyt} 
\label{eq:052}
\end{gather}
\begin{gather}
	\frac{\mathrm{d}{\bar{u}}^{xy}\left(t\right)}{\mathrm{d}t}={\bar{S}}_{1i,j,\ell}^{xyt} 	
\label{eq:053}
	\end{gather}
	\begin{gather}
	\frac{1}{Re}\frac{\mathrm{d}^2{\bar{v}}^{yt}(x)}{\mathrm{d}x^2}-{\bar{u}}_{i,j,\ell}^0\frac{\mathrm{d}{\bar{v}}^{yt}(x)}{\mathrm{d}x}={\bar{S}}_{6i,j,\ell}^{xyt}  
\label{eq:054}
\end{gather}
	\begin{gather}
	\frac{1}{Re}\frac{\mathrm{d}^2{\bar{v}}^{xt}\left(y\right)}{\mathrm{d}y^2}-{\bar{v}}_{i,j,\ell}^0\frac{\mathrm{d}{\bar{v}}^{xt}(y)}{\mathrm{d}y}={\bar{S}}_{5i,j,\ell}^{xyt}
\label{eq:055}
\end{gather}
	\begin{gather}
	\frac{\mathrm{d}{\bar{v}}^{xy}\left(t\right)}{\mathrm{d}t}={\bar{S}}_{4i,j,\ell}^{xyt}
\label{eq:056}
\end{gather}
The three ODEs (\eqns\ref{eq:051} - \ref{eq:053}) are derived from the $x$-momentum equation (\eqn\ref{eq:047}), while the subsequent three ODEs (\eqns\ref{eq:054} - \ref{eq:056}) are obtained from the $y$-momentum equation (\eqn\ref{eq:048}), both using TIP as discussed in \sect\ref{sec:2.1.1}. Furthermore, all the assumptions and definitions outlined in \sect\ref{sec:2.1} for developing the one-dimensional scheme are equally applicable to the two-dimensional case. 
The definitions of pseudo-source terms ({${\bar{S}}_{1i,j,\ell}^{xyt}$, ${\bar{S}}_{2i,j,\ell}^{xyt}$, ${\bar{S}}_{3i,j,\ell}^{xyt}$, ${\bar{S}}_{4i,j,\ell}^{xyt}$, ${\bar{S}}_{5i,j,\ell}^{xyt}$, and ${\bar{S}}_{6i,j,\ell}^{xyt}$) are given as
\begin{gather}
	\begin{split}
		{\bar{S}}_{3i,j,\ell}^{xyt}=\frac{1}{2\tau_{\ell}}\int_{-\tau_{\ell}}^{+\tau_{\ell}}{\frac{\mathrm{d}{\bar{u}}^{xy}\left(t\right)}{\mathrm{d}t}\mathrm{d}t}-\frac{1}{2b_j}\int_{-b_j}^{+b_j}\left(\frac{1}{Re}\frac{\mathrm{d}^2{\bar{u}}^{xt}\left(y\right)}{\mathrm{d}y^2}-{\bar{v}}_{i,j,\ell}^0\frac{\mathrm{d}{\bar{u}}^{xt}\left(y\right)}{\mathrm{d}y}\right)\mathrm{d}y-{\bar{f}}_{xi,j,\ell}^{xyt}
	\end{split}
\label{eq:057}
\end{gather}
\begin{gather}
	\begin{split}
		{\bar{S}}_{2i,j,\ell}^{xyt}=\frac{1}{2\tau_{\ell}}\int_{-\tau_{\ell}}^{+\tau_{\ell}}{\frac{\mathrm{d}{\bar{u}}^{xy}\left(t\right)}{\mathrm{d}t}\mathrm{d}t}-\frac{1}{2a_i}\int_{-a_i}^{+a_i}\left(\frac{1}{Re}\frac{\mathrm{d}^2{\bar{u}}^{yt}\left(x\right)}{\mathrm{d}x^2}-{\bar{u}}_{i,j,\ell}^0\frac{\mathrm{d}{\bar{u}}^{yt}\left(x\right)}{\mathrm{d}x}\right)\mathrm{d}x-{\bar{f}}_{xi,j,\ell}^{xyt} 
	\end{split}
\label{eq:058}
\end{gather}
\begin{gather}
	\begin{split}
		{\bar{S}}_{1i,j,\ell}^{xyt}=\frac{1}{2a_i}\int_{-a_i}^{+a_i}\left(\frac{1}{Re}\frac{\mathrm{d}^2{\bar{u}}^{yt}\left(x\right)}{\mathrm{d}x^2}-{\bar{u}}_{i,j,\ell}^0\frac{\mathrm{d}{\bar{u}}^{yt}\left(x\right)}{\mathrm{d}x}\right)\mathrm{d}x+\frac{1}{2b_j}\int_{-b_j}^{+b_j}\left(\frac{1}{Re}\frac{\mathrm{d}^2{\bar{u}}^{xt}\left(y\right)}{\mathrm{d}y^2}-{\bar{v}}_{i,j,\ell}^0\frac{\mathrm{d}{\bar{u}}^{xt}\left(y\right)}{\mathrm{d}y}\right)\mathrm{d}y+{\bar{f}}_{xi,j,\ell}^{xyt}
	\end{split}
\label{eq:059}
\end{gather}
\begin{gather}
	\begin{split}
		{\bar{S}}_{6i,j,\ell}^{xyt}=\frac{1}{2\tau_{\ell}}\int_{-\tau_{\ell}}^{+\tau_{\ell}}{\frac{\mathrm{d}{\bar{v}}^{xy}\left(t\right)}{\mathrm{d}t}\mathrm{d}t}-\frac{1}{2b_j}\int_{-b_j}^{+b_j}\left(\frac{1}{Re}\frac{\mathrm{d}^2{\bar{v}}^{xt}\left(y\right)}{\mathrm{d}y^2}-{\bar{v}}_{i,j,\ell}^0\frac{\mathrm{d}{\bar{v}}^{xt}\left(y\right)}{\mathrm{d}y}\right)\mathrm{d}y-{\bar{f}}_{yi,j,\ell}^{xyt}
	\end{split}
\label{eq:060}
\end{gather}
\begin{gather}
	\begin{split}
		{\bar{S}}_{5i,j,\ell}^{xyt}=\frac{1}{2\tau_{\ell}}\int_{-\tau_{\ell}}^{+\tau_{\ell}}{\frac{\mathrm{d}{\bar{v}}^{xy}\left(t\right)}{\mathrm{d}t}\mathrm{d}t}-\frac{1}{2a_i}\int_{-a_i}^{+a_i}\left(\frac{1}{Re}\frac{\mathrm{d}^2{\bar{v}}^{yt}\left(x\right)}{\mathrm{d}x^2}-{\bar{u}}_{i,j,\ell}^0\frac{\mathrm{d}{\bar{v}}^{yt}\left(x\right)}{\mathrm{d}x}\right)\mathrm{d}x-{\bar{f}}_{yi,j,\ell}^{xyt}
	\end{split}
\label{eq:061}
\end{gather}
\begin{gather}
	\begin{split}
		{\bar{S}}_{4i,j,\ell}^{xyt}=\frac{1}{2a_i}\int_{-a_i}^{+a_i}\left(\frac{1}{Re}\frac{\mathrm{d}^2{\bar{v}}^{yt}\left(x\right)}{\mathrm{d}x^2}-{\bar{u}}_{i,j,\ell}^0\frac{\mathrm{d}{\bar{v}}^{yt}\left(x\right)}{\mathrm{d}x}\right)\mathrm{d}x+\frac{1}{2b_j}\int_{-b_j}^{+b_j}\left(\frac{1}{Re}\frac{\mathrm{d}^2{\bar{v}}^{xt}\left(y\right)}{\mathrm{d}y^2}-{\bar{v}}_{i,j,\ell}^0\frac{\mathrm{d}{\bar{v}}^{xt}\left(y\right)}{\mathrm{d}y}\right)\mathrm{d}y+{\bar{f}}_{yi,j,\ell}^{xyt} 
	\end{split}
\label{eq:062}
\end{gather}
Notably, the above-defined pseudo-source terms are already averaged over the entire node and are presented in a manner that facilitates the straightforward derivation of the constraint equations. Upon rearranging the pseudo-source terms, the six constraint equations are formulated in terms of surface-averaged quantities as follows.
\begin{gather}
	\begin{split}
		{\bar{S}}_{3i,j,\ell}^{xyt}=\frac{{\bar{u}}_{i,j,\ell}^{xy}-{\bar{u}}_{i,j,\ell-1}^{xy}}{2\tau_{\ell}}+\frac{{\bar{J}}_{xi,j,\ell}^{xt}-{\bar{J}}_{xi,j-1,\ell}^{xt}}{2b_j}+{\bar{v}}_{i,j,\ell}^0\frac{{\bar{u}}_{i,j,\ell}^{xt}-{\bar{u}}_{i,j-1,\ell}^{xt}}{2b_j}-{\bar{f}}_{xi,j,\ell}^{xyt}
	\end{split}
\label{eq:063}
\end{gather}
\begin{gather}
	\begin{split}
		{\bar{S}}_{2i,j,\ell}^{xyt}=\frac{{\bar{u}}_{i,j,\ell}^{xy}-{\bar{u}}_{i,j,\ell-1}^{xy}}{2\tau_{\ell}}+\frac{{\bar{J}}_{xi,j,\ell}^{yt}-{\bar{J}}_{xi-1,j,\ell}^{yt}}{2a_i}+{\bar{u}}_{i,j,\ell}^0\frac{{\bar{u}}_{i,j,\ell}^{yt}-{\bar{u}}_{i-1,j,\ell}^{yt}}{2a_i}-{\bar{f}}_{xi,j,\ell}^{xyt}
	\end{split}
\label{eq:064}
\end{gather}
\begin{gather}
	\begin{split}
		{\bar{S}}_{1i,j,\ell}^{xyt}=-\frac{{\bar{J}}_{xi,j,\ell}^{yt}-{\bar{J}}_{xi-1,j,\ell}^{yt}}{2a_i}-{\bar{u}}_{i,j,\ell}^0\frac{{\bar{u}}_{i,j,\ell}^{yt}-{\bar{u}}_{i-1,j,\ell}^{yt}}{2a_i}-\frac{{\bar{J}}_{xi,j,\ell}^{xt}-{\bar{J}}_{xi,j-1,\ell}^{xt}}{2b_j}-{\bar{v}}_{i,j,\ell}^0\frac{{\bar{u}}_{i,j,\ell}^{xt}-{\bar{u}}_{i,j-1,\ell}^{xt}}{2b_j}+{\bar{f}}_{xi,j,\ell}^{xyt} 	
	\end{split}
\label{eq:065}
\end{gather}
\begin{gather}
	\begin{split}
		{\bar{S}}_{6i,j,\ell}^{xyt}=\frac{{\bar{v}}_{i,j,\ell}^{xy}-{\bar{v}}_{i,j,\ell-1}^{xy}}{2\tau_{\ell}}+\frac{{\bar{J}}_{yi,j,\ell}^{xt}-{\bar{J}}_{yi,j-1,\ell}^{xt}}{2b_j}+{\bar{v}}_{i,j,\ell}^0\frac{{\bar{v}}_{i,j,\ell}^{xt}-{\bar{v}}_{i,j-1,\ell}^{xt}}{2b_j}-{\bar{f}}_{yi,j,\ell}^{xyt}
	\end{split}
\label{eq:066}
\end{gather}
\begin{gather}
	\begin{split}
		{\bar{S}}_{5i,j,\ell}^{xyt}=\frac{{\bar{v}}_{i,j,\ell}^{xy}-{\bar{v}}_{i,j,\ell-1}^{xy}}{2\tau_{\ell}}+\frac{{\bar{J}}_{yi,j,\ell}^{yt}-{\bar{J}}_{yi-1,j,\ell}^{yt}}{2a_i}+{\bar{u}}_{i,j,\ell}^0\frac{{\bar{v}}_{i,j,\ell}^{yt}-{\bar{v}}_{i-1,j,\ell}^{yt}}{2a_i}-{\bar{f}}_{yi,j,\ell}^{xyt}
	\end{split}
\label{eq:067}
\end{gather}
\begin{gather}
	\begin{split}
		{\bar{S}}_{4i,j,\ell}^{xyt}=-\frac{{\bar{J}}_{yi,j,\ell}^{yt}-{\bar{J}}_{yi-1,j,\ell}^{yt}}{2a_i}-{\bar{u}}_{i,j,\ell}^0\frac{{\bar{v}}_{i,j,\ell}^{yt}-{\bar{v}}_{i-1,j,\ell}^{yt}}{2a_i}-\frac{{\bar{J}}_{yi,j,\ell}^{xt}-{\bar{J}}_{yi,j-1,\ell}^{xt}}{2b_j}-{\bar{v}}_{i,j,\ell}^0\frac{{\bar{v}}_{i,j,\ell}^{xt}-{\bar{v}}_{i,j-1,\ell}^{xt}}{2b_j}+{\bar{f}}_{yi,j,\ell}^{xyt}
	\end{split}
\label{eq:068}
\end{gather}
Furthermore, by averaging the original Burgers’ equations (\eqns\ref{eq:047} and \ref{eq:048}) over the node ($i,j,\ell$) using the integral operator defined as
\begin{gather*}
	I_{xyt}(g,l,h) = \frac{1}{8glh}\int_{-g}^{+g}\int_{-l}^{+l}\int_{-h}^{+h} F \mathrm{d}t \mathrm{d}y \mathrm{d}x,\quad\text{where}\quad	g = a_i, l = b_j, h = \tau_{\ell}
\end{gather*} 
and realising the definition of pseudo source terms (\eqns\eqref{eq:051} - \eqref{eq:056}), the two additional constraint equations are derived as follows.
\begin{gather}
	{\bar{S}}_{1i,j,\ell}^{xyt}={\bar{S}}_{2i,j,\ell}^{xyt}+{\bar{S}}_{3i,j,\ell}^{xyt}+{\bar{f}}_{xi,j,\ell}^{xyt} 
\label{eq:069}
\end{gather}
\begin{gather}
	{\bar{S}}_{4i,j,\ell}^{xyt}={\bar{S}}_{5i,j,\ell}^{xyt}+{\bar{S}}_{6i,j,\ell}^{xyt}+{\bar{f}}_{yi,j,\ell}^{xyt}  	
\label{eq:070}
\end{gather}
To eliminate the surface-averaged terms in \eqns\eqref{eq:063} - \eqref{eq:070}, we derive expressions for each surface-averaged variable in terms of cell-centered values. This requires solving the six ODEs (\eqns\ref{eq:051} - \ref{eq:056}) analytically. The solution procedure for these ODEs follows the methodology outlined in \sect\ref{sec:2.1.2}. Specifically, the space-averaged ODEs (\eqns\ref{eq:053} and \ref{eq:056}) are solved using an approach analogous to that for \eqn\eqref{eq:010} in the one-dimensional case, yielding expressions for ${\bar{u}}_{i,j,\ell}^{xy}$ and ${\bar{v}}_{i,j,\ell}^{xy}$. Similarly, the remaining four space-time averaged ODEs (\eqns\ref{eq:051}, \ref{eq:052}, \ref{eq:054}, \ref{eq:055}) are solved using a process similar to the solution of \eqn\eqref{eq:011} in the one-dimensional case. 

However, it is important to note a key distinction in boundary conditions between the one-dimensional and two-dimensional scenarios for the solution of these remaining four ODEs. In the one-dimensional case, line-averaged velocities (${\bar{u}}_{i,\ell}^t$) and flux (${\bar{J}}_{i,\ell}^t$) serve as local boundary conditions. In contrast, the two-dimensional case requires plane-averaged velocities (${\bar{u}}_{i,j,\ell}^{xt}$, ${\bar{u}}_{i,j,\ell}^{yt}$, ${\bar{v}}_{i,j,\ell}^{xt}$, ${\bar{v}}_{i,j,\ell}^{yt}$) and fluxes (${\bar{J}}_{xi,j,\ell}^{xt}$, ${\bar{J}}_{xi,j,\ell}^{yt}$, ${\bar{J}}_{yi,j,\ell}^{xt}$, ${\bar{J}}_{yi,j,\ell}^{yt}$) as local boundary conditions for solving these ODEs analytically. Subsequently, the continuity condition, as discussed in \sect\ref{sec:2.1.4}, is applied in a manner analogous to the one-dimensional case to ensure consistency at the shared edges of adjacent nodes. This process allows for the derivation of expressions for surface-averaged velocities (${\bar{u}}_{i,j,\ell}^{xy}$, ${\bar{u}}_{i,j,\ell}^{xt}$, ${\bar{u}}_{i,j,\ell}^{yt}$, ${\bar{v}}_{i,j,\ell}^{xy}$, ${\bar{v}}_{i,j,\ell}^{xt}$, ${\bar{v}}_{i,j,\ell}^{yt}$) and fluxes (${\bar{J}}_{xi,j,\ell}^{xt}$, ${\bar{J}}_{xi,j,\ell}^{yt}$, ${\bar{J}}_{yi,j,\ell}^{xt}$, ${\bar{J}}_{yi,j,\ell}^{yt}$) in the two-dimensional framework, as detailed below.
\begin{gather}
	\begin{split}
		{\bar{u}}_{i,j,\ell}^{xy}={\bar{u}}_{i,j,\ell}^{xyt}+\tau_{\ell}{\bar{S}}_{1i,j,\ell}^{xyt} 
	\end{split}
\label{eq:071}
\end{gather}
\begin{gather}
	\begin{split}
			{\bar{u}}_{i,j,\ell}^{yt}=\frac{A_{32}{\bar{S}}_{3i,j,\ell}^{xyt}-\ A_{52,i+1}{\bar{S}}_{3i+1,j,\ell}^{xyt}+\ A_{31}{\bar{u}}_{i,j,\ell}^{xyt}-\ A_{51,i+1}{\bar{u}}_{i+1,j,\ell}^{xyt}}{A_{31}-A_{51,i+1}} 
	\end{split}
\label{eq:072}
\end{gather}
\begin{gather}
	\begin{split}
			{\bar{u}}_{i,j,\ell}^{xt}=\frac{B_{32}{\bar{S}}_{2i,j,\ell}^{xyt}-\ B_{52,j+1}{\bar{S}}_{2i,j+1,\ell}^{xyt}+\ B_{31}{\bar{u}}_{i,j,\ell}^{xyt}-\ B_{51,j+1}{\bar{u}}_{i,j+1,\ell}^{xyt}}{B_{31}-B_{51,j+1}} 
	\end{split}
\label{eq:073}
\end{gather}
\begin{gather}
\begin{split}
	{\bar{J}}_{xi,j,\ell}^{yt}=\frac{-A_{32}A_{51,i+1}{\bar{S}}_{3i,j,\ell}^{xyt}+\ {A_{31}A}_{52,i+1}{\bar{S}}_{3i+1,j,\ell}^{xyt}+\ A_{31}A_{51,i+1}({\bar{u}}_{i+1,j,\ell}^{xyt}-{\bar{u}}_{i,j,\ell}^{xyt})}{A_{31}-A_{51,i+1}} 
\end{split}
\label{eq:074}
\end{gather}
\begin{gather}
\begin{split}
	{\bar{J}}_{xi,j,\ell}^{xt}=\frac{-B_{32}B_{51,j+1}{\bar{S}}_{2i,j,\ell}^{xyt}+\ {B_{31}B}_{52,j+1}{\bar{S}}_{2i,j+1,\ell}^{xyt}+\ B_{31}B_{51,j+1}({\bar{u}}_{i,j+1,\ell}^{xyt}-{\bar{u}}_{i,j,\ell}^{xyt})}{B_{31}-B_{51,j+1}} 
\end{split}
\label{eq:075}
\end{gather}
\begin{gather}
\begin{split}
		{\bar{v}}_{i,j,\ell}^{xy}={\bar{v}}_{i,j,\ell}^{xyt}+\tau_{\ell}{\bar{S}}_{4i,j,\ell}^{xyt}
\end{split}
\label{eq:076}
\end{gather}
\begin{gather}
	\begin{split}
			{\bar{v}}_{i,j,\ell}^{yt}=\frac{A_{32}{\bar{S}}_{6i,j,\ell}^{xyt}-\ A_{52,i+1}{\bar{S}}_{6i+1,j,\ell}^{xyt}+\ A_{31}{\bar{v}}_{i,j,\ell}^{xyt}-\ A_{51,i+1}{\bar{v}}_{i+1,j,\ell}^{xyt}}{A_{31}-A_{51,i+1}}
	\end{split}
\label{eq:077}
\end{gather}
\begin{gather}
	\begin{split}
			{\bar{v}}_{i,j,\ell}^{xt}=\frac{B_{32}{\bar{S}}_{5i,j,\ell}^{xyt}-\ B_{52,j+1}{\bar{S}}_{5i,j+1,\ell}^{xyt}+\ B_{31}{\bar{v}}_{i,j,\ell}^{xyt}-\ B_{51,j+1}{\bar{v}}_{i,j+1,\ell}^{xyt}}{B_{31}-B_{51,j+1}}
	\end{split}
\label{eq:078}
\end{gather}
\begin{gather}
	\begin{split}
			{\bar{J}}_{yi,j,\ell}^{yt}=\frac{-A_{32}A_{51,i+1}{\bar{S}}_{6i,j,\ell}^{xyt}+\ {A_{31}A}_{52,i+1}{\bar{S}}_{6i+1,j,\ell}^{xyt}+\ A_{31}A_{51,i+1}({\bar{v}}_{i+1,j,\ell}^{xyt}-{\bar{v}}_{i,j,\ell}^{xyt})}{A_{31}-A_{51,i+1}}
	\end{split}
\label{eq:079}
\end{gather}
\begin{gather}
	\begin{split}
			{\bar{J}}_{yi,j,\ell}^{xt}=\frac{-B_{32}B_{51,j+1}{\bar{S}}_{5i,j,\ell}^{xyt}+\ {B_{31}B}_{52,j+1}{\bar{S}}_{5i,j+1,\ell}^{xyt}+\ B_{31}B_{51,j+1}({\bar{v}}_{i,j+1,\ell}^{xyt}-{\bar{v}}_{i,j,\ell}^{xyt})}{B_{31}-B_{51,j+1}}
	\end{split}
\label{eq:080}
\end{gather}
The expressions for velocities and fluxes (${\bar{u}}_{i,j,\ell-1}^{xy}$, ${\bar{v}}_{i,j,\ell-1}^{xy}$, ${\bar{u}}_{i-1,j,\ell}^{yt}$, ${\bar{u}}_{i,j-1,\ell}^{xt}$, ${\bar{v}}_{i-1,j,\ell}^{yt}$, ${\bar{v}}_{i,j-1,\ell}^{xt}$, ${\bar{J}}_{xi-1,j,\ell}^{yt}$, ${\bar{J}}_{xi,j-1,\ell}^{xt}$, ${\bar{J}}_{yi-1,j,\ell}^{yt}$, ${\bar{J}}_{yi,j-1,\ell}^{xt}$) can be obtained by modifying the indexing in the relevant equations as follows: changing ($i,j,\ell$) to ($i,j,\ell-1$) in \eqns\eqref{eq:071} and \eqref{eq:076}, ($i,j,\ell$) to ($i-1,j,\ell$) in \eqns\eqref{eq:072}, \eqref{eq:074}, \eqref{eq:077} and \eqref{eq:079}, and ($i,j,\ell$) to ($i,j-1,\ell$) in \eqns\eqref{eq:073}, \eqref{eq:075}, \eqref{eq:078} and \eqref{eq:080}.
By substituting all these surface-averaged expressions derived from \eqns\eqref{eq:071} - \eqref{eq:080} into the constraint equations (\eqns\ref{eq:063} - \ref{eq:070}), we obtain eight algebraic equations per node in terms of cell-centered variables, as follows.
\begin{gather}
	\begin{split}
	{\bar{S}}_{3i,j,\ell}^{xyt}&=\frac{{\bar{S}}_{1i,j,\ell}^{xyt}-{\bar{S}}_{1i,j,\ell-1}^{xyt}}{2}+\frac{{\bar{u}}_{i,j,\ell}^{xyt}-{\bar{u}}_{i,j,\ell-1}^{xyt}}{2\tau_{\ell}}+F_{51}{\bar{S}}_{2i,j,\ell}^{xyt}+F_{52}{\bar{S}}_{2i,j-1,\ell}^{xyt}+F_{53}{\bar{S}}_{2i,j+1,\ell}^{xyt}\\&+F_{54}{\bar{u}}_{i,j,\ell}^{xyt}+F_{55}{\bar{u}}_{i,j-1,\ell}^{xyt}+F_{56}{\bar{u}}_{i,j+1,\ell}^{xyt}-{\bar{f}}_{xi,j,\ell}^{xyt} 
	\end{split}
\label{eq:081}
\end{gather}
\begin{gather}
	\begin{split}
	{\bar{S}}_{2i,j,\ell}^{xyt}&=\frac{{\bar{S}}_{1i,j,\ell}^{xyt}-{\bar{S}}_{1i,j,\ell-1}^{xyt}}{2}+\frac{{\bar{u}}_{i,j,\ell}^{xyt}-{\bar{u}}_{i,j,\ell-1}^{xyt}}{2\tau_{\ell}}+F_{31}{\bar{S}}_{3i,j,\ell}^{xyt}+F_{32}{\bar{S}}_{3i-1,j,\ell}^{xyt}+F_{33}{\bar{S}}_{3i+1,j,\ell}^{xyt}\\&+F_{34}{\bar{u}}_{i,j,\ell}^{xyt}+F_{35}{\bar{u}}_{i-1,j,\ell}^{xyt}+F_{36}{\bar{u}}_{i+1,j,\ell}^{xyt}-{\bar{f}}_{xi,j,\ell}^{xyt} 
	\end{split}
\label{eq:082}
\end{gather}
\begin{gather}
	\begin{split}
	{\bar{S}}_{1i,j,\ell}^{xyt}&=-\left(F_{31}{\bar{S}}_{3i,j,\ell}^{xyt}+F_{32}{\bar{S}}_{3i-1,j,\ell}^{xyt}+F_{33}{\bar{S}}_{3i+1,j,\ell}^{xyt}+F_{34}{\bar{u}}_{i,j,\ell}^{xyt}+F_{35}{\bar{u}}_{i-1,j,\ell}^{xyt}+F_{36}{\bar{u}}_{i+1,j,\ell}^{xyt}\right)\\&-\left(F_{51}{\bar{S}}_{2i,j,\ell}^{xyt}+F_{52}{\bar{S}}_{2i,j-1,\ell}^{xyt}+F_{53}{\bar{S}}_{2i,j+1,\ell}^{xyt}+F_{54}{\bar{u}}_{i,j,\ell}^{xyt}+F_{55}{\bar{u}}_{i,j-1,\ell}^{xyt}+F_{56}{\bar{u}}_{i,j+1,\ell}^{xyt}\right)+{\bar{f}}_{xi,j,\ell}^{xyt} 
	\end{split}
\label{eq:083}
\end{gather}
\begin{gather}
	\begin{split}
{\bar{u}}_{i,j,\ell}^{xyt}&=M_{51}{\bar{S}}_{3i,j,\ell}^{xyt}+M_{52}{\bar{S}}_{3i-1,j,\ell}^{xyt}+M_{53}{\bar{S}}_{3i+1,j,\ell}^{xyt}+M_{54}{\bar{u}}_{i-1,j,\ell}^{xyt}+M_{55}{\bar{u}}_{i+1,j,\ell}^{xyt}+N_{51}{\bar{S}}_{2i,j,\ell}^{xyt}+N_{52}{\bar{S}}_{2i,j-1,\ell}^{xyt}\\&+N_{53}{\bar{S}}_{2i,j+1,\ell}^{xyt}+N_{54}{\bar{u}}_{i,j-1,\ell}^{xyt}+N_{55}{\bar{u}}_{i,j+1,\ell}^{xyt}+L_{51}{\bar{S}}_{1i,j,\ell}^{xyt}-{L_{52}\bar{S}}_{1i,j,\ell-1}^{xyt}+L_{53}{\bar{u}}_{i,j,\ell-1}^{xyt}+L_{54}{\bar{f}}_{xi,j,\ell}^{xyt} 
	\end{split}
\label{eq:084}
\end{gather}
\begin{gather}
	\begin{split}
	{\bar{S}}_{6i,j,\ell}^{xyt}&=\frac{{\bar{S}}_{4i,j,\ell}^{xyt}-{\bar{S}}_{4i,j,\ell-1}^{xyt}}{2}+\frac{{\bar{v}}_{i,j,\ell}^{xyt}-{\bar{v}}_{i,j,\ell-1}^{xyt}}{2\tau_{\ell}}+F_{51}{\bar{S}}_{5i,j,\ell}^{xyt}+F_{52}{\bar{S}}_{5i,j-1,\ell}^{xyt}+F_{53}{\bar{S}}_{5i,j+1,\ell}^{xyt}\\&+F_{54}{\bar{v}}_{i,j,\ell}^{xyt}+F_{55}{\bar{v}}_{i,j-1,\ell}^{xyt}+F_{56}{\bar{v}}_{i,j+1,\ell}^{xyt}-{\bar{f}}_{yi,j,\ell}^{xyt} 
	\end{split}
\label{eq:085}
\end{gather}
\begin{gather}
	\begin{split}
	{\bar{S}}_{5i,j,\ell}^{xyt}&=\frac{{\bar{S}}_{4i,j,\ell}^{xyt}-{\bar{S}}_{4i,j,\ell-1}^{xyt}}{2}+\frac{{\bar{v}}_{i,j,\ell}^{xyt}-{\bar{v}}_{i,j,\ell-1}^{xyt}}{2\tau_{\ell}}+F_{31}{\bar{S}}_{6i,j,\ell}^{xyt}+F_{32}{\bar{S}}_{6i-1,j,\ell}^{xyt}+F_{33}{\bar{S}}_{6i+1,j,\ell}^{xyt}\\&+F_{34}{\bar{v}}_{i,j,\ell}^{xyt}+F_{35}{\bar{v}}_{i-1,j,\ell}^{xyt}+F_{36}{\bar{v}}_{i+1,j,\ell}^{xyt}-{\bar{f}}_{yi,j,\ell}^{xyt} 
	\end{split}
\label{eq:086}
\end{gather}
\begin{gather}
	\begin{split}
{\bar{S}}_{4i,j,\ell}^{xyt}&=-\left(F_{31}{\bar{S}}_{6i,j,\ell}^{xyt}+F_{32}{\bar{S}}_{6i-1,j,\ell}^{xyt}+F_{33}{\bar{S}}_{6i+1,j,\ell}^{xyt}+F_{34}{\bar{v}}_{i,j,\ell}^{xyt}+F_{35}{\bar{v}}_{i-1,j,\ell}^{xyt}+F_{36}{\bar{v}}_{i+1,j,\ell}^{xyt}\right)\\&-\left(F_{51}{\bar{S}}_{5i,j,\ell}^{xyt}+F_{52}{\bar{S}}_{5i,j-1,\ell}^{xyt}+F_{53}{\bar{S}}_{5i,j+1,\ell}^{xyt}+F_{54}{\bar{v}}_{i,j,\ell}^{xyt}+F_{55}{\bar{v}}_{i,j-1,\ell}^{xyt}+F_{56}{\bar{v}}_{i,j+1,\ell}^{xyt}\right)+{\bar{f}}_{yi,j,\ell}^{xyt} 
	\end{split}
\label{eq:087}
\end{gather}
\begin{gather}
	\begin{split}
{\bar{v}}_{i,j,\ell}^{xyt}&=M_{51}{\bar{S}}_{6i,j,\ell}^{xyt}+M_{52}{\bar{S}}_{6i-1,j,\ell}^{xyt}+M_{53}{\bar{S}}_{6i+1,j,\ell}^{xyt}+M_{54}{\bar{v}}_{i-1,j,\ell}^{xyt}+M_{55}{\bar{v}}_{i+1,j,\ell}^{xyt}+N_{51}{\bar{S}}_{5i,j,\ell}^{xyt}+N_{52}{\bar{S}}_{5i,j-1,\ell}^{xyt}\\&+N_{53}{\bar{S}}_{5i,j+1,\ell}^{xyt}+N_{54}{\bar{v}}_{i,j-1,\ell}^{xyt}+N_{55}{\bar{v}}_{i,j+1,\ell}^{xyt}+L_{51}{\bar{S}}_{4i,j,\ell}^{xyt}-{L_{52}\bar{S}}_{4i,j,\ell-1}^{xyt}+L_{53}{\bar{v}}_{i,j,\ell-1}^{xyt}+L_{54}{\bar{f}}_{yi,j,\ell}^{xyt} 
	\end{split}
\label{eq:088}
\end{gather}
It is important to note that every coefficient in the final set of algebraic equations (\eqns\ref{eq:081} - \ref{eq:088}) is dependent upon the velocities at the current time step (${\bar{u}}_{i,j,\ell}^0$ and ${\bar{v}}_{i,j,\ell}^0$), introducing non-linearity to the entire system. Detailed definitions for all the coefficients ($A$’s, $B$’s, $F$’s, $G$’s, $L$’s, $M$’s, $N$’s) for two-dimensional case are listed in \ref{app:A}. 
\subsubsection{Approximation of non-linear convective velocity}
\label{sec:2.2.1}
The definition of the convective velocities in MNIM for the two-dimensional case is given as
\begin{gather}
	\begin{split}
		{\bar{u}}_{i,j,\ell}^0=\frac{{\bar{u}}_{i,j,\ell}^{yt}+{\bar{u}}_{i-1,j,\ell}^{yt}+{\bar{u}}_{i,j,\ell}^{xt}+{\bar{u}}_{i,j-1,\ell}^{xt}}{4}
	\end{split}
\label{eq:089}
\end{gather}
\begin{gather}
	\begin{split}
	{\bar{v}}_{i,j,\ell}^0=\frac{{\bar{v}}_{i,j,\ell}^{yt}+{\bar{v}}_{i-1,j,\ell}^{yt}+{\bar{v}}_{i,j,\ell}^{xt}+{\bar{v}}_{i,j-1,\ell}^{xt}}{4}
\end{split}
\label{eq:090}
\end{gather}
Similar to the one-dimensional case, the simple definition of ${\bar{u}}_{i,j,\ell}^0$ and ${\bar{v}}_{i,j,\ell}^0$ can be written as
\begin{gather}
	\begin{split}
		{\bar{u}}_{i,j,\ell}^0=\frac{{\bar{u}}_{i-1,j,\ell}^{xyt}+{\bar{u}}_{i,j-1,\ell}^{xyt}\add{+}{\bar{u}}_{i,j,\ell}^{xyt}+{\bar{u}}_{i+1,j,\ell}^{xyt}+{\bar{u}}_{i,j+1,\ell}^{xyt}}{5} 
	\end{split}
\label{eq:091}
\end{gather}
\begin{gather}
	\begin{split}
	{\bar{v}}_{i,j,\ell}^0=\frac{{\bar{v}}_{i-1,j,\ell}^{xyt}+{\bar{v}}_{i,j-1,\ell}^{xyt}\add{+}{\bar{v}}_{i,j,\ell}^{xyt}+{\bar{v}}_{i+1,j,\ell}^{xyt}+{\bar{v}}_{i,j+1,\ell}^{xyt}}{5}
	\end{split}
\label{eq:092}
\end{gather}
The modified definition of ${\bar{u}}_{i,j,\ell}^0$ and ${\bar{v}}_{i,j,\ell}^0$ can be derived by using \eqns\eqref{eq:071} - \eqref{eq:080}, which establish the relationship between the surface-averaged and cell-centered values, and substituting the expressions from these equations into \eqn\eqref{eq:089}, we obtain:
\begin{gather}
	\begin{split}
		{\bar{u}}_{i,j,\ell}^0 & =F_{71}{\bar{S}}_{3i-1,j,\ell}^{ytx}+F_{72}{\bar{S}}_{3i,j,\ell}^{ytx}+F_{73}{\bar{S}}_{3i+1,j,\ell}^{ytx}+F_{74}{\bar{u}}_{i-1,j,\ell}^{xyt}+F_{75}{\bar{u}}_{i,j,\ell}^{xyt}+F_{76}{\bar{u}}_{i+1,j,\ell}^{xyt}\\&+G_{71}{\bar{S}}_{2i,j-1,\ell}^{xty}+G_{72}{\bar{S}}_{2i,j,\ell}^{xty}+G_{73}{\bar{S}}_{2i,j+1,\ell}^{xty}+G_{74}{\bar{u}}_{i,j-1,\ell}^{xyt}+G_{75}{\bar{u}}_{i,j,\ell}^{xyt}+G_{76}{\bar{u}}_{i,j+1,\ell}^{xyt} 
	\end{split}
\label{eq:093}
\end{gather}
All the $F$’s and $G$’s coefficients are listed in \ref{app:A}. The other equation for ${\bar{v}}_{i,j,\ell}^0$ can be obtained in a similar manner.
\section{Results and discussion}
\label{sec:3}
This section presents results from applying the developed MCCNIM scheme to four Burgers' problems, including three one-dimensional cases and one two-dimensional \add{case}, each with a known analytical solution.
\add{The first example presents a detailed comparison of present results and errors with the classical NIM approach and other schemes, highlighting the effectiveness of the proposed MCCNIM method. 
	The next two one-dimensional examples provide detailed error analysis and examine the spatial and temporal convergence of the scheme. Additionally, the performance of MCCNIM with traditional MNIM, focusing on iteration counts and CPU runtime,  for these examples has been established.} To validate the scheme in higher dimensions, we then solve a two-dimensional problem, with the final example confirming the effectiveness of the scheme in such cases. Notably, the algebraic equations resulting from the discretization using the developed MCCNIM scheme are iteratively solved using the Picard-based method, for which the algorithm for the one-dimensional case is \add{illustrated in Algorithm \ref{Algorithm1}.}

\begin{algorithm}
\add{
\setlength{\baselineskip}{18pt} 
\caption{Picard iterative solver for one-dimensional MCCNIM} \label{Algorithm1}
\begin{algorithmic}[1]
\State \textbf{Inputs:} Total time ($T$), time step ($\Delta t$), number of grid points ($n_x$), grid spacing ($\Delta x$), Reynolds number ($Re$), convergence criterion ($\epsilon$) 
\State \textbf{Initialize:} Set $t = 0$, initialize velocity field ($\bar{u}^{xt}$)  
\While{$t < T$} \Comment{Time-stepping loop}
    \State Apply boundary conditions  
    \State \textbf{Picard Iteration:} Nonlinear fixed-point iteration \Comment{Nonlinear iteration count}  
    \Repeat
        \State Compute $(\bar{u}^0)^k$, using \eqn\eqref{eq:039} or \eqref{eq:040}, incorporating the latest available estimates of $\bar{u}^{xt}$  
        \State Calculate the coefficient matrix $A^k$ using updated values of $(\bar{u}^0)^k$  
        \State Compute right-hand-side vector $b^k$  
        \State \textbf{Solve linear system:} Using Krylov subspace method   \Comment{Linear iteration count}
        \State Solve $A^k x^{k+1} = b^k$ using GMRES (Krylov solver) 
        \State Update velocity field: $(\bar{u}^0)^{k+1} \leftarrow (\bar{u}^0)^{k}$  
    \Until{convergence criterion $\epsilon$ is met}  
    \State Update solution: $x^{n+1} \leftarrow x^{n}$  
    \State Advance time: $t + \Delta t \leftarrow t$  
\EndWhile
\end{algorithmic}
}
\end{algorithm}

\add{The proposed Picard iterative solver for MCCNIM (Algorithm \ref{Algorithm1}) is designed to solve the non-linear system resulting from the discretization of the one-dimensional governing equations (\eqns\ref{eq:035} - \ref{eq:037}). 
The solver employs a time-stepping procedure where, at each time step, a fixed-point iteration (Picard iteration) is used to handle the non-linearity. Within each Picard iteration, a linear system is formulated and solved using the Krylov (GMRES) method. Convergence is assessed based on a prescribed tolerance $\epsilon$, which is set to $1 \times 10^{-6}$ for all calculations in this study. This procedure (Algorithm \ref{Algorithm1}) has further been extended and developed  for two-dimensional problems, and utilized in solving the last test case (Example 4: two-dimensional propagating shock problem) in this work.} 
Moreover, in the context of the approximation of the averaged convective velocity (${\bar{u}}_{i,\ell}^0$), it should be noted that all calculations and results are based on the modified definition provided in \eqn\eqref{eq:040}. Instances where the initial approximation of ${\bar{u}}_{i,\ell}^0$ from \eqn\eqref{eq:039} is used for comparison and other purposes are explicitly indicated. Therefore, there should be no confusion regarding the general definition of ${\bar{u}}_{i,\ell}^0$, which is consistently defined by \eqn\eqref{eq:040}.

A distinctive feature of the MCCNIM is the averaging of discrete variables over nodes, in contrast to finite-difference methods, which perform computations at grid points. As a result, precise quantitative error analysis in MCCNIM requires the node-averaged exact solution. Although computing the node-averaged exact solution analytically is straightforward in some cases, it poses challenges in most others. 
\add{Therefore, we employ numerical integration using MATLAB's built-in functions: \texttt{integral2} for one-dimensional cases to compute $\bar{u}^{xt}$ and \texttt{integral3} for two-dimensional cases to compute $\bar{u}^{xyt}$, ensuring an accurate comparison for our numerical scheme. In the one-dimensional case of evaluating $\bar{u}^{xt}$, the exact solution is averaged by performing numerical integration using \texttt{integral2} in MATLAB over the node's width ($2a_i$) (from $x_{c}-a_i$ to $x_{c}+a_i$) and height $2\tau_\ell$ (from $t_f - 2\tau_\ell$ to $t_f$), followed by division by $4a_i\tau_\ell$. 	Similarly, in the two-dimensional case to compute $\bar{u}^{xyt}$, the exact solution is averaged using \texttt{integral3} over the node's width $2a_i$ (from $x_c-a_i$ to $x_c+a_i$), length $2b_j$ (from $y_c-b_j$ to $y_c+b_j$), and height $2\tau_\ell$ (from $t_f - 2\tau_\ell$ to $t_f$), followed by division by $8a_ib_j\tau_\ell$. Here, $x_c$ and $y_c$ denote the centroid coordinates of the corresponding node in the $x$- and $y$-directions, respectively, while $t_f$ represents the current time step.}
The root mean square (RMS) error has been computed using the following expression. 
\begin{gather}
		\mathrm{RMS}=\sqrt{\frac{1}{n_x}\sum_{i=1}^{n_x}\left|{\bar{u}}_{i,\ell}^{xt}-{\bar{u}}_{i,\ell}^{exact}\right|^2} 
		\quad\text{or}\quad
		\mathrm{RMS}=
		\sqrt{\frac{1}{n_xn_y}\sum_{i=1}^{n_x}\sum_{j=1}^{n_y}\left|{\bar{u}}_{i,j,\ell}^{xt}-{\bar{u}}_{i,j,\ell}^{exact}\right|^2} 
\label{eq:094}
\end{gather}
\add{where} $n_x$ \add{and $n_y$ denotes the number of} spatial grid points \add{in $x$- and $y$-directions}. \add{The indices $i$ and $j$ are used to denote the spatial grid points, and $\ell$ denotes the time index corresponding to the current time step ($t_f$). Notably, \eqn\eqref{eq:094} shows that the RMS error is computed at the current time step ($t_f$), ensuring consistency across all test problems throughout the paper.}

In this study, for computational setup, MATLAB R2024b has been utilized as the integrated development environment and compiler for code generation and execution. \add{The simulations are performed on an Intel Core i9-12900K processor with 16 cores (8P + 8E), 24 threads, and a base clock speed of 3.2 GHz (P-Cores) / 2.4 GHz (E-Cores), with a maximum turbo frequency of up to 5.2 GHz.}
\subsection{Numerical tests}
\label{sec:3.1}
\subsubsection{Example 1: One-Dimensional Propagating Shock Problem}
\label{sec:3.1.2}
The governing equation for the {one-dimensional} propagating shock problem \citep{23Rizwan_uddin_1997} is given by
\begin{gather}
	\begin{split}
	\frac{\partial u(x,t)}{\partial t}+u(x,t)\frac{\partial u(x,t)}{\partial x}=\frac{1}{Re}\frac{\partial^2u(x,t)}{\partial x^2}      \qquad \text{for}\quad            x\in[-2,2]
	\end{split}
\label{eq:098}
\end{gather}
subject to the following initial and boundary conditions.
\begin{gather}
	\begin{split}
		u\left(x,t=0\right)=\frac{1}{2}\left[1-\tanh\left(\frac{xRe}{4}\right)\right]
	\end{split}
\label{eq:099}
\end{gather}
%
%
%
\begin{gather}
	\begin{split}
	u\left(x =\pm 2,t\right)=\frac{1}{2}\left[1-\tanh\left(\pm \frac{Re}{2}-\frac{tRe}{8}\right)\right]
	\end{split}
\label{eq:101}
\end{gather}
The exact solution to the problem (\eqn\ref{eq:098}) is expressed as follows.
\begin{gather}
	\begin{split}
		u\left(x,t\right)=\frac{1}{2}\left[1-\tanh\left(\frac{xRe}{4}-\frac{tRe}{8}\right)\right]  
	\end{split}
\label{eq:102}
\end{gather}
\add{\fig\ref{fig:08} compares the exact solution (\eqn\ref{eq:102}) with the solution obtained using the proposed MCCNIM scheme for two cases: a low Reynolds number ($Re = 10$) in panel (a) and a high Reynolds number ($Re = 10^6$) in panel (b). For the low Re case, using a large timestep ($\Delta t = 0.1$) and a coarse grid size ($\Delta x = 0.1$), the scheme accurately captures the propagation front, as seen in panel (a) of \fig\ref{fig:08}. The robustness of the scheme has been further been established for the high $Re$ case, employing a finer timestep ($\Delta t = 0.0125$) and grid size ($\Delta x = 0.00625$), , where the scheme continues to maintain a smooth propagation front, as demonstrated in panel (b) of \fig\ref{fig:08}. These results underscore the effectiveness of the proposed scheme in handling problems with strong non-linearities.}
\begin{figure}[!b]%
	\centering
	\subfigure[$Re=10$, $\Delta x=0.1$, and $\Delta t=0.1$] {\includegraphics[width=0.48\linewidth]{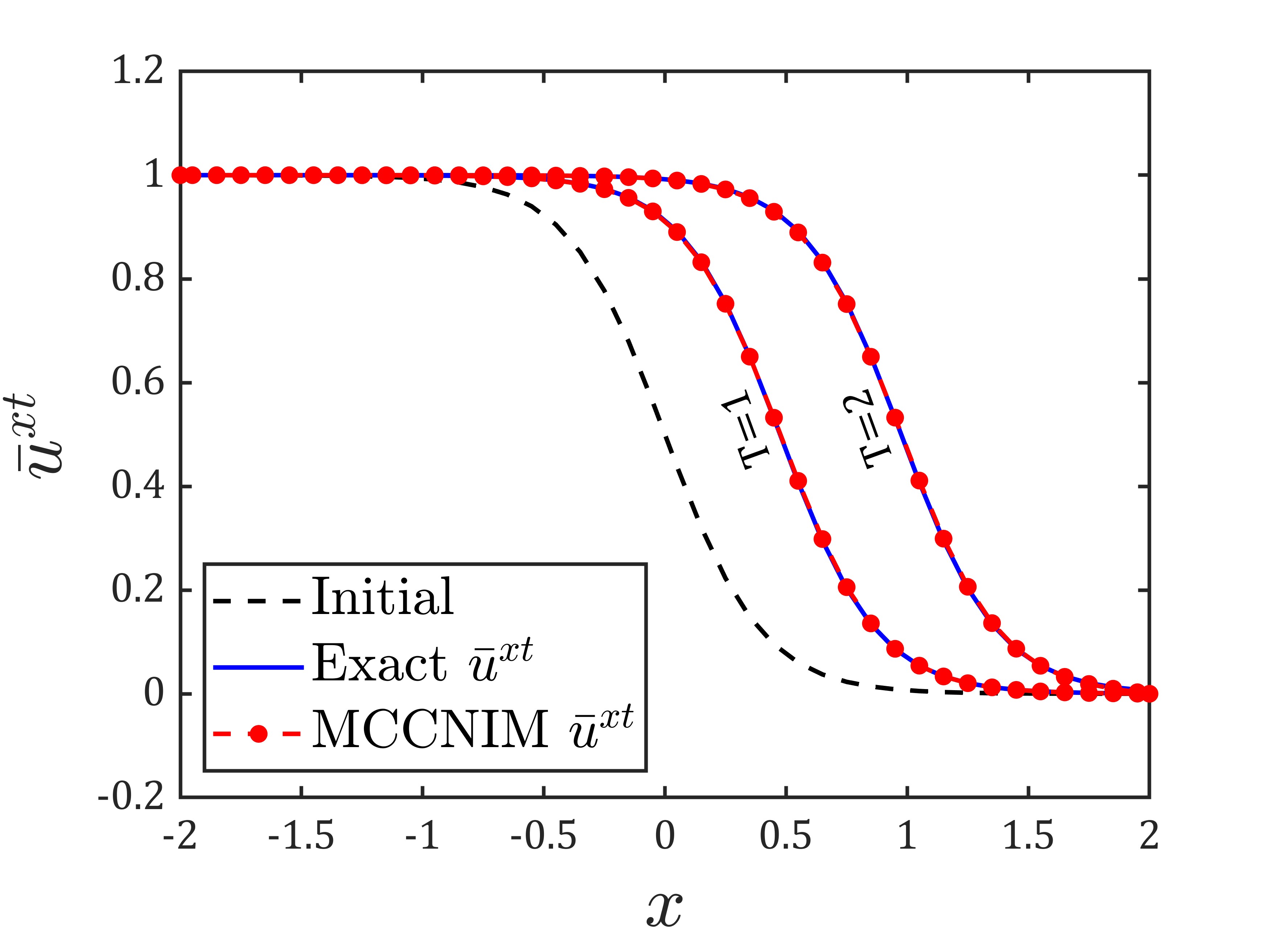}\label{fig8a}}
	\subfigure[$Re={10}^6$, $\Delta x=0.00625$, and $\Delta t=0.0125$] {\includegraphics[width=0.48\linewidth]{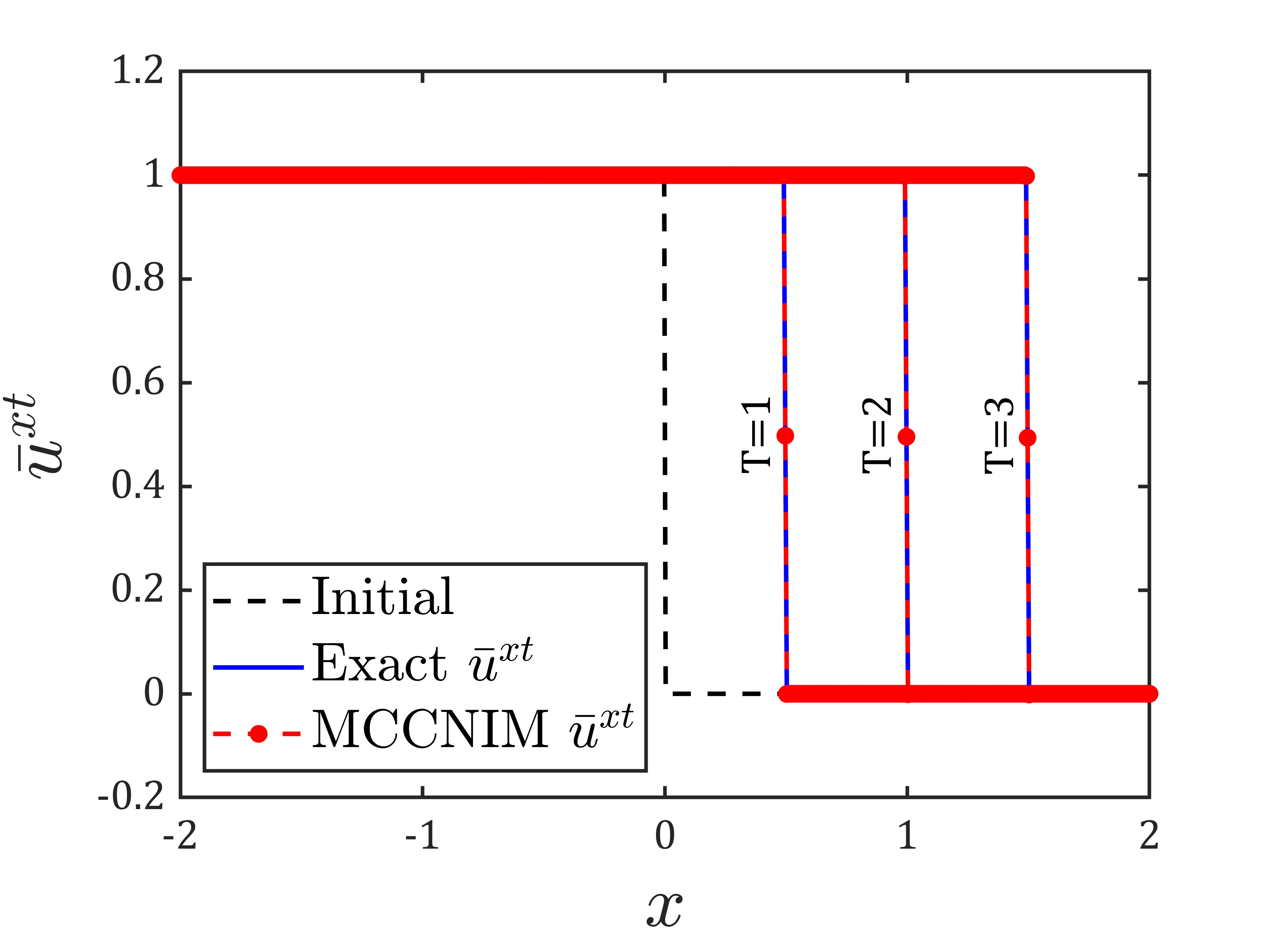}\label{fig8b}}
	\caption{Comparison of numerical and exact solutions for Example 1 at different times using MCCNIM scheme.}%
	\label{fig:08}%
\end{figure}
It is critical to understand that the analytical solution used to calculate RMS error differs for MCCNIM and MNIM. For instance, the analytical solution in MNIM \citep{23Rizwan_uddin_1997} is averaged over the surface, either spatially or temporally; on the other hand, the analytical solution in MCCNIM is averaged across the whole cell, including both temporal and spatial dimensions.
As a result, although it offers a reasonable comparison, precise accuracy may not be obtained when comparing the cell-centered value in MCCNIM with the surface-averaged value in MNIM.  

The chosen test problem has also been used in the literature \citep{23Rizwan_uddin_1997} for validation of the modified NIM (or MNIM) with another scheme called the Crank-Nicolson 4-Point Upwind (CN-4PU) scheme. When comparing the RMS errors of the two schemes \citep{23Rizwan_uddin_1997}, it was shown that MNIM performs noticeably better than CN-4PU. \add{To ensure a reliable and meaningful performance of the presently developed MCCNIM scheme, in comparison to the MNIM scheme \citep{23Rizwan_uddin_1997}, this study has compared the numerical and exact solutions for three representative cases, covering all relevant scenarios as follows.}
\begin{itemize}
    \item \add{\textbf{Case 1}: Solution at $T = 1$ and $T = 3$ with $Re = 50$, $n_x = 20$, and $\Delta t = 0.01$, based on the reference case \citep[pages 130--131,][]{23Rizwan_uddin_1997}.}
    \item \add{\textbf{Case 2}: Solution at $T = 1$ and $T = 3$ with $Re = 100$, $n_x = 20$, and $\Delta t = 0.01$, based on the reference case \citep[pages 132--133,][]{23Rizwan_uddin_1997}.}
    \item \add{\textbf{Case 3}: Solution at $T = 1.5$ with $Re = 300$, $n_x = 100$, and $\Delta t = 0.05$, based on the reference case \citep[pages 128,][]{23Rizwan_uddin_1997}.}
\end{itemize}
\begin{figure}[!b]%
	\centering
	\subfigure[$\bar{u}^x$ vs. $x$ using MNIM] {\includegraphics[width=0.48\linewidth]{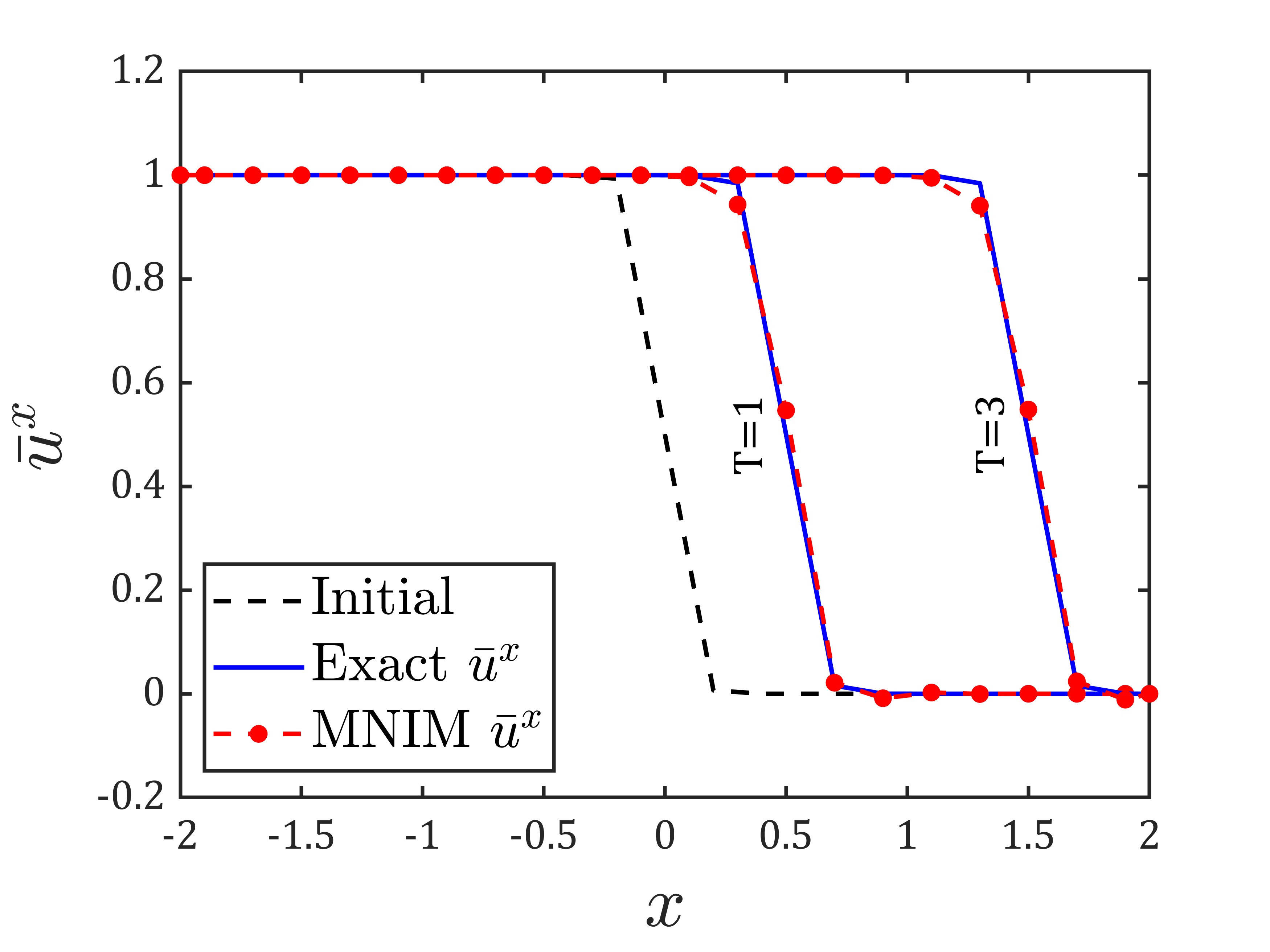}\label{fig13a}}
	\subfigure[$\bar{u}^t$ vs. $x$ using MNIM] {\includegraphics[width=0.48\linewidth]{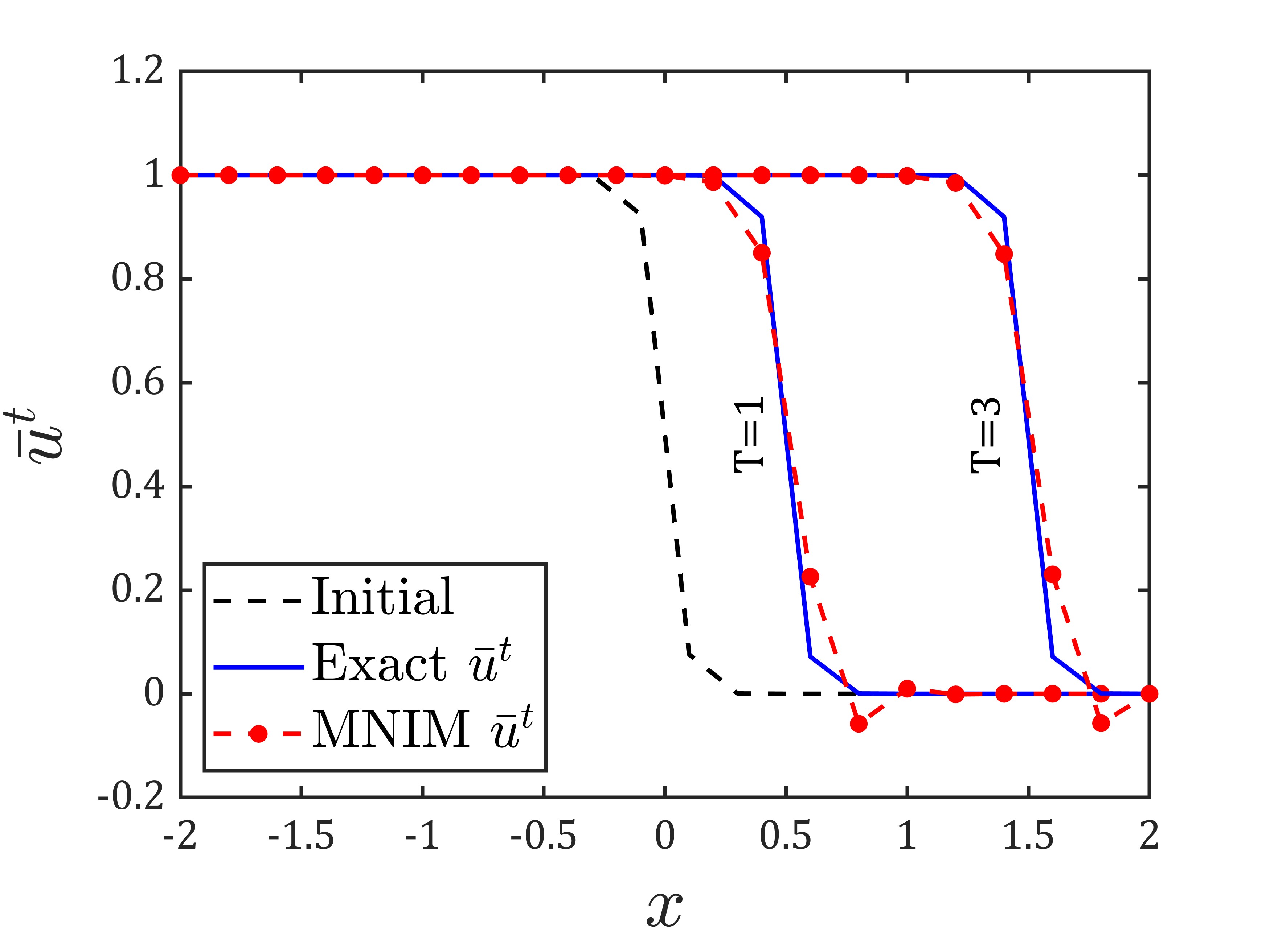}\label{fig13b}}\\%
	\subfigure[$\bar{u}^{xt}$ vs. $x$ using MCCNIM] {\includegraphics[width=0.48\linewidth]{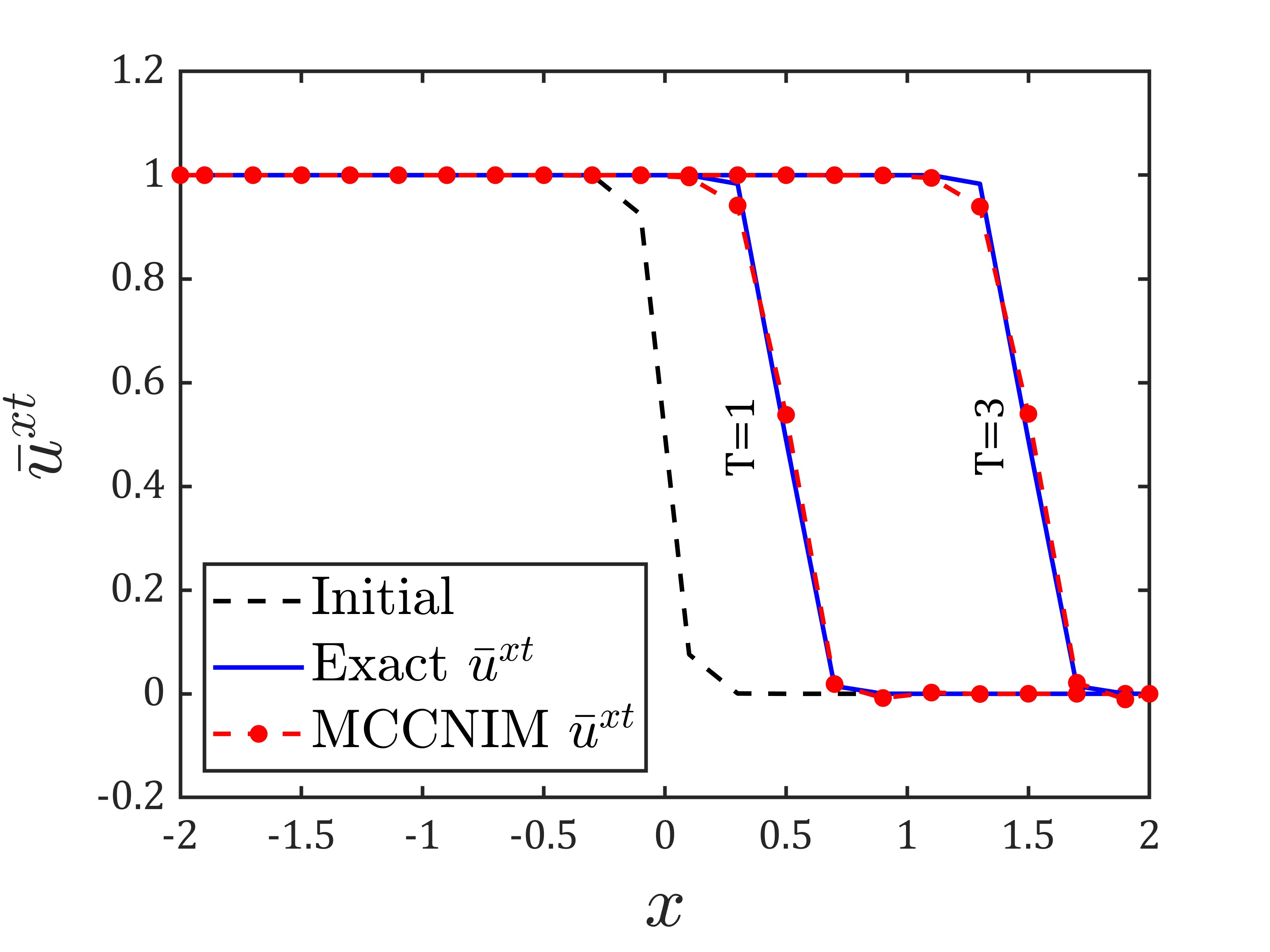}\label{fig13c}}%
	\caption{\add{Comparions of the results for \textbf{Case 1}. The RMS errors are obtained as $1.406 \times 10^{-2}$ at ($T=1$) and $1.484 \times 10^{-2}$ at ($T=3$) for MNIM using $\bar{u}^x$; $4.119 \times 10^{-2}$ at ($T=1$) and $4.217 \times 10^{-2}$ at ($T=3$) for MNIM using $\bar{u}^t$; $3.047 \times 10^{-2}$ at ($T=1$) and $3.129 \times 10^{-2}$ at ($T=3$) for MNIM using both $\bar{u}^x$ and $\bar{u}^t$; $1.383 \times 10^{-2}$ at ($T=1$) and $1.456 \times 10^{-2}$ at ($T=3$) for MCCNIM using $\bar{u}^{xt}$.}}%
	\label{fig:Case1}%
\end{figure}
\begin{figure}[!b]%
	\centering
	\subfigure[$\bar{u}^x$ vs. $x$ using MNIM] {\includegraphics[width=0.48\linewidth]{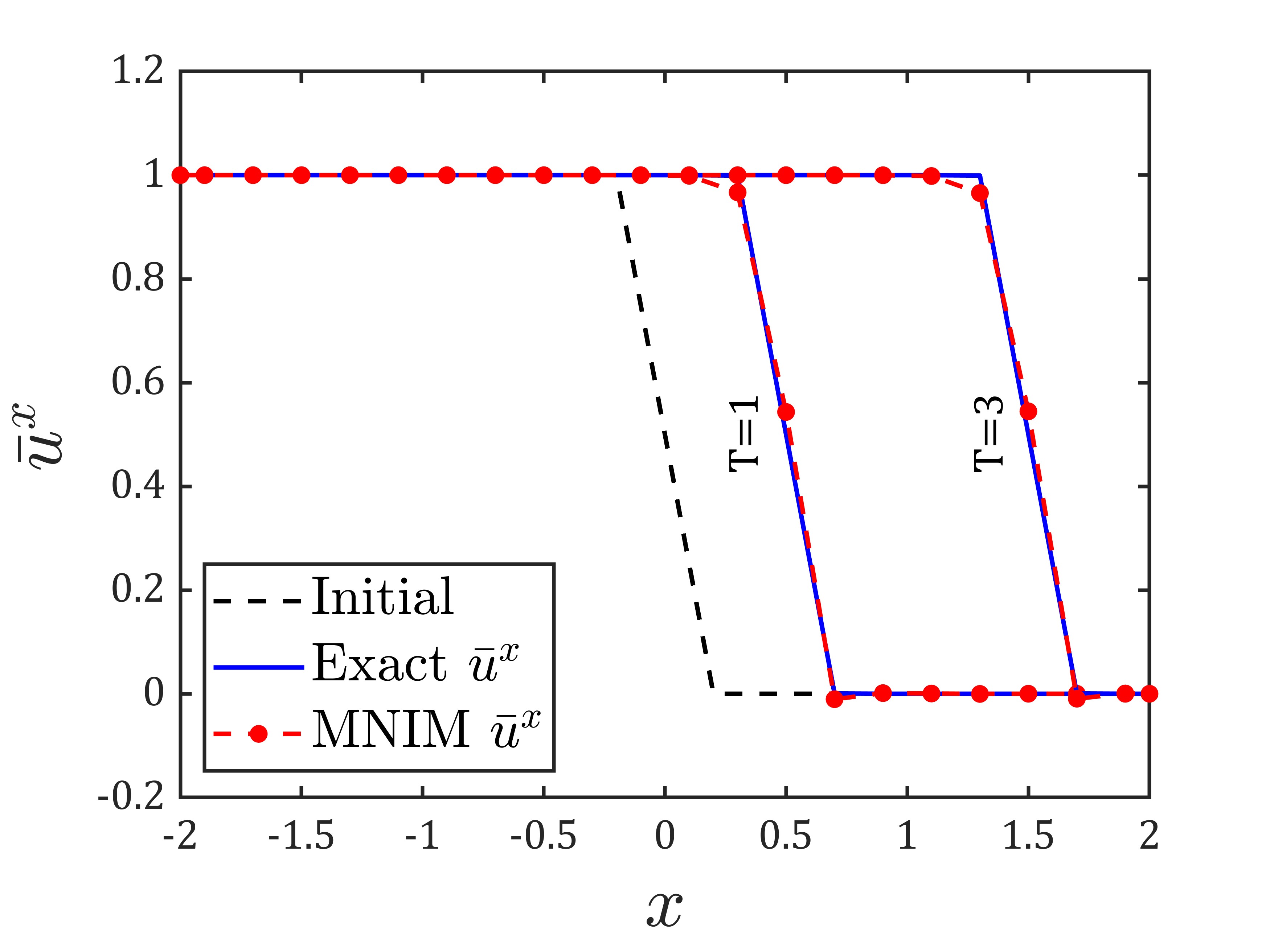}\label{fig13a}}
	\subfigure[$\bar{u}^t$ vs. $x$ using MNIM] {\includegraphics[width=0.48\linewidth]{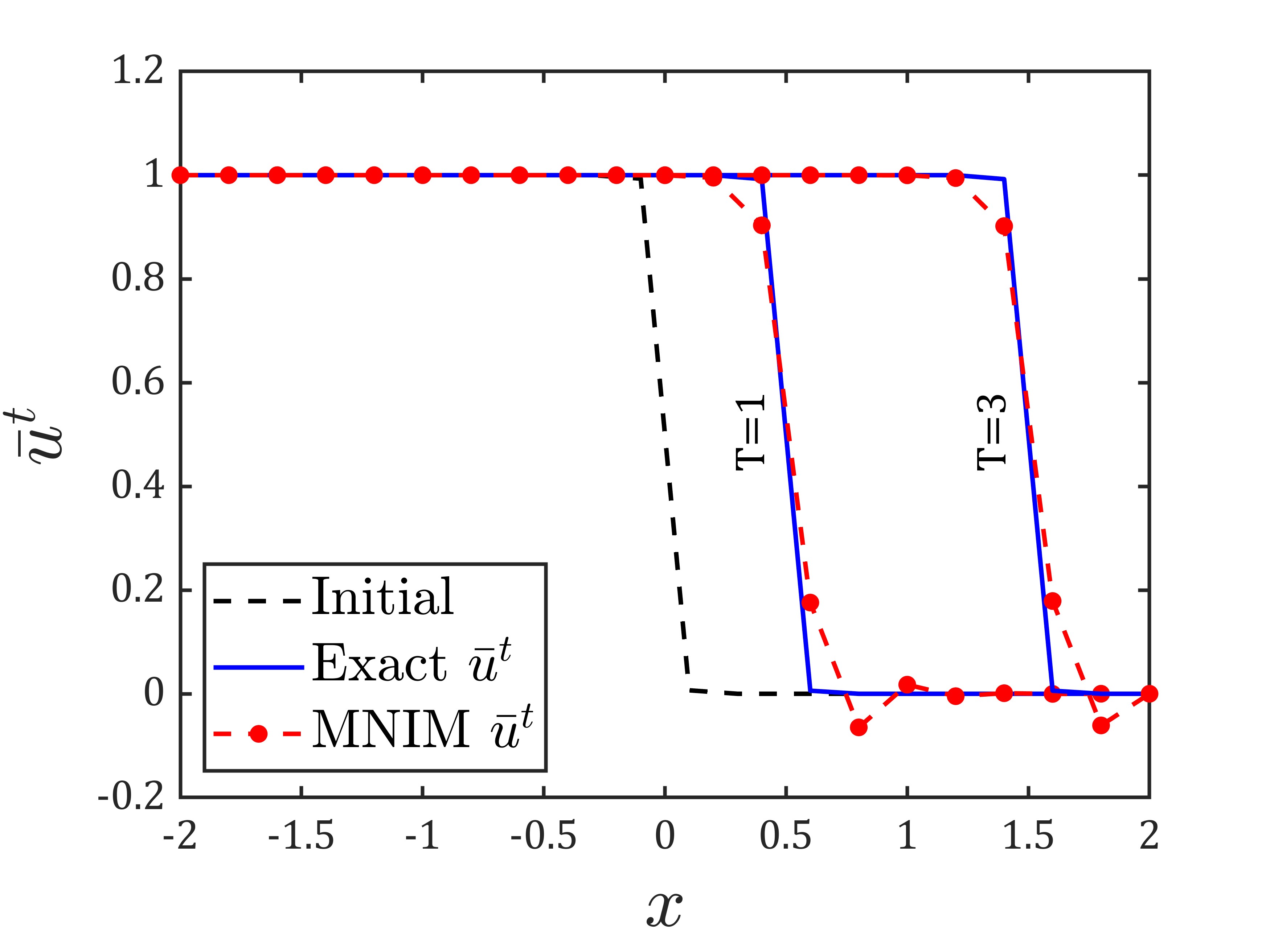}\label{fig13b}}\\%
	\subfigure[$\bar{u}^{xt}$ vs. $x$ using MCCNIM] {\includegraphics[width=0.48\linewidth]{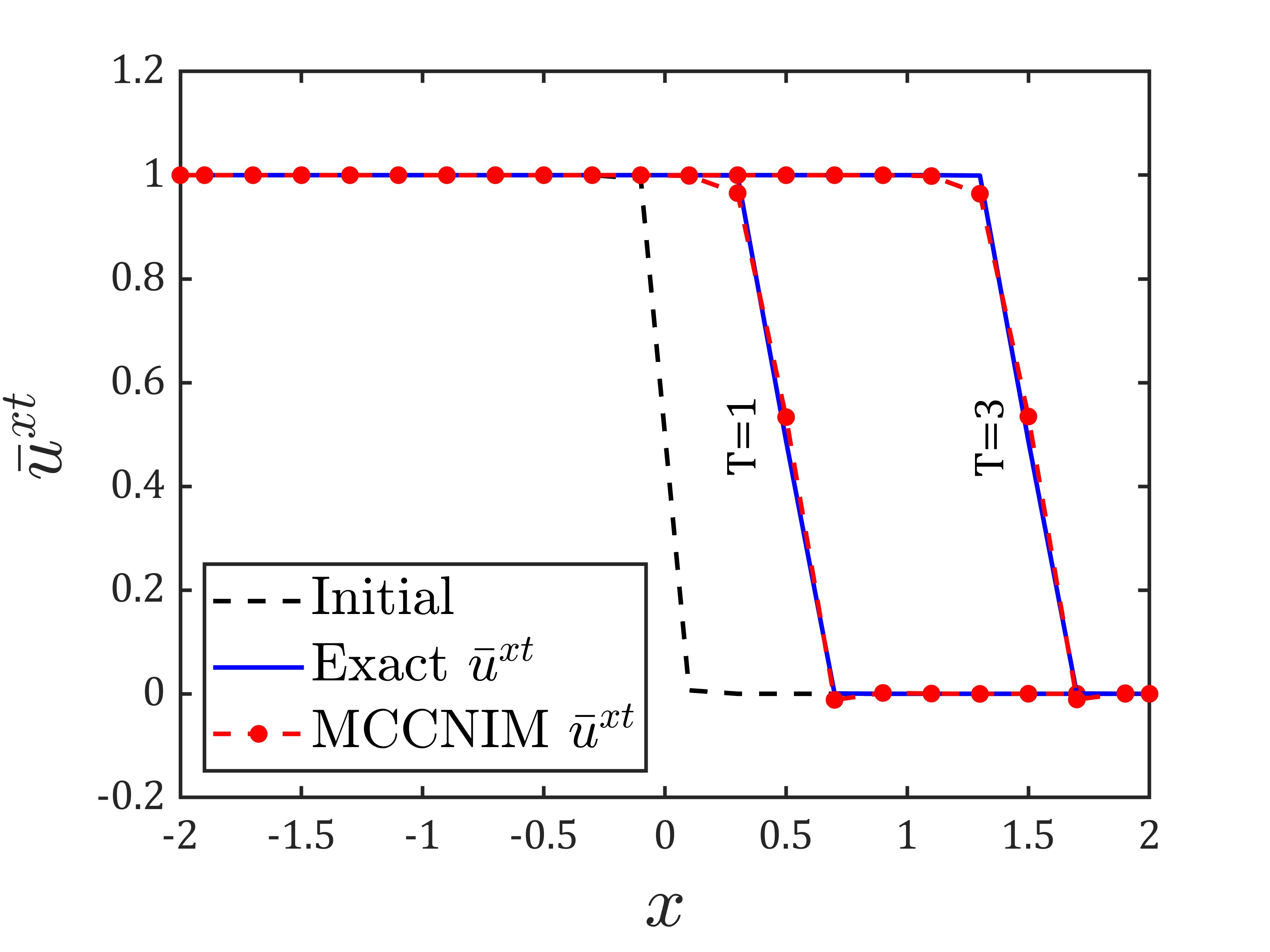}\label{fig13c}}%
	\caption{\add{Comparions of the results for \textbf{Case 2}. The RMS errors are obtained as $1.239 \times 10^{-2}$ at ($T=1$) and $1.273 \times 10^{-2}$ at ($T=3$) for MNIM using $\bar{u}^x$; $4.666 \times 10^{-2}$ at ($T=1$) and $4.694 \times 10^{-2}$ at ($T=3$) for MNIM using $\bar{u}^t$; $3.375 \times 10^{-2}$ at ($T=1$) and $3.401 \times 10^{-2}$ at ($T=3$) for MNIM using both $\bar{u}^x$ and $\bar{u}^t$; $1.241 \times 10^{-2}$ at ($T=1$) and $1.275 \times 10^{-2}$ at ($T=3$) for MCCNIM using $\bar{u}^{xt}$.}}%
	\label{fig:Case2}%
\end{figure}
\begin{figure}[!b]%
	\centering
	\subfigure[$\bar{u}^x$ vs. $x$ using MNIM] {\includegraphics[width=0.48\linewidth]{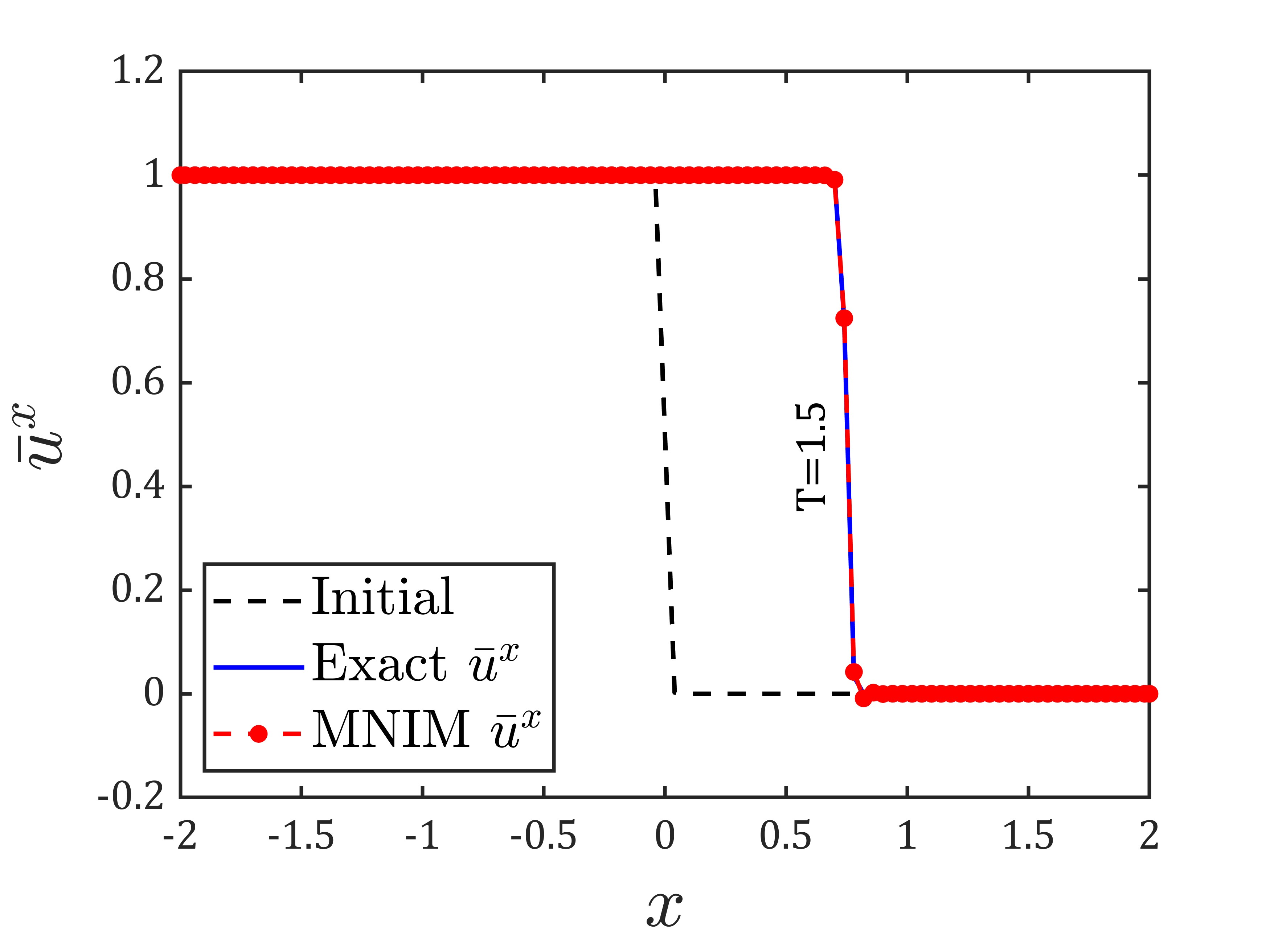}\label{fig13a}}
	\subfigure[$\bar{u}^t$ vs. $x$ using MNIM] {\includegraphics[width=0.48\linewidth]{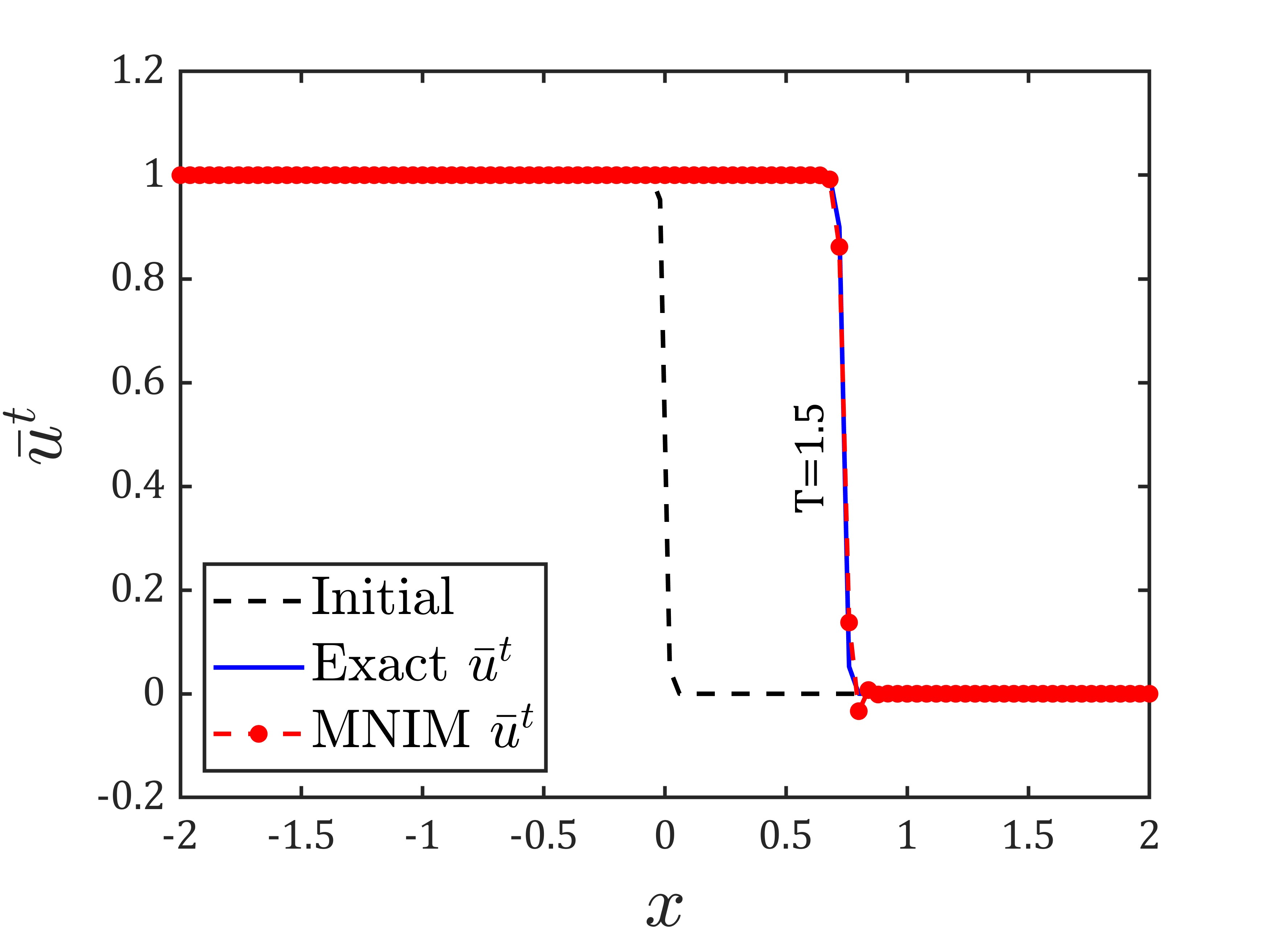}\label{fig13b}}\\%
	\subfigure[$\bar{u}^{xt}$ vs. $x$ using MCCNIM] {\includegraphics[width=0.48\linewidth]{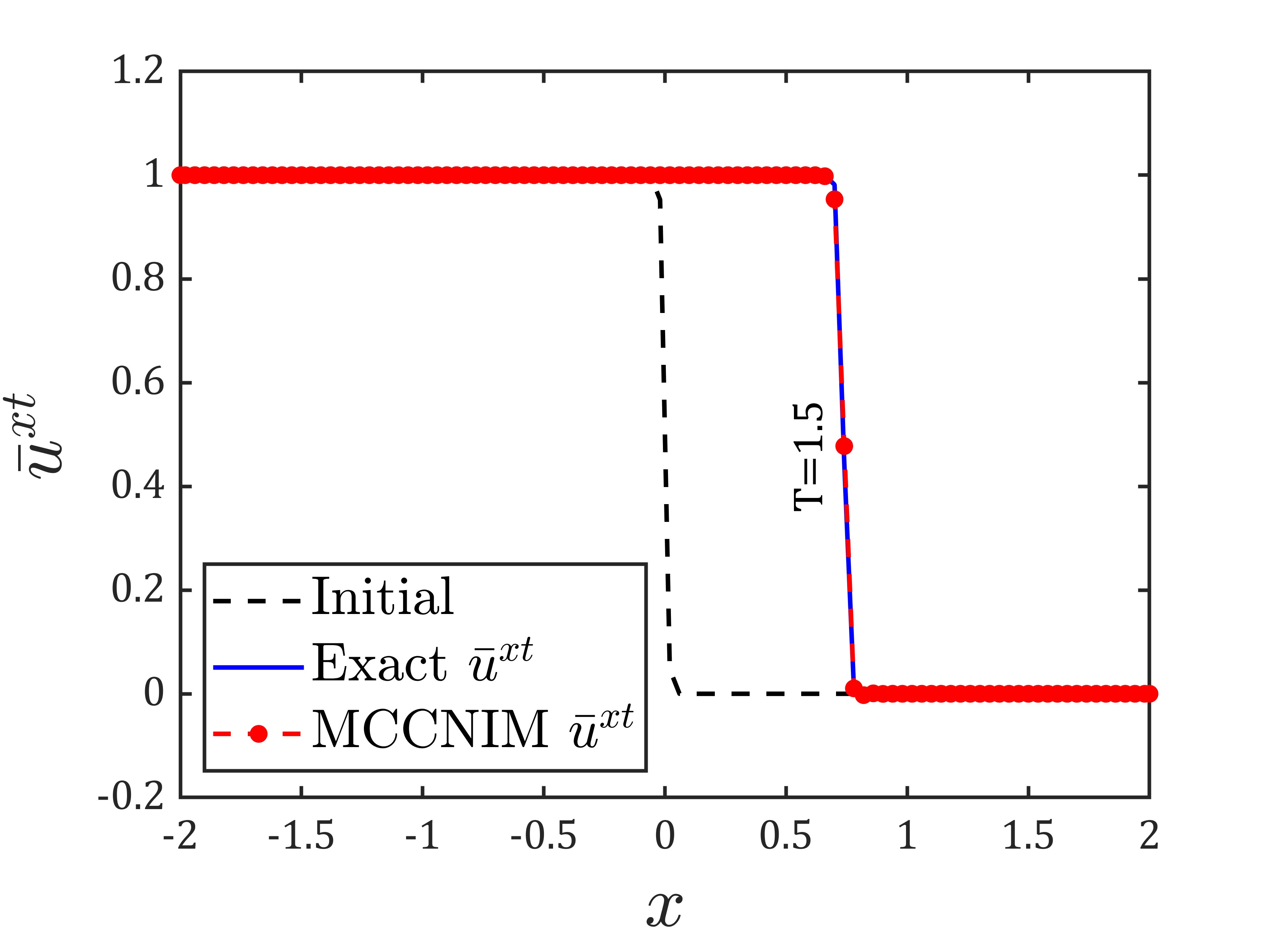}\label{fig13c}}%
	\caption{\add{Comparions of the results for \textbf{Case 3}. The RMS errors are obtained as $1.562 \times 10^{-3}$ for MNIM using $\bar{u}^x$; $1.002 \times 10^{-2}$ for MNIM using $\bar{u}^t$; $7.155 \times 10^{-3}$ MNIM using both $\bar{u}^x$ and $\bar{u}^t$; $4.148 \times 10^{-3}$ for MCCNIM using $\bar{u}^{xt}$.}}%
	\label{fig:Case3}%
\end{figure}
\add{A direct comparison has been ensured by independently developed in-house source code for the MNIM scheme by strictly following the methodology described by \citet{23Rizwan_uddin_1997}. \figs\ref{fig:Case1} to \ref{fig:Case3} illustrate the comparisons presenting the exact solution alongside the results from both the MNIM and MCCNIM schemes for Cases 1 to 3, respectively. The results are validated by computing the RMS errors for all the cases, showing exact agreement with the literature \citep{23Rizwan_uddin_1997}.} 
Due to distinct errors associated with the variables, \figs\ref{fig:Case1} to \ref{fig:Case3} present two separate panels for MNIM; Panel (a) of each figure illustrates the comparison of the solution of ${\bar{u}}^x$ with the exact ${\bar{u}}^x$, while panel (b) depicts the solution of ${\bar{u}}^t$ compared to the exact ${\bar{u}}^t$; in the panel (c), the solution of the developed scheme (MCCNIM) for ${\bar{u}}^{xt}$ is contrasted with the exact ${\bar{u}}^{xt}$ solution. \add{As evident in each figure, the MCCNIM scheme, shown in panel (c), provides the closest match to the exact solution, demonstrating its superior accuracy compared to the traditional MNIM scheme.}
For example,  the error discrepancy between ${\bar{u}}^x$ (space-averaged surface), ${\bar{u}}^t$ (time-averaged surface), and the average of  ${\bar{u}}^x$ and ${\bar{u}}^t$ values are substantially smaller in the Case 3 of $Re=300$. It would be preferable to use the averaged error of both surface-averaged values (i.e., the average of both ${\bar{u}}^x$ and ${\bar{u}}^t$) when comparing MCCNIM with NIM \citep{22Rizwan_uddin_1997}. 
\add{The RMS errors computed for case 3  (shown in \fig\ref{fig:Case3}) are obtained as $0.1562 \times 10^{-2}$,  $0.1002 \times 10^{-1}$ and $0.7155 \times 10^{-2}$ using spatial-averaged velocity ($\bar{u}^x$),  time-averaged velocity ($\bar{u}^t$) and combined space-time averaged velocity ($\bar{u}^x$, $\bar{u}^t$), respectively,  which perfectly match the literature results \citep{23Rizwan_uddin_1997}, thereby confirming the accuracy and validity of the MNIM implementation.}
%
%
%
In contrast, the RMS error in MCCNIM for the identical case is $0.4148\times{10}^{-2}$, which shows a decrease. 

\add{Furthermore, one additional case ($Re = 10^3$) has been illustrated in \fig\ref{fig:09} to compare the exact and numerical solutions for both schemes to establish the performance of MCCNIM scheme at higher Reynolds numbers.} 
%
%
\begin{figure}[!b]%
	\centering
	\subfigure[$\bar{u}^x$ vs. $x$ using MNIM] {\includegraphics[width=0.48\linewidth]{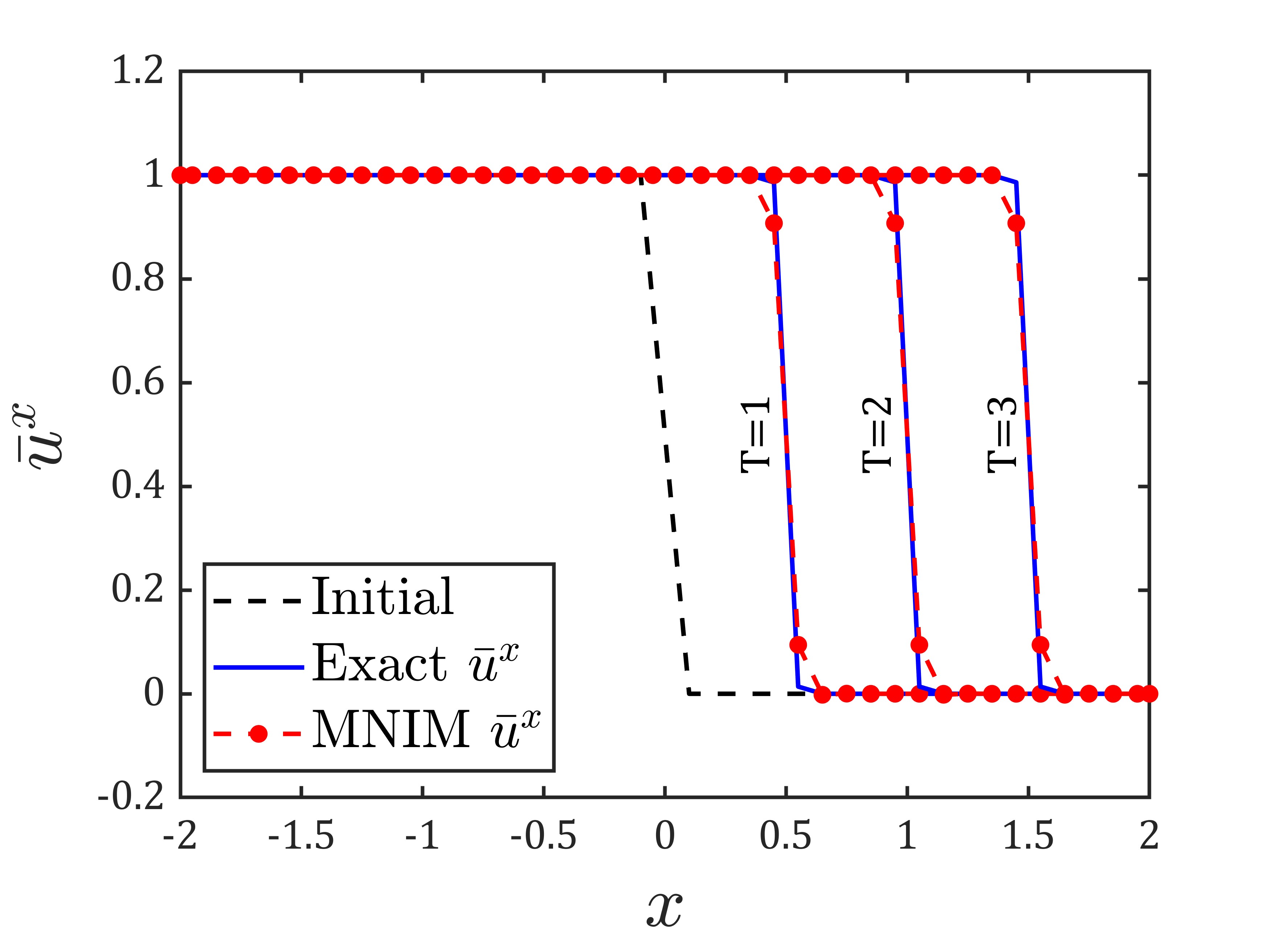}\label{fig9a}}
	\subfigure[$\bar{u}^t$ vs. $x$ using MNIM] {\includegraphics[width=0.48\linewidth]{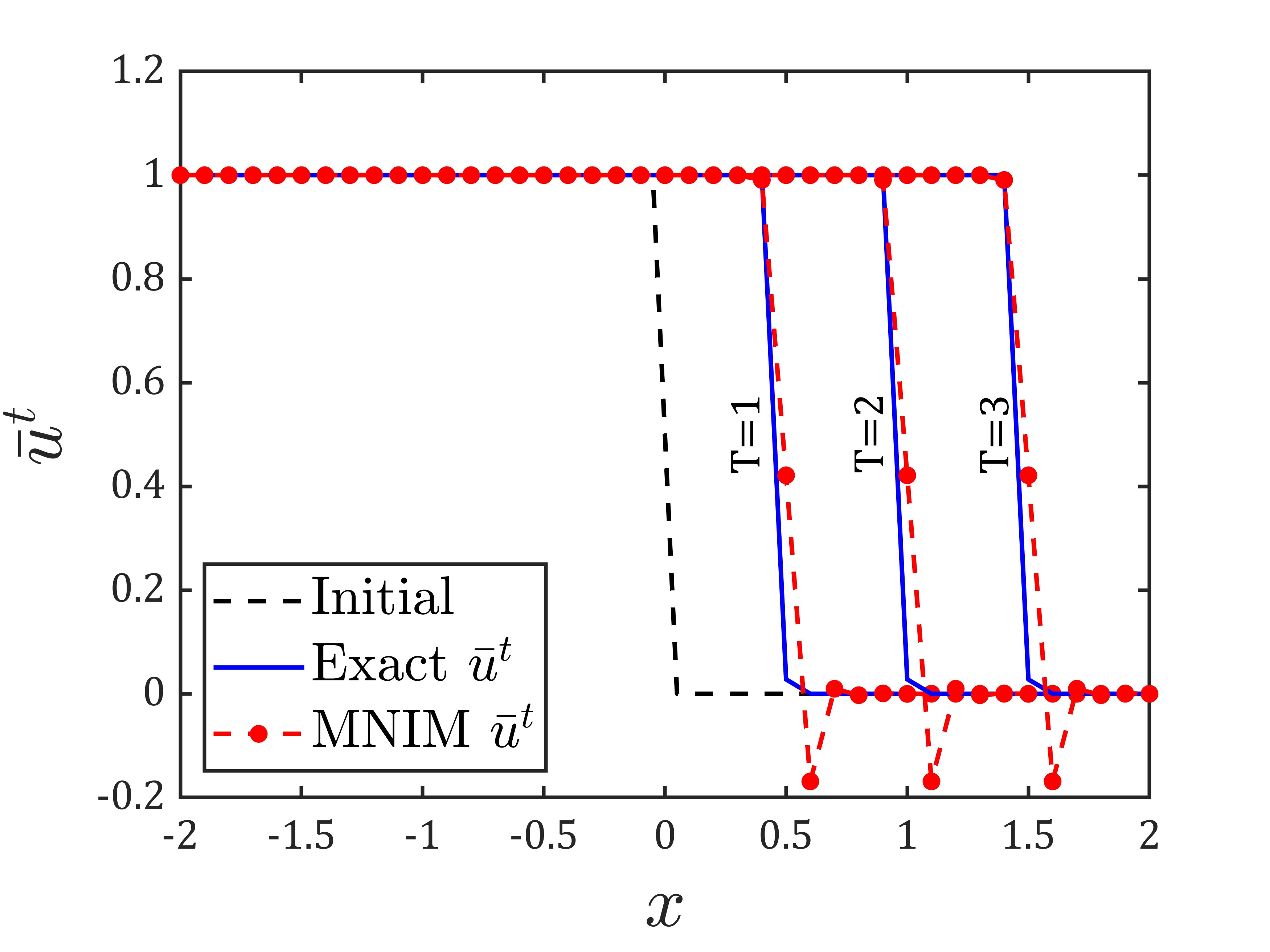}\label{fig9b}}\\
	\subfigure[$\bar{u}^{xt}$ vs. $x$ using MCCNIM] {\includegraphics[width=0.48\linewidth]{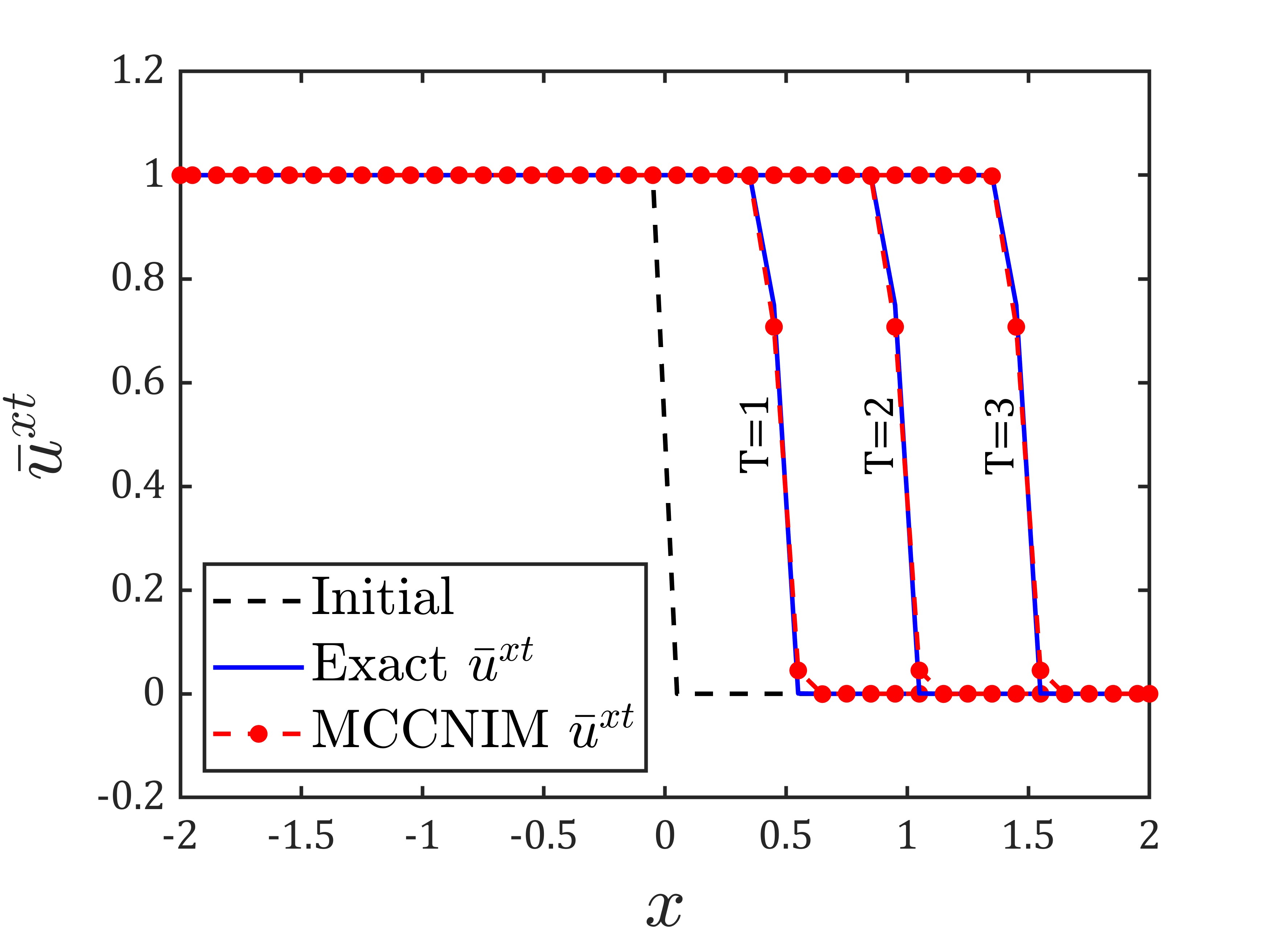}\label{fig9c}}
	\caption{Comparison of numerical and analytical solutions of Example 1 using MNIM and MCCNIM at \(T = 1\), \(2\), and \(3\) with \(Re = 10^3\), \(n_x = 40\), and \(\Delta t = 0.1\). The RMS errors at $T=1$ are obtained as $1.781\times{10}^{-3}$ for MNIM using $\bar{u}^x$; $6.862\times{10}^{-2}$ for MNIM using $\bar{u}^t$; $4.985\times{10}^{-2}$ MNIM using both $\bar{u}^x$ and $\bar{u}^t$; $9.415\times{10}^{-3}$ for MCCNIM using $\bar{u}^{xt}$.}%
	\label{fig:09}%
\end{figure}
%
%
%
As depicted in panel (c) of \fig\ref{fig:09}, the MCCNIM demonstrates highly accurate tracking of the wave front even at high Reynolds numbers, with an RMS error of $9.415\times{10}^{-3}$. This precision is achieved with minimal oscillation and without the necessity for artificial dissipation. In contrast, for MNIM, the comparison of ${\bar{u}}^x$ with the exact ${\bar{u}}^x$ is notably accurate, featuring a low RMS error of $1.781\times{10}^{-3}$. However, when assessing ${\bar{u}}^t$ against the exact ${\bar{u}}^t$, a higher RMS error is observed, measuring $6.862\times{10}^{-2}$, as illustrated in panel (b) of \fig\ref{fig:09}. 
The combined RMS error, incorporating both ${\bar{u}}^x$ and ${\bar{u}}^t$, is $4.985\times{10}^{-2}$, which exceeds the overall error in MNIM. Additional information regarding the quantitative comparison of RMS errors between MCCNIM and MNIM, along with errors from various other schemes  \citep{23Rizwan_uddin_1997}, is presented in \tab\ref{tab:1}.
\begin{table}[!t]
	\begin{center}\renewcommand{\arraystretch}{1.5}
		\caption{Comparison of RMS errors of MCCNIM (using the two definitions of ${\bar{u}}_{i,\ell}^0$ given by \eqns\ref{eq:039} and \ref{eq:040}) vs. CN-4PU, CNIM, and MNIM \citep{23Rizwan_uddin_1997} for $\Delta t=0.1$, {and} 20 nodes.}\label{tab:1}
		\resizebox{0.6\textwidth}{!}{		
			\begin{tabular}{|l|c|c|c|c|}
				\hline
		Method & \multicolumn{4}{c|}{RMS errors ($\times 10^{-2}$)} \\ \cline{2-5}
		 & \multicolumn{2}{c|}{$Re = 50$} & \multicolumn{2}{c|}{$Re = 100$} \\ \cline{2-5}
		& $T = 1$ & $T = 3$ & $T = 1$ & $T = 3$ \\ \hline
		CN-4PU & 3.652 & 3.524 & 4.546 & 4.367 \\ \hline
		CNIM   & 3.091 & 3.068 & 6.558 & 6.376 \\ \hline
		MNIM ($\bar{u}^x$) & 1.338 & 1.382 & 1.172 & 1.185 \\ \hline
		MNIM ($\bar{u}^t$) & 2.721 & 2.775 & 2.666 & 2.678 \\ \hline
		MNIM ($\bar{u}^x$ and $\bar{u}^t$) & 2.127 & 2.175 & 2.041 & 2.053 \\ \hline
		MCCNIM using \eqn\eqref{eq:039} ($\bar{u}^{xt}$) & 1.933 & 2.170 & 2.059 & 2.184 \\ \hline
		MCCNIM using \eqn\eqref{eq:040} ($\bar{u}^{xt}$) & 1.649 & 1.722 & 1.648 & 1.688 \\ \hline
		\end{tabular}}
	\end{center}
\end{table}

It is important to note that we have derived two definitions for the node-averaged velocities (${\bar{u}}_{i,\ell}^0$) and compared their associated errors using the formulations provided in \eqns\eqref{eq:039} and \eqref{eq:040}. As shown in \tab\ref{tab:1}, the MCCNIM scheme employing the simple definition of ${\bar{u}}_{i,\ell}^0$ from \eqn\eqref{eq:039} demonstrates accuracy comparable to the traditional MNIM. However, when ${\bar{u}}_{i,\ell}^0$ is defined according to \eqn\eqref{eq:040}, the scheme achieves better accuracy in all cases.
\begin{table}[!t]
	\begin{center}\renewcommand{\arraystretch}{1.5}
		\caption{Comparison of RMS errors of developed MCCNIM (using the two definitions of ${\bar{u}}_{i,\ell}^0$ given by \eqns\ref{eq:039} and \ref{eq:040}) with MNIM and M{$^2$}NIM \citep{25Wescott_2001} for $Re=100$ and $T=2$.}\label{tab:2}
		\resizebox{0.7\textwidth}{!}{		
			\begin{tabular}{|c|c|c|c|c|c|}
				\hline
		$\Delta x$ & $\Delta t$ & \multicolumn{4}{c|}{RMS errors ($\times 10^{-2}$)} \\ \cline{3-6}
		&  & MNIM & M{$^2$}NIM & MCCNIM (\eqn\ref{eq:039}) & MCCNIM (\eqn\ref{eq:040}) \\ \hline		
		0.2 & 0.02  & 4.88 & 4.75 & 3.74 & 2.61 \\ \cline{2-6}
		& 0.01  & 4.61 & 4.56 & 3.78 & 2.67 \\ \cline{2-6}
		& 0.005 & 4.48 & 4.45 & 3.80 & 2.69 \\ \hline
		0.1 & 0.02  & 2.15 & 1.94 & 2.33 & 1.44 \\ \cline{2-6}
		& 0.01  & 2.02 & 1.95 & 2.35 & 1.49 \\ \cline{2-6}
		& 0.005 & 1.95 & 1.92 & 2.36 & 1.52 \\ \hline
		0.05 & 0.02  & 0.759 & 0.601 & 0.989 & 0.615 \\ \cline{2-6}
		& 0.01  & 0.757 & 0.679 & 0.985 & 0.644 \\ \cline{2-6}
		& 0.005 & 0.752 & 0.730 & 0.992 & 0.653 \\ \hline
		\end{tabular}}
	\end{center}
\end{table}
\begin{table}[!t]
	\begin{center}\renewcommand{\arraystretch}{1.5}
		\caption{Comparison of RMS errors of developed MCCNIM (using the two definitions of ${\bar{u}}_{i,\ell}^0$ given by \eqns\ref{eq:039} and \ref{eq:040}) with MNIM and M{$^2$}NIM \citep{25Wescott_2001} for $Re=200$ and $T=2$.}\label{tab:3}
		\resizebox{0.7\textwidth}{!}{		
			\begin{tabular}{|c|c|c|c|c|c|}
				\hline
		$\Delta x$ & $\Delta t$ & \multicolumn{4}{c|}{RMS errors ($\times 10^{-2}$)} \\ \cline{3-6}
		 &  & MNIM & M{$^2$}NIM & MCCNIM (\eqn\ref{eq:039}) & MCCNIM (\eqn\ref{eq:040}) \\ \hline		
		0.2 & 0.02 & 6.36 & 6.59 & 4.39 & 3.01 \\ \cline{2-6}
		& 0.01 & 6.20 & 6.12 & 4.47 & 3.08 \\ \cline{2-6}
		& 0.005 & 5.90 & 5.87 & 4.51 & 3.11 \\ \hline
		0.1 & 0.02 & 3.77 & 3.45 & 2.99 & 1.79 \\ \cline{2-6}
		& 0.01 & 3.43 & 3.33 & 3.03 & 1.89 \\ \cline{2-6}
		& 0.005 & 3.24 & 3.20 & 3.04 & 1.93 \\ \hline
		0.05 & 0.02 & 1.66 & 1.29 & 2.49 & 0.92 \\ \cline{2-6}
		& 0.01 & 1.51 & 1.33 & 2.50 & 1.03 \\ \cline{2-6}
		& 0.005 & 1.42 & 1.37 & 2.51 & 1.07 \\ \hline
		\end{tabular}}
	\end{center}
\end{table}
A comprehensive analysis of RMS errors has also been presented in \tabs\ref{tab:2} and \ref{tab:3}, drawing comparisons between MCCNIM (using the two definitions of ${\bar{u}}_{i,\ell}^0$ given by \eqns\ref{eq:039} and \ref{eq:040}), MNIM, and M$^2$NIM across a spectrum of scenarios involving varying Reynolds numbers, node sizes, and time steps \citep{25Wescott_2001}. The results unequivocally demonstrate that, even in the case of the non-linear Burgers' equation, MCCNIM consistently outperforms other traditional NIMs developed thus far for fluid flow applications.
\begin{figure}[!b]%
	\centering
	\subfigure[$n_x = 80$]{\includegraphics[width=0.48\linewidth]{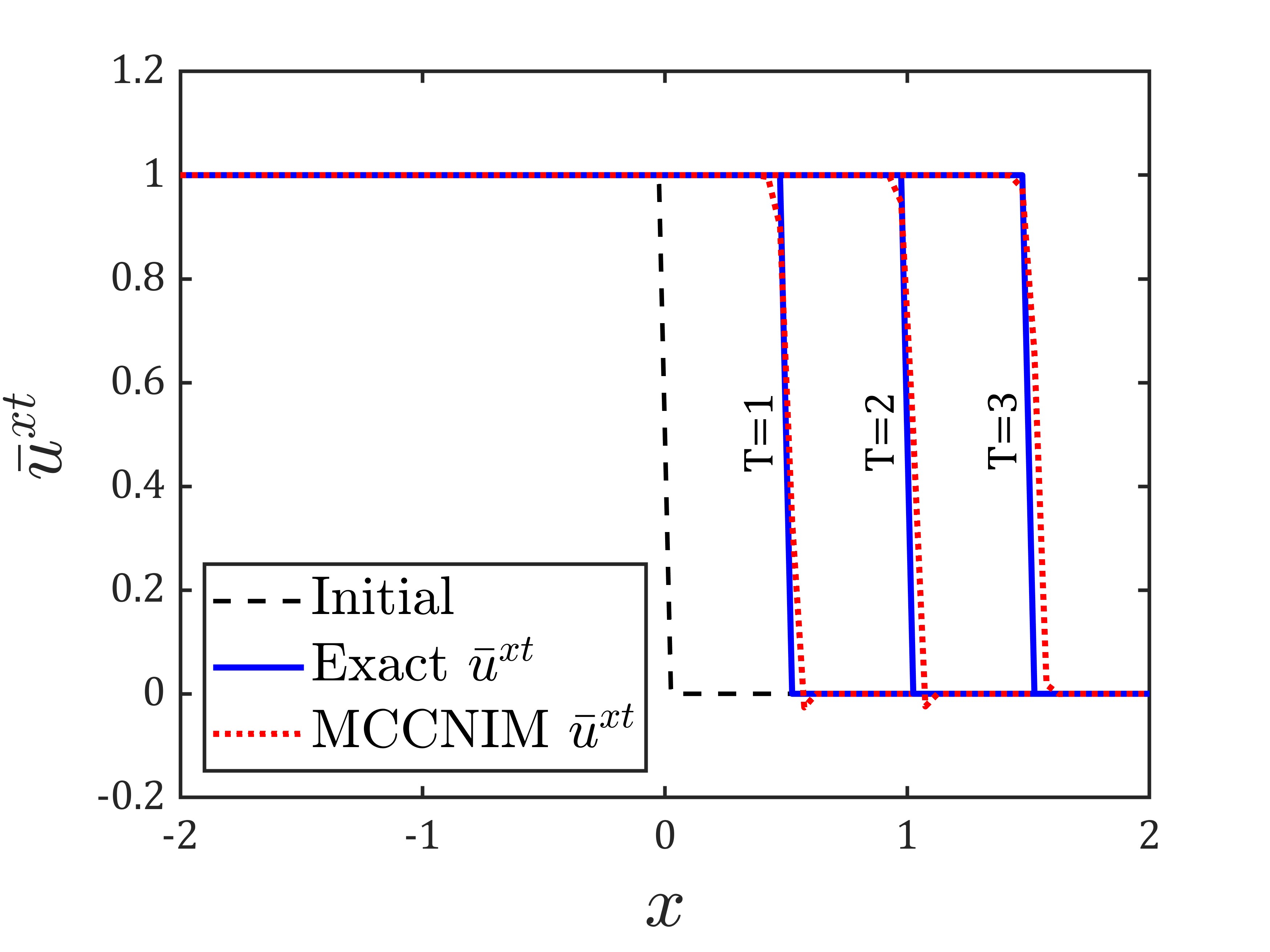}\label{fig_nx_80}}
	\subfigure[$n_x = 160$] {\includegraphics[width=0.48\linewidth]{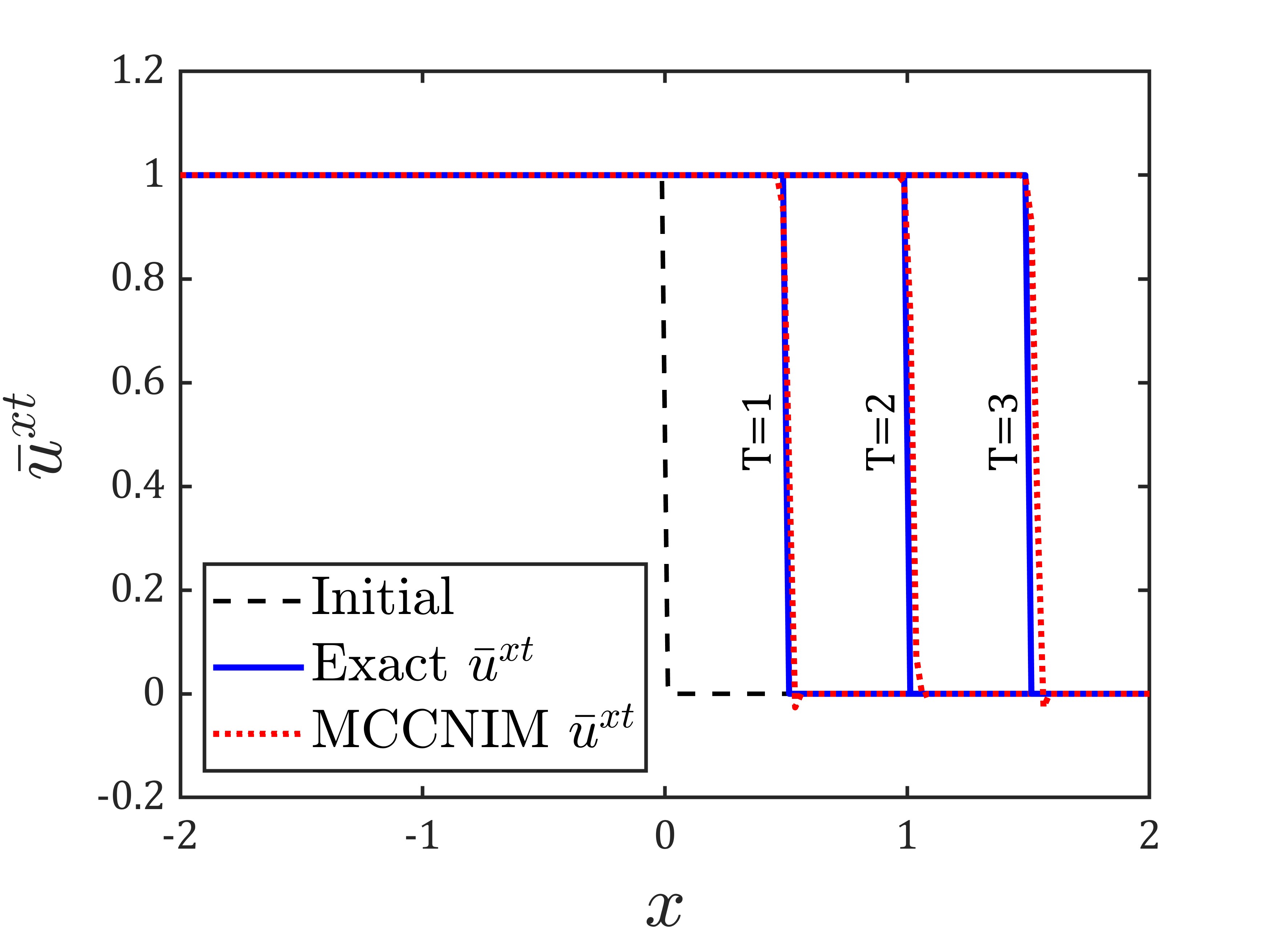}\label{fig_nx_160}}\\
	\subfigure[$n_x = 320$] {\includegraphics[width=0.48\linewidth]{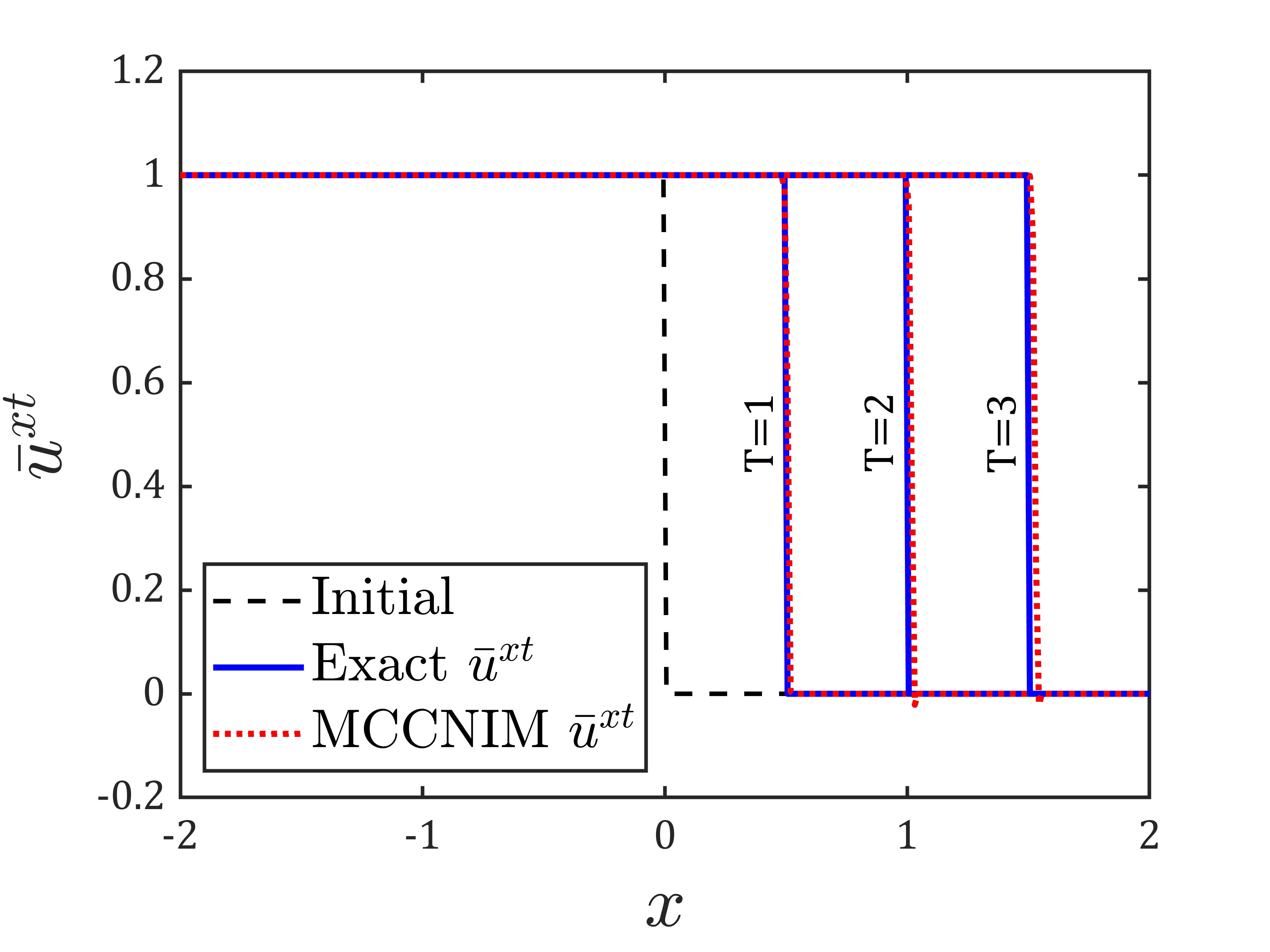}\label{fig_nx_320}}%
	\subfigure[$n_x = 640$] {\includegraphics[width=0.48\linewidth]{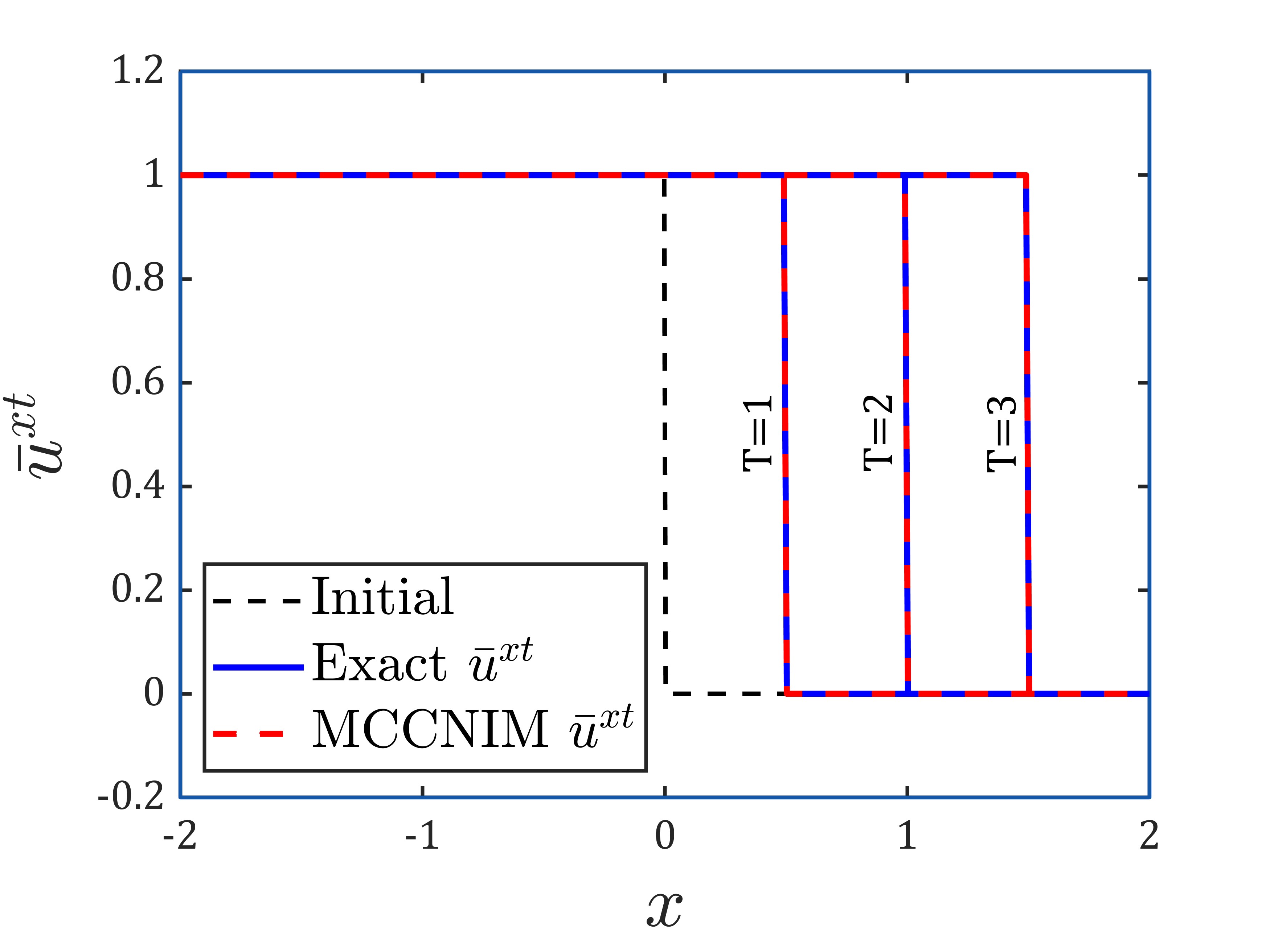}\label{fig_nx_640}}%
        \caption{\add{Comparison of the numerical solution obtained using the MCCNIM scheme with \eqn\eqref{eq:039} at a high Reynolds number ($Re=10^9$) and a fixed time step ($\Delta t = 0.0125$) for Example 1 against the exact solution, for different grid sizes.}}%
	\label{fig:High_Re}%
\end{figure}
  
\add{Further, the performance of the developed MCCNIM scheme has also been evaluated in the inviscid limit ($Re \to \infty$) by solving the considered problem (Example 1) at $Re = 10^9$, effectively reaching a nearly inviscid regime. Notably, the convective velocity (${\bar{u}}_{i,\ell}^0$) defined by \eqn\eqref{eq:040} fails to yield a converged solution at the high Reynolds numbers (i.e., $Re>10^3$).  It is due to the fact that, in this particular problem, the `diffusion-dominated' regime transitions to a `convection-dominated' regime for $Re > 10^3$, and the abrupt drop, propagating shock, from 1 to 0 leads to divergence. This divergence primarily caused by the amplification of numerical instabilities inherent in the coefficient formulations, particularly $F_{71}$, $F_{72}$, and $F_{73}$, associated with the convective velocity defined by \eqn\eqref{eq:040}. Consequently, \eqn\eqref{eq:039} is adopted to test this specific case  in the inviscid limit ($Re \to \infty$).  
The solutions computed for different grid sizes using a fixed time step ($\Delta t = 0.0125$) and the convective velocity (${\bar{u}}_{i,\ell}^0$) defined by \eqn\eqref{eq:039} are shown in \fig\ref{fig:High_Re}. 
Although a high Reynolds number case is also considered in the next example, where a stable solution is obtained, for this particular problem, numerical divergence is observed at large $Re$. However, \eqn\eqref{eq:039} has been rigorously tested across multiple cases, demonstrating accuracy comparable to the traditional MNIM, as confirmed in \tabs\ref{tab:1} to \ref{tab:3}. This consistency suggests that the scheme retains its accuracy even at high $Re$. Since viscosity is nearly negligible in this case, making the problem behave similarly to an inviscid flow, the concept of shock velocity can be used for further validation of the numerical results.}

\add{Furthermore, the shock velocity in the inviscid Burgers' equation is determined using the Rankine-Hugoniot condition \citep{kundu2024fluid,evans2022partial} as follows.}
\begin{gather}
	\begin{split}
		\add{S = \frac{u_L + u_R}{2}}
	\end{split}
\label{eq:Analytical_Shock}
\end{gather}
\add{where $u_L$ and $u_R$ are the left and right states of the shock. The viscous Burgers' equation (\eqn\ref{eq:001}), develops a thin transition layer instead of a true discontinuity (shock) for a finitely higher Reynolds number ($Re$). 
However, for sufficiently large vale of $Re$ ($=10^9$), the shock velocity remains approximately the same as in the inviscid case. In numerical simulations, the instantaneous shock location $x_s(t)$ is tracked, and the shock velocity ($S_{\text{num}}$) is computed as follows}
\begin{gather}
	\begin{split}
		\add{S_{\text{num}} = \frac{x_s(t+\Delta t) - x_s(t)}{\Delta t}}
	\end{split}
\label{eq:Numerical_Shock}
\end{gather}
\add{where $x_s(t)$ is the estimated shock position at time $t$. The shock position is identified by the location of the maximum gradient $|u_x|$.}  
\begin{table}[!t]
        \add{
	\begin{center}\renewcommand{\arraystretch}{1.5}
		\caption{Comparison of theoretical and numerical shock velocities at high Reynolds number ($Re=10^9$) for Example 1 using a fixed time step ($\Delta t = 0.0125$) for different grid sizes.}\label{tab:shock_velocity}
		\resizebox{0.35\textwidth}{!}{		
			\begin{tabular}{|c|c|c|c|}
        \hline
        $n_x$ & $S$ & $S_{\text{num}}$ & Error (\%) \\
        \hline 
        80  & 0.5 & 0.52500   & 5         \\
        160 & 0.5 & 0.51250  & 2.5       \\
        320 & 0.5 & 0.50625 & 1.25      \\
        640 & 0.5 & 0.50000     & 0  \\
        \hline
    \end{tabular}}
	\end{center}
    }
\end{table}
\add{To assess the accuracy of the numerical scheme, the theoretical shock velocity ($S$) is compared with the computed numerical shock velocity ($S_{\text{num}}$) in \tab\ref{tab:shock_velocity}. Additionally, the error between $S$ and $S_{\text{num}}$ is quantified to further validate the accuracy of the scheme. The close agreement between the two values demonstrates that the proposed MCCNIM method effectively captures shock propagation.  This comparison highlights the robustness of the MCCNIM method in resolving shock dynamics with minimal numerical dissipation.
This rigorous analysis establishes the reliability and accuracy of the results as compared with previous studies and MNIM implementation,  and strengthens confidence in the proposed MCCNIM method. 
}
%
%
\subsubsection{Example 2: Propagation and Diffusion of an Initial Sinusoidal Wave}
\label{sec:3.1.1}
The second test problem \citep{24Chen_2010} is concerned with the one-dimensional Burgers’ equation which involves an initial sinusoidal wave that propagates and diffuses within a confined flow domain along the $x$-direction \add{as follows}: 
\begin{gather}
		\begin{split}
		\frac{\partial u(x,t)}{\partial t}+u(x,t)\frac{\partial u(x,t)}{\partial x}=\frac{1}{Re}\frac{\partial^2u(x,t)}{\partial x^2}      \qquad \text{for}\quad            x\in[0,1], t\in[0,T]
	\end{split}
	\label{eq:095a}
\end{gather}
It is subject to the initial condition, and homogeneous Dirichlet boundary conditions as follows.
\begin{gather}
	\begin{split}
		u\left(x,0\right)=\sin(\pi x), \qquad    0<x<1
	\end{split}
\label{eq:095}
\end{gather}
%
%
\begin{gather}
	\begin{split}
		u\left(0,t\right)=u\left(1,t\right)=0,  \qquad     0\le t\le T
	\end{split}
\label{eq:096}
\end{gather}
The exact Fourier solution to this problem, represented by an infinite series, is given as follows.
\begin{gather}
		u\left(x,t\right) = \left(\frac{2\pi}{Re}\right)
		\frac{\sum_{k=1}^{\infty} C_k ke^{\alpha} \add{\sin(k\pi x)}}
		{C_0 + \sum_{k=1}^{\infty} C_k e^{\alpha} \add{\cos(k\pi x)}}
	\label{eq:097}
\end{gather}
where 
\begin{gather*}
	\begin{split}
		\alpha = -\left(\frac{k^2\pi^2t}{Re}\right), \qquad \beta = {-\left(\frac{Re}{2\pi}\right)(1-\cos(\pi x))}, \qquad
	C_0=\int_{0}^{1}e^{\beta}\mathrm{d}x \\
	C_k=2\int_{0}^{1}e^{\beta}\add{\cos\left(k\pi x\right)}\mathrm{d}x, \qquad k=1,2,3,\ldots   \qquad \qquad    
	\end{split}
\label{eq:097a}
\end{gather*}
\begin{figure}[!b]%
	\centering
	\subfigure[$Re=1$, $n_x=10$]{\includegraphics[width=0.48\linewidth]{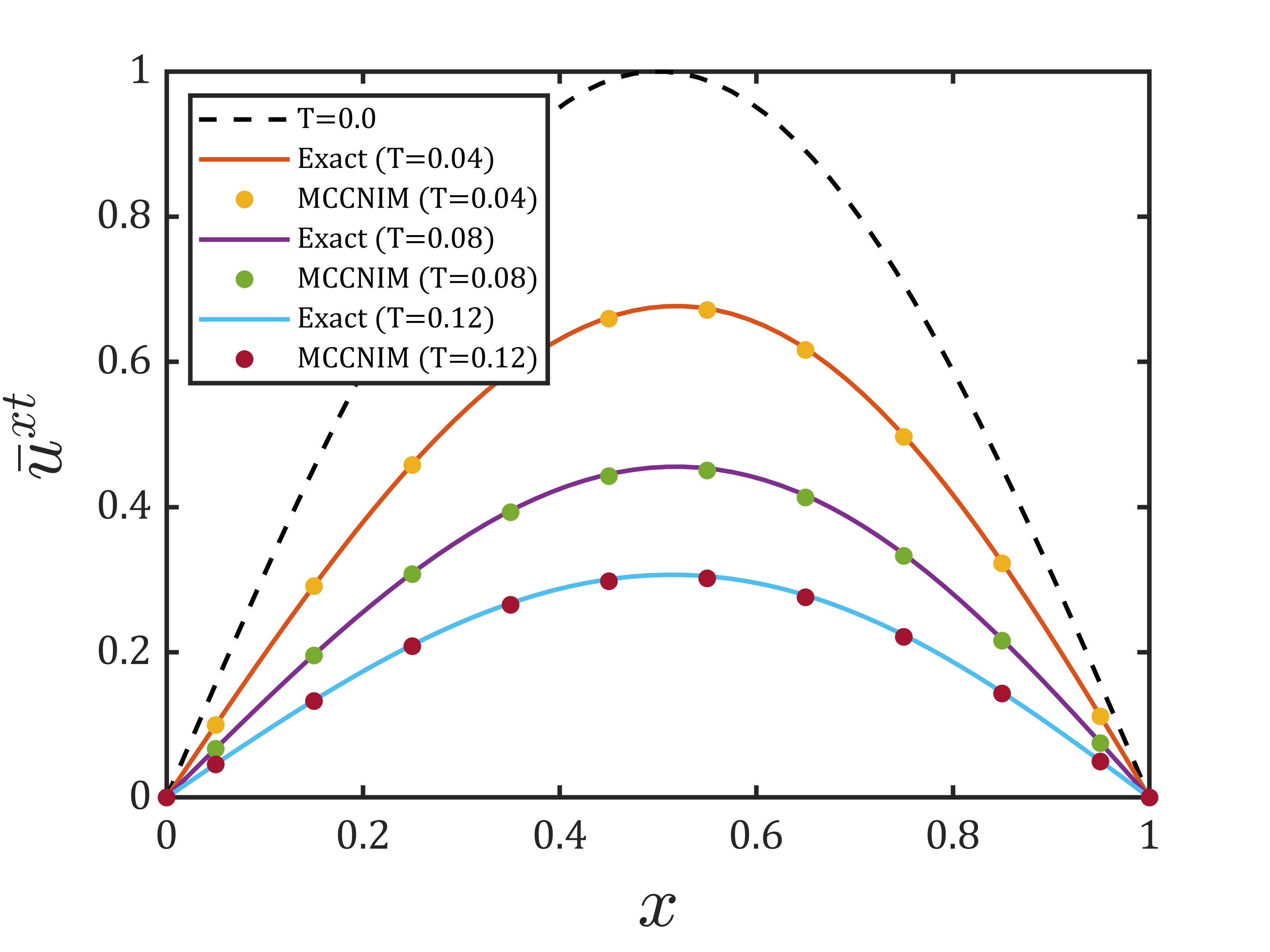}\label{fig5a}}
	\subfigure[$Re=10$, $n_x=16$] {\includegraphics[width=0.48\linewidth]{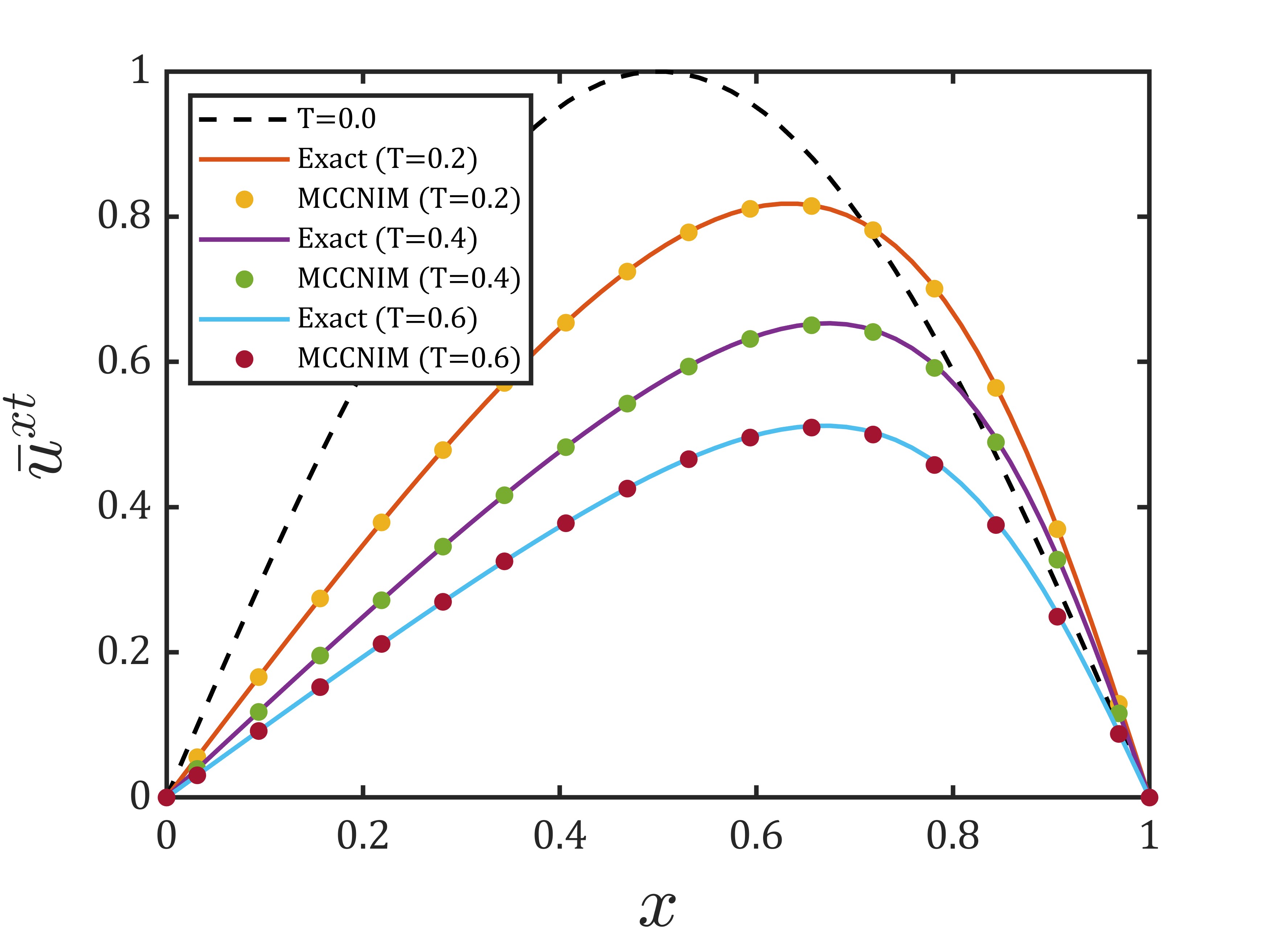}\label{fig5b}}\\%
	\subfigure[$Re=100$, $n_x=24$] {\includegraphics[width=0.48\linewidth]{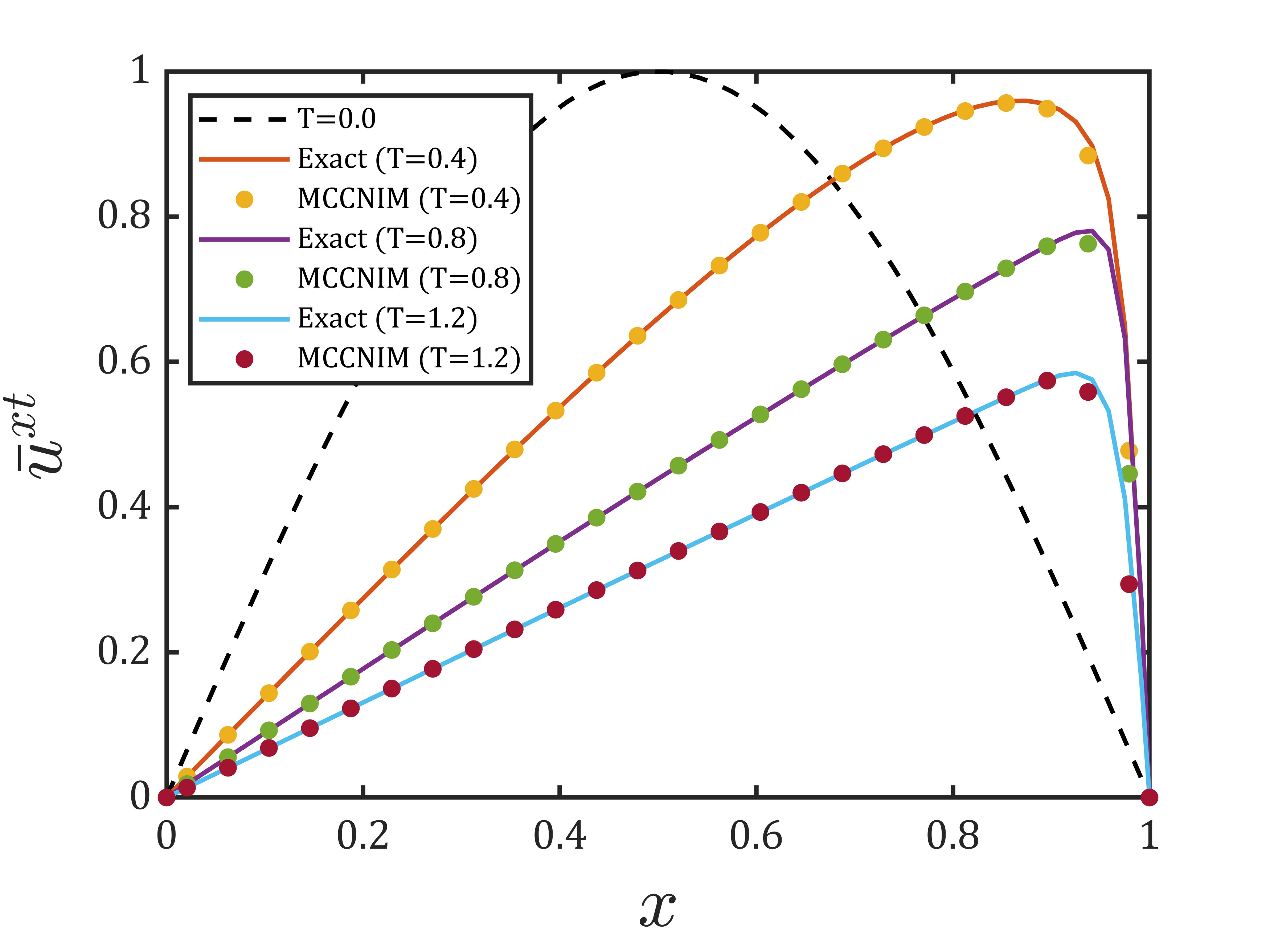}\label{fig5c}}%
	\caption{Comparison of the numerical solution obtained using the MCCNIM scheme with $\Delta t=0.001$ for Example 2 against the exact solution, evaluated at interval of $T/3$ for three different conditions.}%
	\label{fig:05}%
\end{figure}
\begin{figure}[!b]%
	\centering
	\subfigure[$Re=10^2$]{\includegraphics[width=0.48\linewidth]{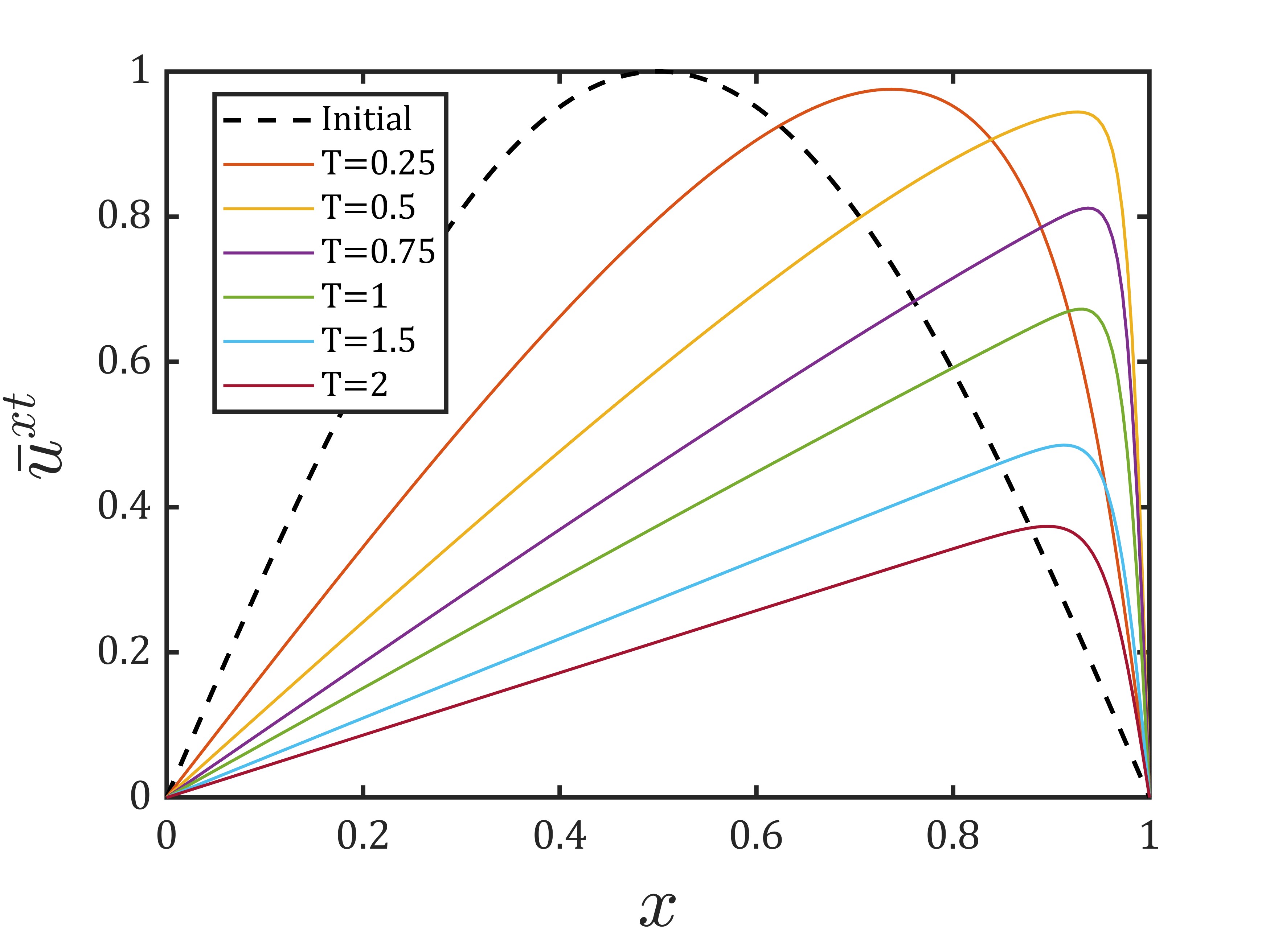}\label{fig7a}}
	\subfigure[$Re=10^3$] {\includegraphics[width=0.48\linewidth]{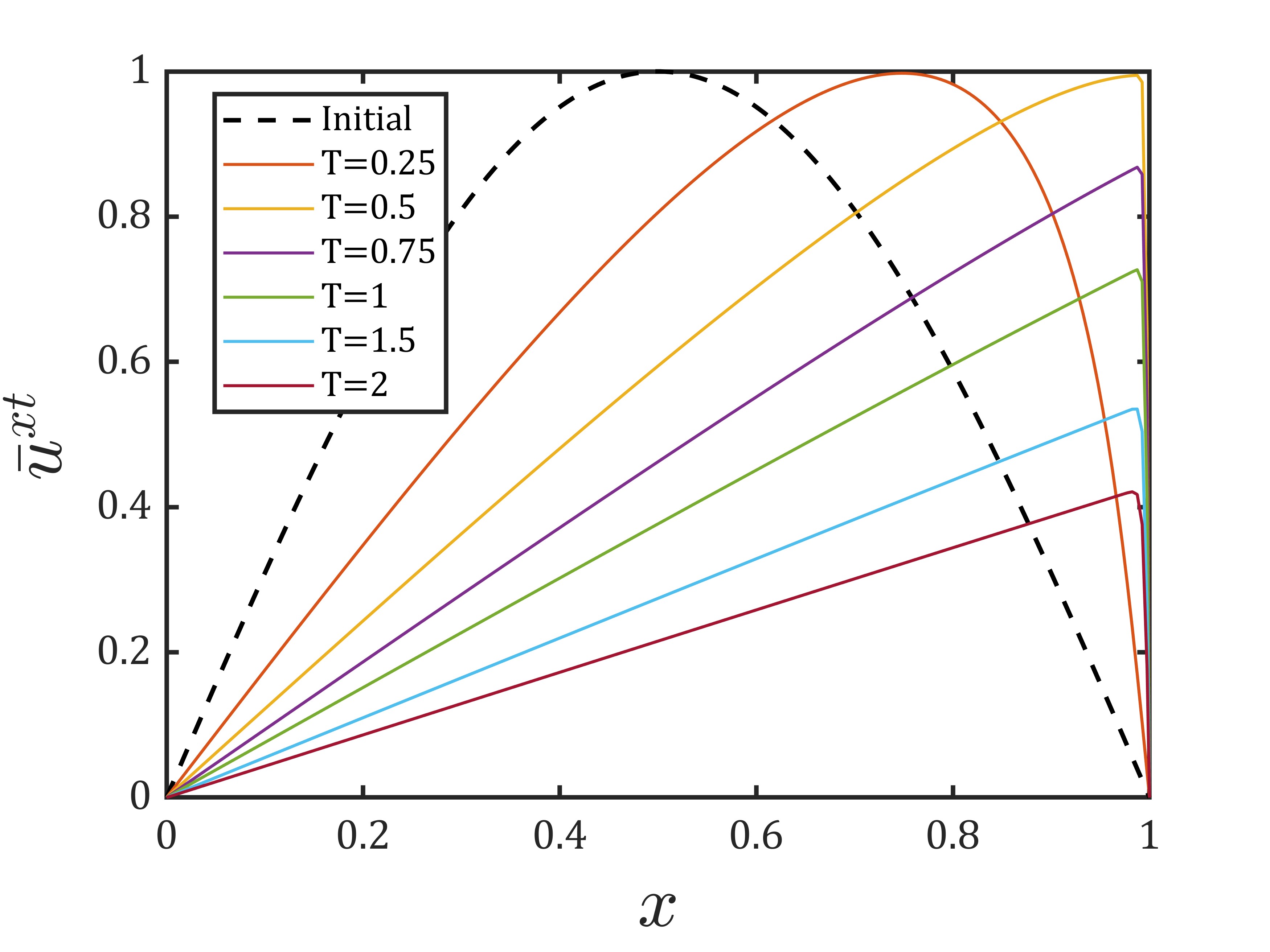}\label{fig7b}}\\
	\subfigure[$Re=10^6$] {\includegraphics[width=0.48\linewidth]{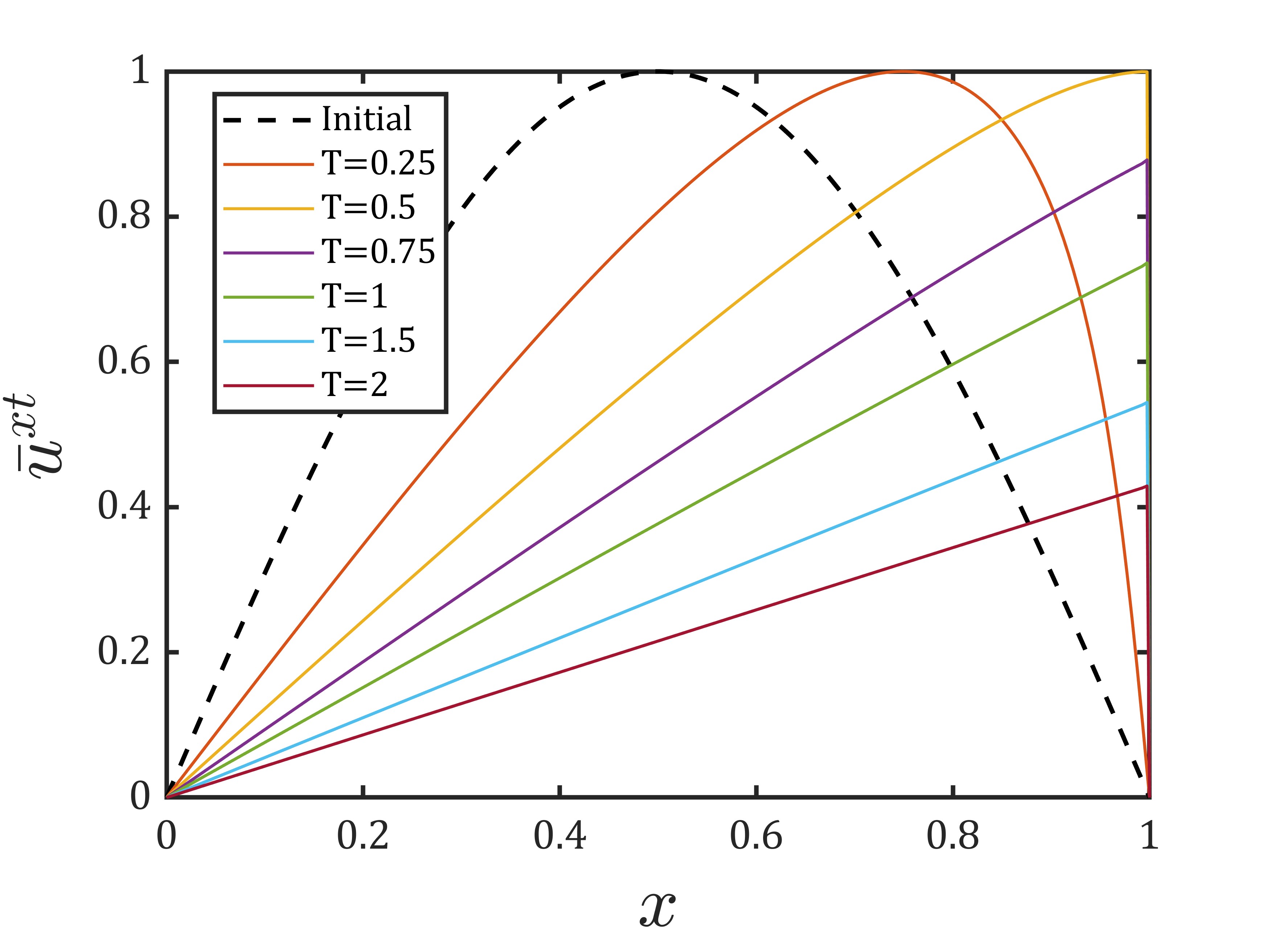}\label{fig7c}}
	\caption{Numerical solutions for Example 2 using $\Delta x = 0.0125$ and $\Delta t = 0.001$ at higher Reynolds numbers.}%
	\label{fig:07}%
\end{figure}
\begin{figure}[!b]%
	\centering
	\subfigure[Spatial rate of convergence] {\includegraphics[width=0.48\linewidth]{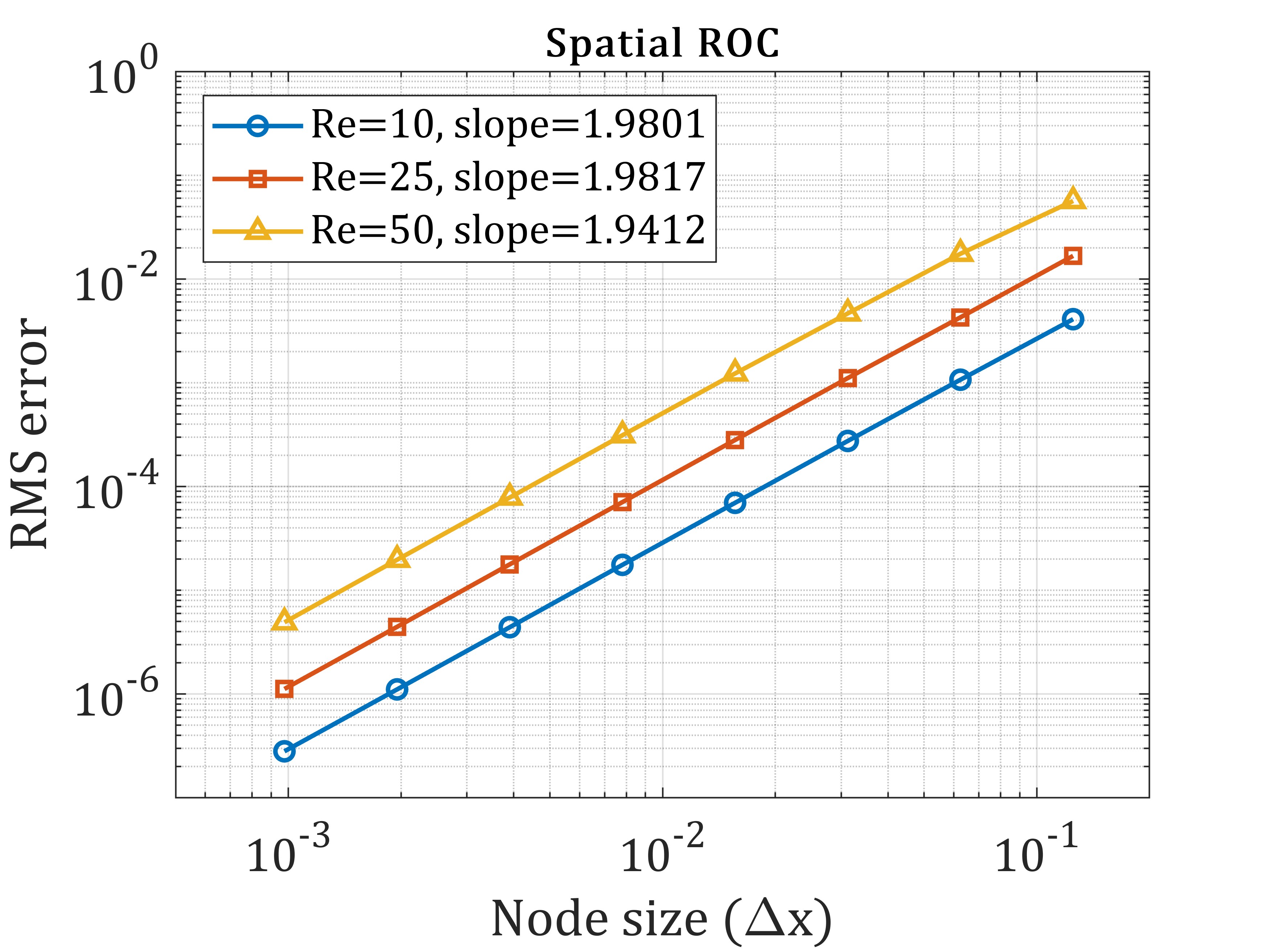}\label{fig6a}}
	\subfigure[Temporal rate of convergence] {\includegraphics[width=0.48\linewidth]{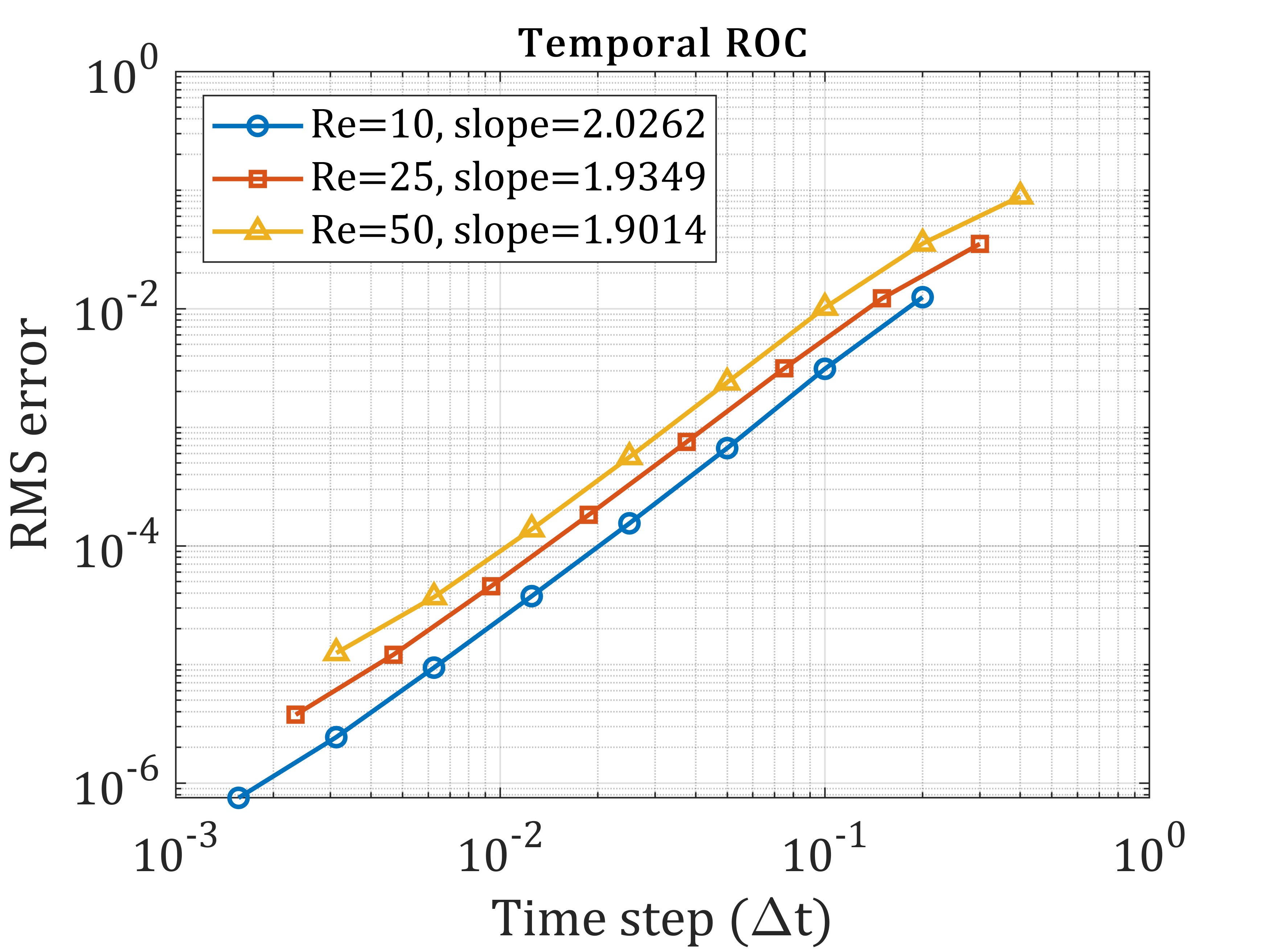}\label{fig6b}}
	\caption{Rate of convergence (ROC) for Example 2 for a range of Reynolds number ($Re$).}%
	\label{fig:06}%
\end{figure}
To validate the proposed MCCNIM scheme, the numerical results are compared with the exact solution for three different Reynolds numbers ($Re=1$, $10$, and $100$), as illustrated in \fig\ref{fig:05}. The simulations are performed up to $T=0.12$, $0.6$, and $1.2$ for $Re=1$, $10$, and $100$, respectively. In each case, the numerical (symbols) and exact (lines) solutions are compared at intervals of $T/3$, showing close alignment in all cases. \add{This demonstrates that the developed MCCNIM scheme effectively preserves the inherent accuracy of coarse-mesh methods, producing accurate solutions even on very coarse grids.} For instance, in the case of $Re=10$, as shown in panel (a) of \fig\ref{fig:05}, only $n_x=10$ grid points are sufficient to achieve an accurate match with the exact solution. 
\begin{table}[!t]
    \add{
	\begin{center}\renewcommand{\arraystretch}{1.5}
        \caption{Comparison of RMS errors between the developed MCCNIM and MNIM schemes for different node sizes ($\Delta x$) in Example 2, using a smaller time step ($\Delta t = 0.0001$).}\label{tab:4}
		\resizebox{1\textwidth}{!}{		
		  \begin{tabular}{|c|c|c|c|c|c|c|}
		      \hline
			    $\Delta x$ & \multicolumn{2}{c|}{$Re = 10$ at $T = 0.2$} & \multicolumn{2}{c|}{$Re = 25$ at $T = 0.3$} & \multicolumn{2}{c|}{$Re = 50$ at $T = 0.4$} \\ \cline{2-7}
			     & MNIM ($\bar{u}^{x}$)  & MCCNIM ($\bar{u}^{xt}$) & MNIM ($\bar{u}^{x}$)  & MCCNIM ($\bar{u}^{xt}$) & MNIM ($\bar{u}^{x}$)  & MCCNIM ($\bar{u}^{xt}$) \\ 
			  \hline
			    0.125   & \(4.59 \times 10^{-3}\)  & \(4.10 \times 10^{-3}\)  & \(1.87 \times 10^{-2}\)  & \(1.67 \times 10^{-2}\)  & \(6.32 \times 10^{-2}\)  & \(5.65 \times 10^{-2}\)  \\
                0.0625    & \(1.14 \times 10^{-3}\)  & \(1.07 \times 10^{-3}\)  & \(4.53 \times 10^{-3}\)  & \(4.27 \times 10^{-3}\)  & \(1.85 \times 10^{-2}\)  & \(1.75 \times 10^{-2}\)  \\
                0.03125   & \(2.83 \times 10^{-4}\)  & \(2.75 \times 10^{-4}\)  & \(1.14 \times 10^{-3}\)  & \(1.11 \times 10^{-3}\)  & \(4.81 \times 10^{-3}\)  & \(4.66 \times 10^{-3}\)  \\
                0.015625  & \(7.08 \times 10^{-5}\)  & \(6.97 \times 10^{-5}\)  & \(2.86 \times 10^{-4}\)  & \(2.81 \times 10^{-4}\)  & \(1.25 \times 10^{-3}\)  & \(1.23 \times 10^{-3}\)  \\
                0.0078125 & \(1.77 \times 10^{-5}\)  & \(1.76 \times 10^{-5}\)  & \(7.14 \times 10^{-5}\)  & \(7.08 \times 10^{-5}\)  & \(3.14 \times 10^{-4}\)  & \(3.12 \times 10^{-4}\)  \\
                0.00390625& \(4.43 \times 10^{-6}\)  & \(4.41 \times 10^{-6}\)  & \(1.79 \times 10^{-5}\)  & \(1.78 \times 10^{-5}\)  & \(7.88 \times 10^{-5}\)  & \(7.84 \times 10^{-5}\)  \\
                0.001953125& \(1.10 \times 10^{-6}\) & \(1.11 \times 10^{-6}\)  & \(4.47 \times 10^{-6}\)  & \(4.45 \times 10^{-6}\)  & \(1.97 \times 10^{-5}\)  & \(1.97 \times 10^{-5}\)  \\
                0.0009765625 & \(2.86 \times 10^{-7}\) & \(2.77 \times 10^{-7}\)  & \(1.12 \times 10^{-6}\)  & \(1.12 \times 10^{-6}\)  & \(4.93 \times 10^{-6}\)  & \(4.92 \times 10^{-6}\)  \\
			  \hline
		  \end{tabular}}
	\end{center}
    }
\end{table}
\begin{table}[!t]
    \add{
	\begin{center}\renewcommand{\arraystretch}{1.5}
		\caption{RMS errors for different {time step ($\Delta t$)} in Example 2, {using the MCCNIM scheme with a fine node size ($\Delta x = 0.0001$).}}\label{tab:5}
		\resizebox{0.75\textwidth}{!}
		{
		\begin{tabular}{|c|c|c|c|c|c|} \hline
		\multicolumn{2}{|c|}{ $Re = 10$ at $T = 0.2$} & \multicolumn{2}{c|}{$Re = 25$ at $T = 0.3$} & \multicolumn{2}{c|}{$Re = 50$ at $T = 0.4$} \\ \hline
		$\Delta t$ & RMS errors & $\Delta t$ & RMS errors & $\Delta t$ & RMS errors \\ \hline
		0.2      & 1.25$\times 10^{-2}$     & 0.3      & 3.52$\times 10^{-2}$    & 0.4      & 8.84$\times 10^{-2}$   \\
		0.1      & 3.11$\times 10^{-3}$   & 0.15     & 1.22$\times 10^{-2}$   & 0.2      & 3.54$\times 10^{-2}$  \\
		0.05     & 6.67$\times 10^{-4}$   & 0.075    & 3.14$\times 10^{-3}$   & 0.1      & 1.02$\times 10^{-2}$  \\
		0.025    & 1.55$\times 10^{-4}$   & 0.0375   & 7.52$\times 10^{-4}$   & 0.05     & 2.39$\times 10^{-3}$  \\
		0.0125   & 3.77$\times 10^{-5}$   & 0.01875  & 1.83$\times 10^{-4}$   & 0.025    & 5.61$\times 10^{-4}$  \\
		0.00625  & 9.41$\times 10^{-6}$   & 0.009375 & 4.57$\times 10^{-5}$   & 0.0125   & 1.38$\times 10^{-4}$  \\
		0.003125  & 2.44$\times 10^{-6}$   & 0.0046875  & 1.21$\times 10^{-5}$   & 0.00625  & 3.70$\times 10^{-5}$  \\
		0.0015625  & 7.51$\times 10^{-7}$   & 0.00234375  & 3.77$\times 10^{-6}$   & 0.003125  & 1.25$\times 10^{-5}$  \\ \hline
		\end{tabular}
		}
	\end{center}
    }
\end{table}
\add{However, for higher Reynolds numbers ($Re > 100$), it is known that the Fourier solution fails to converge because of the slow convergence of the infinite series \citep{miller1966predictor,24Chen_2010}. Therefore, the proposed MCCNIM scheme has been validated for problems with high non-linearity (dominating with increasing $Re$) by analyzing the numerical results for the considered problem for the high Reynolds numbers ($Re=10^2 - 10^6$) at different time points, as shown in \fig\ref{fig:07}.}
Evidently, the scheme operates smoothly without any oscillations near the discontinuous region ($x=1$), even at a very high Reynolds numbers ($Re=10^6$). The curves exhibit the expected physical behavior, validating the effectiveness of the developed MCCNIM scheme for problems characterized by high non-linear behavior.

\add{Furthermore, the efficacy of the established method is illustrated by analyzing the spatial and temporal variations of RMS errors for a range of Reynolds number values. \tab\ref{tab:4} depicts the influence of grid sizes ($\Delta x$) on the RMS error for the developed MCCNIM scheme in comparison with the MNIM scheme for a range of Reynolds numbers using the fixed, smaller time step ($\Delta t=0.0001$). Similarly, \tab\ref{tab:5} presents the influence of time step ($\Delta t$) on the RMS error of the developed MCCNIM scheme for different $Re$ values using the fine grid ($\Delta x = 0.0001$).  The tabular data clearly demonstrate excellent convergence of MCCNIM scheme, as the RMS error decreases with refinement in grid ($\Delta x$) and time step ($\Delta t$). This trend holds regardless of Reynolds number ($Re$) values, consistent with the MNIM scheme.
Further, the rate of convergence (ROC) analysis is performed to obtain the convergence order of the developed MCCNIM scheme by plotting the RMS error as a function of grid size ($\Delta x$) and time step ($\Delta t$) for a range of Reynolds numbers ($Re$) in \fig\ref{fig:06}. Interestingly, the relationship of RMS error with both grid size ($\Delta x$) and time step ($\Delta t$) is linear with slope of approxminately $\approx 2$. This crucial demonstration emphasizes that the developed MCCNIM scheme has quadratic convergence, i.e., second-order accurate in both space and time, O[($\Delta x$)$^2$, ($\Delta t$)$^2$], even for higher Reynolds numbers ($Re$), as indicated in panels (a) and (b) of \fig\ref{fig:06} with spatial rate of convergence (SROC $\approx 2$) and temporal rate of convergence (TROC $\approx 2$). In addition to the previous analysis, these results further confirm that MCCNIM excellently matches or outperforms traditional nodal schemes, reinforcing the accuracy and reliability of the MCCNIM scheme.}
\begin{table}[!t]
    \add{
    \begin{center}\renewcommand{\arraystretch}{1.5}
		\caption{Performance evaluation of the proposed MCCNIM scheme for Example 2 at various Reynolds numbers ($Re$) using $T=2$ and $\Delta t=0.1$.}\label{tab:6}
		\resizebox{0.85\textwidth}{!}
		{
		\begin{tabular}{|c|c|c|c|c|c|c|}
        \hline
        \multirow{2}{*}{$\Delta x$} & \multicolumn{2}{c|}{Picard iterations (Nonlinear)} & \multicolumn{2}{c|}{Krylov iterations (Linear)} & \multicolumn{2}{c|}{CPU Runtime (seconds)} \\ \cline{2-7}
        & MNIM & MCCNIM & MNIM & MCCNIM & MNIM & MCCNIM \\
        \hline
        \multicolumn{7}{|c|}{Re = 50} \\
        \hline
        0.1  & 110 & 126 & 1492  & 1908  & 0.01354  & 0.019089 \\
        0.05  & 132 & 151 & 3030  & 2502  & 0.026312 & 0.026295 \\
        0.025  & 147 & 169 & 6151  & 4173  & 0.063205 & 0.050799 \\
        0.0125  & 153 & 177 & 12444 & 7219  & 0.210863 & 0.112289 \\
        0.00625 & 158 & 181 & 25835 & 12400 & 0.950768 & 0.332598 \\
        0.003125 & 162 & 185 & 56182 & 21488 & 6.640158 & 1.196723 \\
        0.0015625 & 201 & 190 & 153380 & 37789 & 65.491313 & 5.818726 \\
        \hline
        \multicolumn{7}{|c|}{Re = $10^3$} \\
        \hline
        0.1  & 108 & 124 & 1678  & 2652  & 0.013859 & 0.026348 \\
        0.05  & 115 & 128 & 3014  & 2977  & 0.02345  & 0.030108 \\
        0.025  & 120 & 136 & 4369  & 3055  & 0.042249 & 0.036187 \\
        0.0125  & 124 & 145 & 7113  & 4338  & 0.101176 & 0.066436 \\
        0.00625 & 152 & 175 & 12695 & 6255  & 0.321908 & 0.130973 \\
        0.003125 & 190 & 221 & 24307 & 10031 & 1.535329 & 0.36623  \\
        0.0015625 & 215 & 251 & 50913 & 17966 & 9.942712 & 1.564128 \\
        \hline
        \multicolumn{7}{|c|}{Re = $10^6$} \\
        \hline
        0.1  & 108 & 124 & 1722  & 2738  & 0.017558 & 0.027129 \\
        0.05  & 115 & 128 & 3443  & 4238  & 0.027329 & 0.043778 \\
        0.025  & 119 & 133 & 6367  & 4976  & 0.06857  & 0.059237 \\
        0.0125  & 124 & 139 & 11364 & 5878  & 0.201471 & 0.08875  \\
        0.00625 & 127 & 144 & 21485 & 8129  & 0.746016 & 0.190523 \\
        0.003125 & 130 & 145 & 42866 & 12760 & 4.474131 & 0.588388 \\
        0.0015625 & 133 & 150 & 87503 & 23352 & 32.029748 & 3.008782 \\
        \hline
    \end{tabular}
		}
	\end{center}
    }
\end{table}
%

\add{To further evaluate the effectiveness of the developed MCCNIM scheme, Table \ref{tab:6} presents a detailed performance of the present MCCNIM scheme in comparison to MNIM scheme for Example 2, evaluating key computational metrics such as Picard iterations (for nonlinear solves), Krylov iterations (for linear solves), and CPU runtime for the wide range of Reynolds number ($Re = 50 - 10^6$) and various spatial resolutions ($\Delta x$) while keeping a fixed time step ($\Delta t = 0.1$) and a final time ($T = 2$). Note that the Picard and Krylov iterations listed in Table \ref{tab:6} represent the total counts of nonlinear and linear iterations, respectively. In contrast, the CPU runtime indicates the overall computational time (in seconds) required for the complete simulations. Here, "total" refers to the cumulative counts from lines 3 to 16 in Algorithm \ref{Algorithm1}. Table \ref{tab:6} indicates that while MCCNIM requires slightly more Picard iterations in some cases, it significantly reduces the number of Krylov iterations, particularly for finer spatial resolutions. This reduction in linear iterations becomes more pronounced as $Re$ increases, indicating improved stability and efficiency of the MCCNIM scheme for high Reynolds number flows. Additionally, the CPU runtime for MCCNIM is consistently lower than that of MNIM for finer spatial resolutions, demonstrating the computational advantage of the proposed scheme. These results further confirm that MCCNIM maintains accuracy and enhances computational efficiency, making it a promising approach for solving the Burger's equation.}

\subsubsection{Example 3: One-Dimensional Simulation of Shock-Like Formation}
\add{This example is commonly used to represent the shock-like formations governed by Burgers' equation subject to the initial conditions and Dirichlet boundary conditions expressed as follows. 
\begin{gather}
	\begin{split}
		\frac{\partial u(x,t)}{\partial t}+u(x,t)\frac{\partial u(x,t)}{\partial x}=\frac{1}{Re}\frac{\partial^2u(x,t)}{\partial x^2}      \qquad \text{for}\quad            x\in[0,1], t\in[0,T]
	\end{split}
	\label{eq:105a}
\end{gather}
\begin{gather}
	\begin{split}
		u\left(x,0\right)=\sin(2\pi x), \qquad    x\in[0,1]
	\end{split}
	\label{eq:105}
\end{gather}
\begin{gather}
	\begin{split}
		u\left(0,t\right)=u\left(1,t\right)=0,  \qquad     t\in[0,T]
	\end{split}
	\label{eq:106}
\end{gather}
The exact Fourier solution to Burgers' equation is provided in \eqn\eqref{eq:097}. In this case, the Fourier coefficients, represented by an infinite series, are given as follows:
\begin{gather*}
	\begin{split}
		C_0=\int_{0}^{1}e^{-\frac{Re}{4\pi}(1-\cos(2\pi x))}dx,\qquad
		C_k=2\int_{0}^{1}e^{-\frac{Re}{4\pi}(1-\cos(2\pi x))}\add{\cos\left(k\pi x\right)}dx,\qquad k=1,2,3,\ldots     
	\end{split}
	\label{eq:104}
\end{gather*}
}
\begin{figure}[!b]%
	\centering
	\subfigure[$Re=1$, $n_x=10$]{\includegraphics[width=0.48\linewidth]{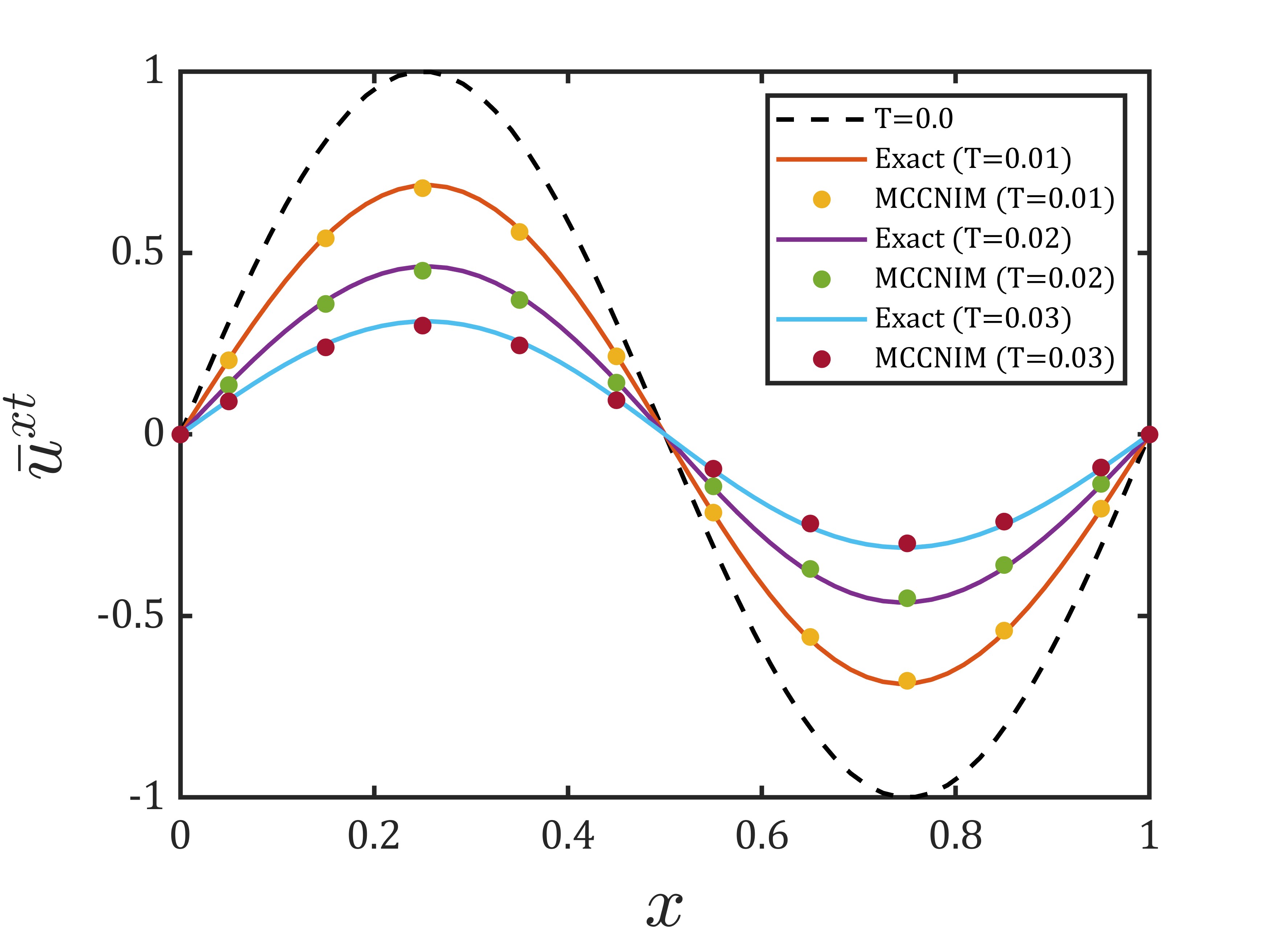}\label{fig13a}}
	\subfigure[$Re=10$, $n_x=16$] {\includegraphics[width=0.48\linewidth]{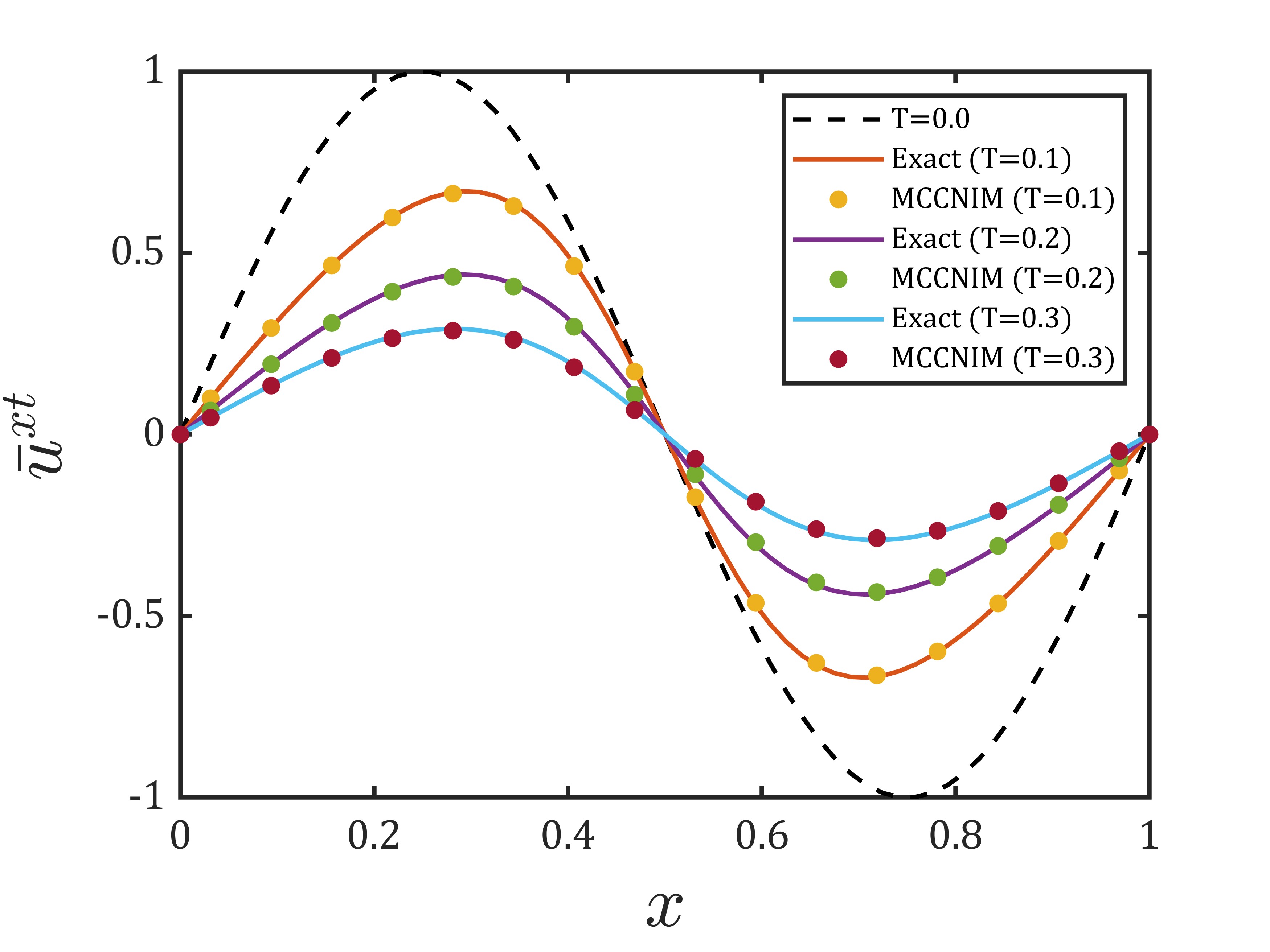}\label{fig13b}}\\%
	\subfigure[$Re=100$, $n_x=24$] {\includegraphics[width=0.48\linewidth]{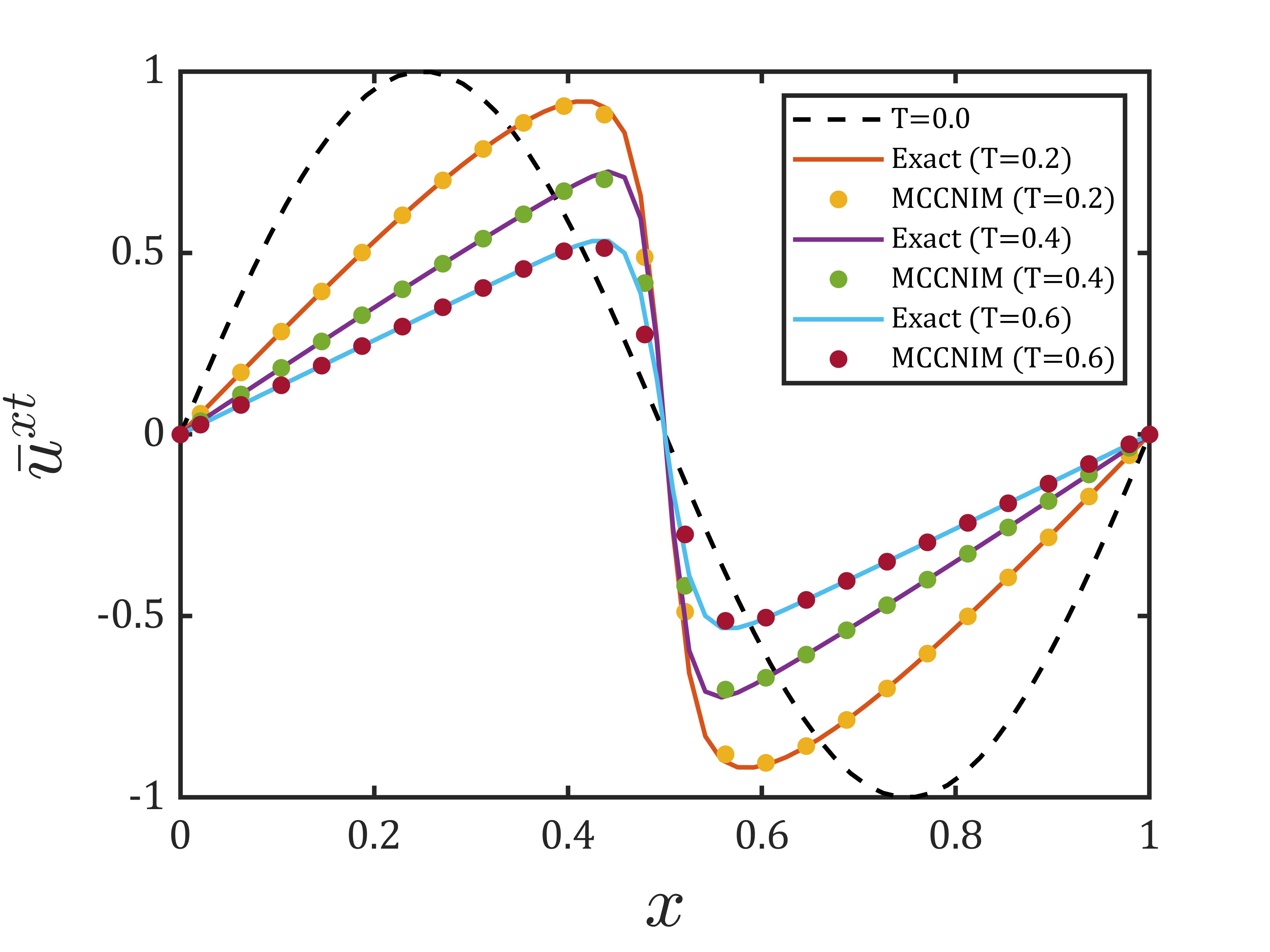}\label{fig13c}}%
	\caption{\add{Comparison of the numerical solution obtained using the MCCNIM scheme with $\Delta t =0.001$ against the exact solution for Example 3, at interval of $T/3$ for three different values of Reynolds numbers.}}%
	\label{fig:13}%
\end{figure}
\begin{figure}[!b]%
	\centering
	\subfigure[$Re=10^2$]{\includegraphics[width=0.48\linewidth]{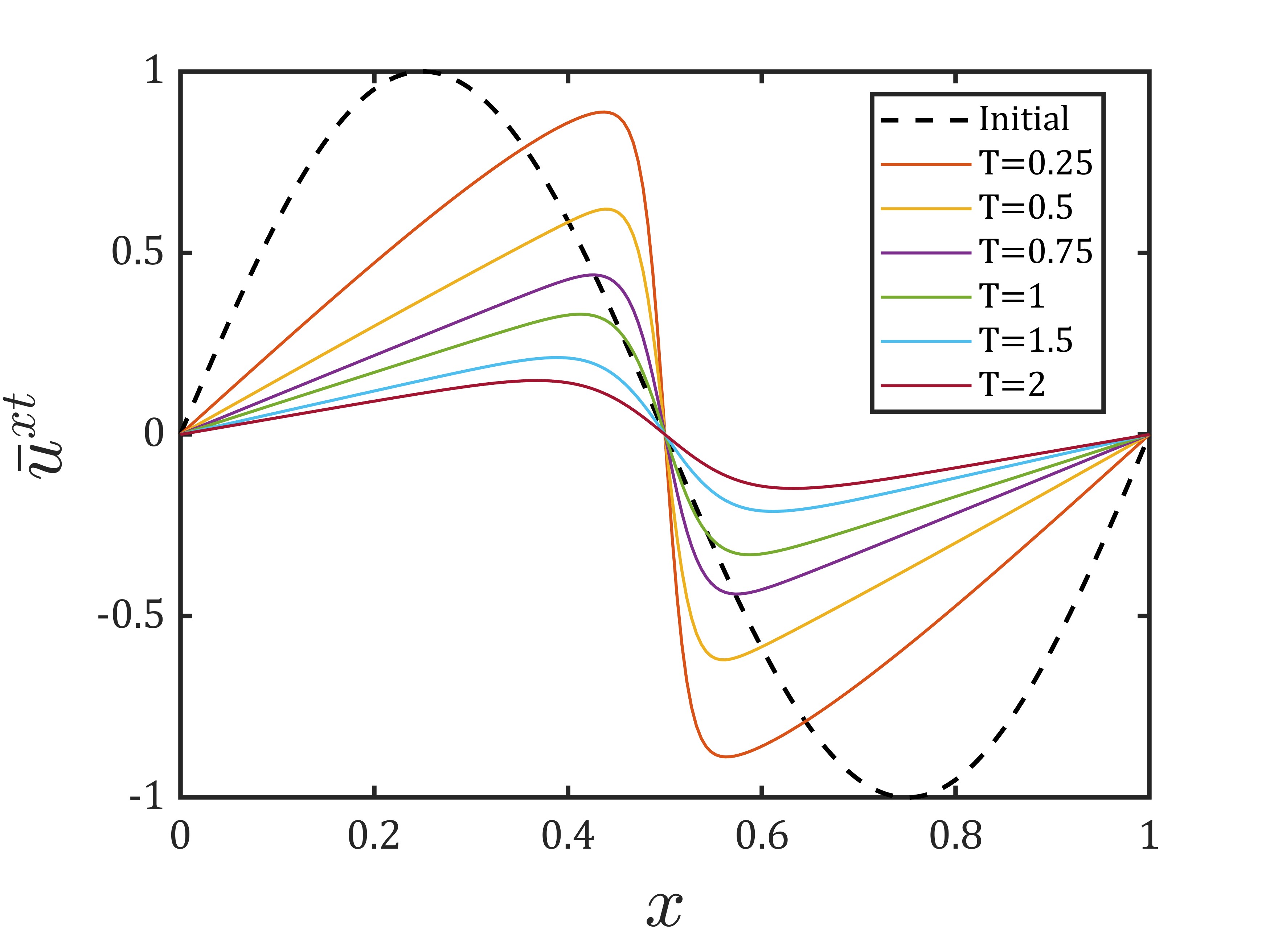}\label{fig14a}}
	\subfigure[$Re=10^3$] {\includegraphics[width=0.48\linewidth]{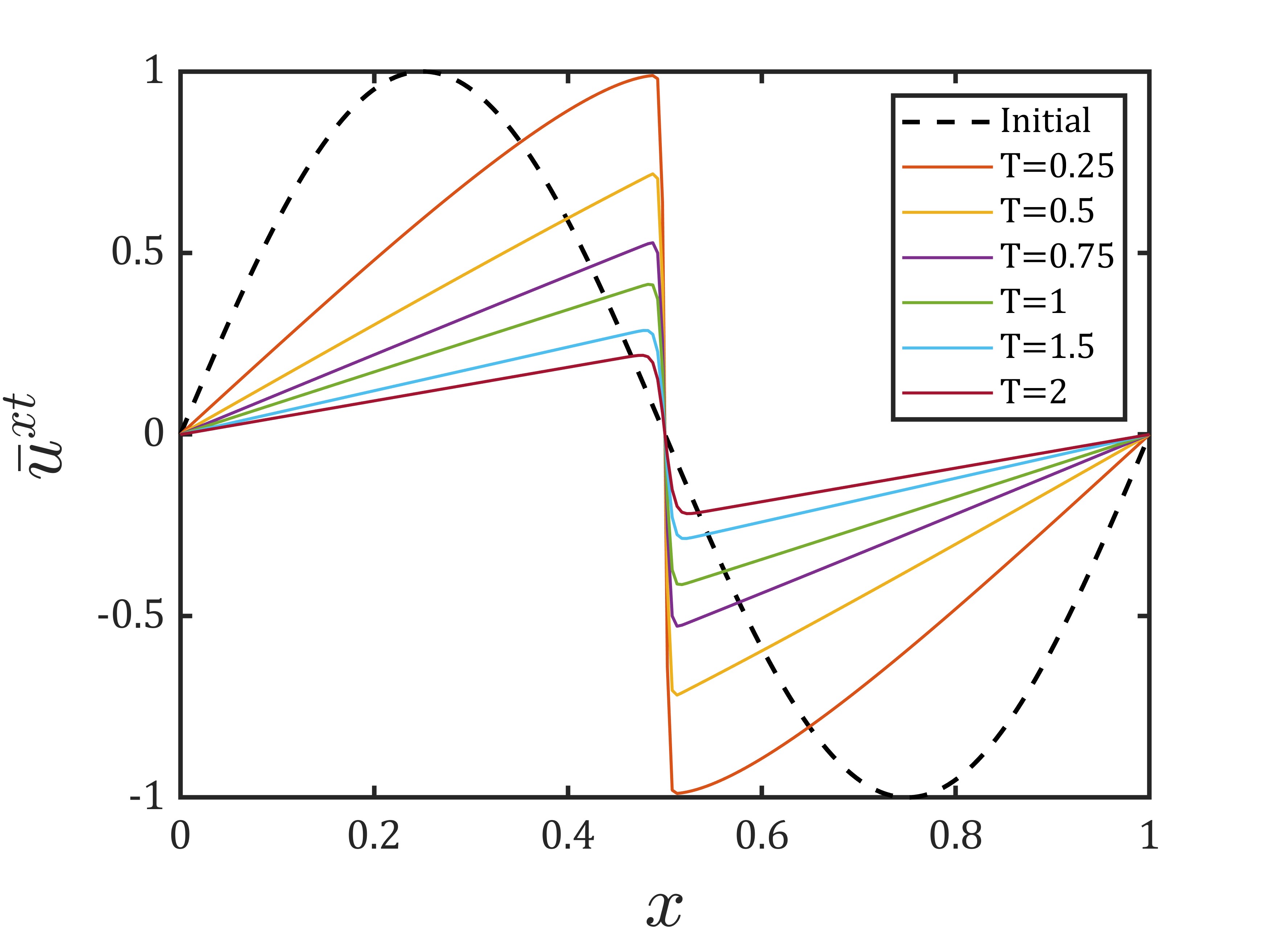}\label{fig14b}}\\
	\subfigure[$Re=10^6$] {\includegraphics[width=0.48\linewidth]{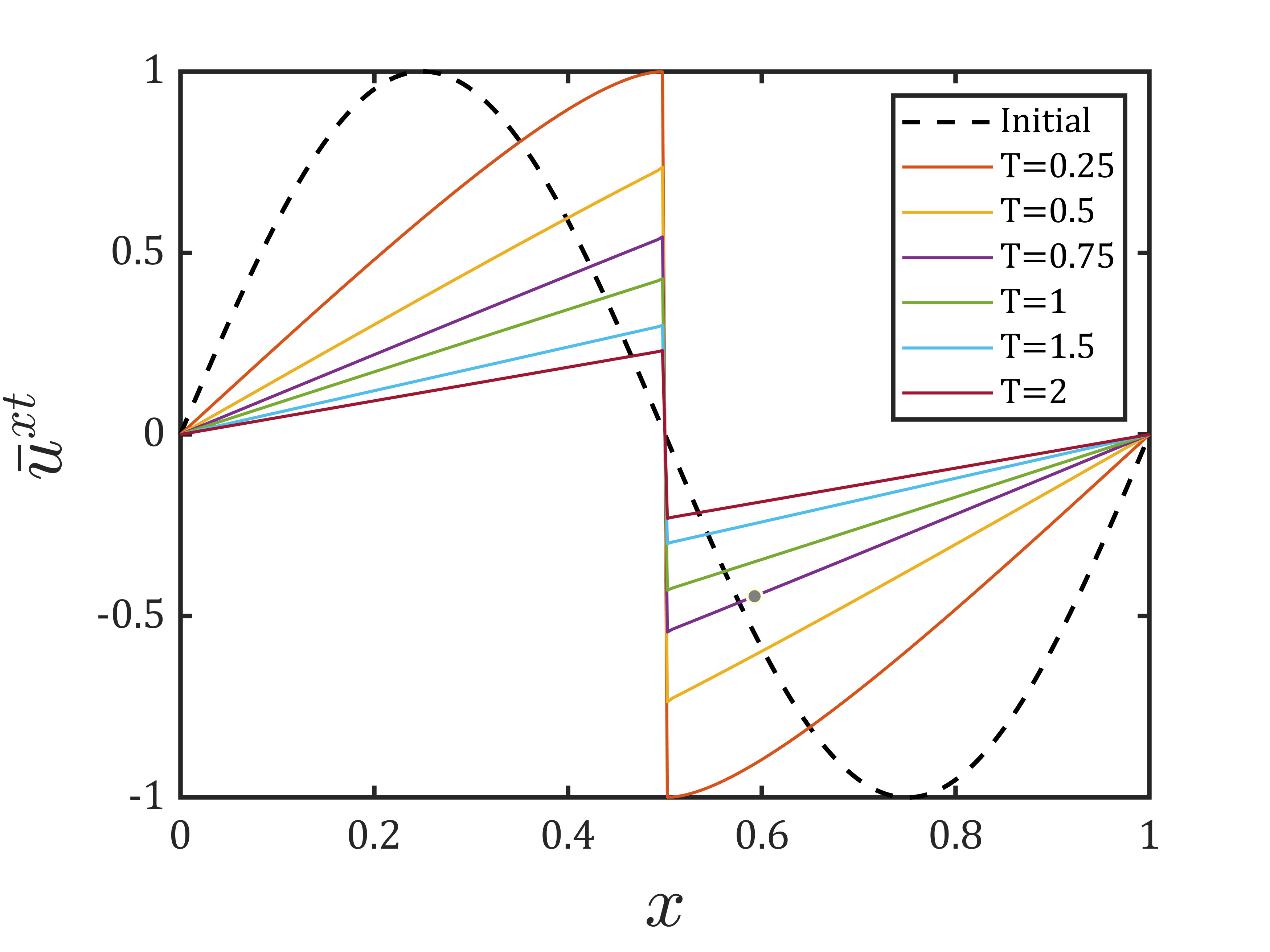}\label{fig14c}}
	\caption{\add{Numerical solutions for Example 3 at higher Reynolds numbers using the MCCNIM scheme  with  $\Delta x = 0.0125$ and $\Delta t = 0.001$.}}%
	\label{fig:14}%
\end{figure}
\add{The numerical experiments for Example 3 are conducted following the same approach as in Example 2. The numerical results are compared with the exact solution for various Reynolds number ($Re$) values and grid sizes ($\Delta x$) to validate the accuracy of the proposed MCCNIM scheme. \fig\ref{fig:13} illustrates that the numerical solutions exhibit strong agreement with the exact solutions, for a range of Reynolds number ($Re=1-100$) and time ($T=0-0.6$), reinforcing the effectiveness of the MCCNIM method.}
\begin{figure}[!b]%
	\centering
	\subfigure[Spatial rate of convergence] {\includegraphics[width=0.48\linewidth]{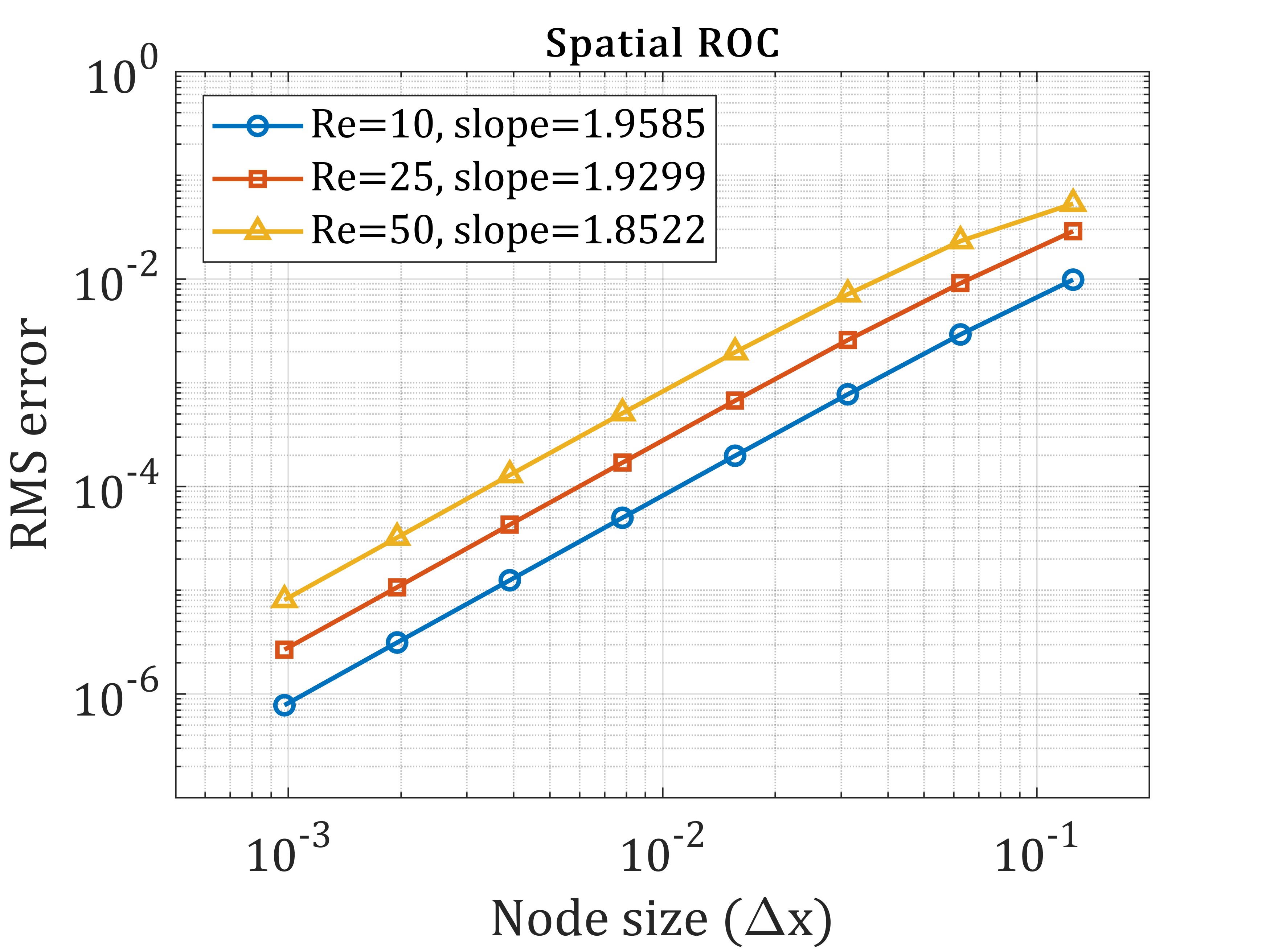}\label{fig15a}}
	\subfigure[Temporal rate of convergence] {\includegraphics[width=0.48\linewidth]{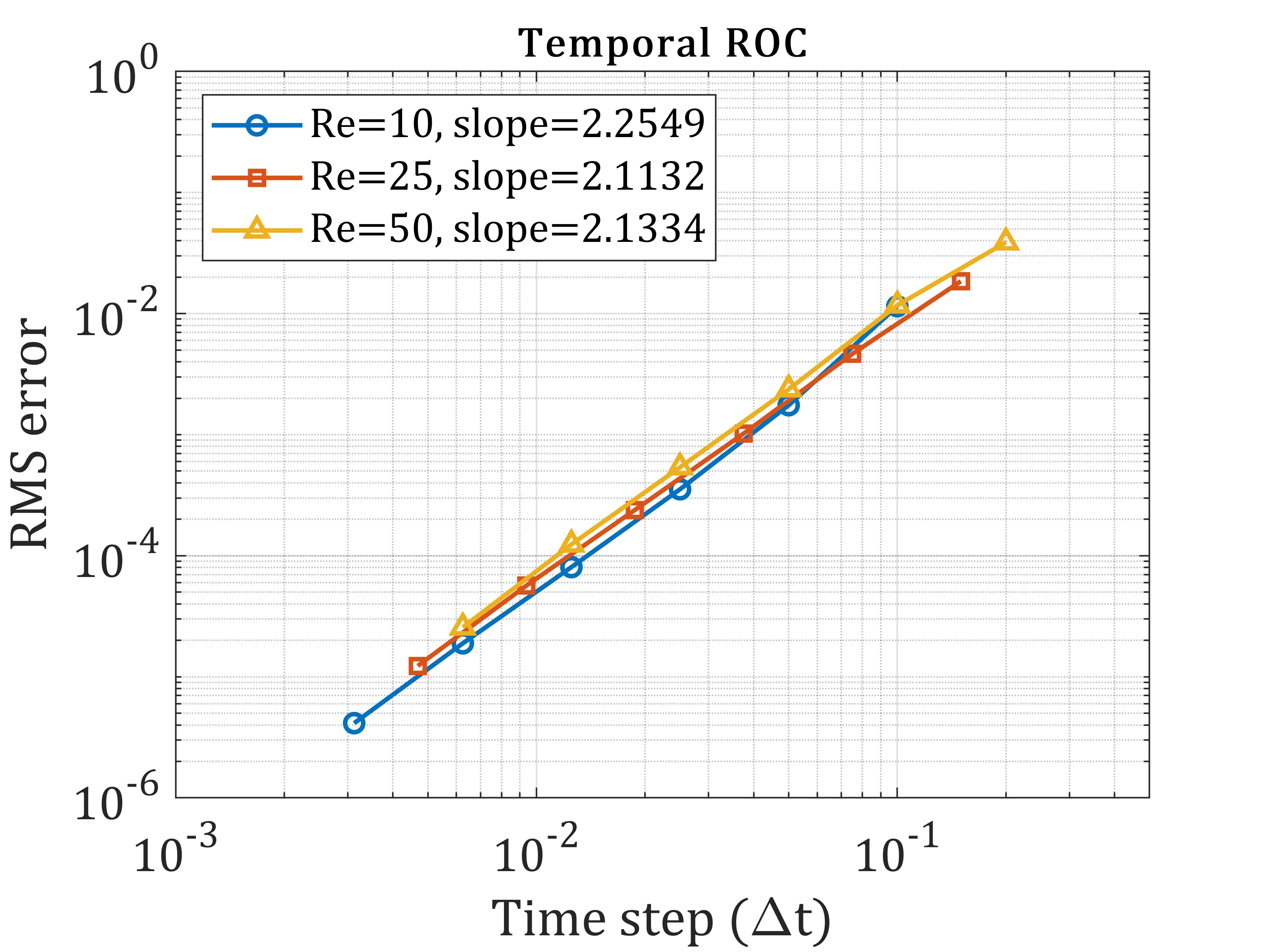}\label{fig15b}}
	\caption{\add{Rate of convergence (ROC) for Example 3 for a range of Reynolds number ($Re$).}}%
	\label{fig:15}%
\end{figure}
\begin{table}[!t]
    \add{
	\begin{center}\renewcommand{\arraystretch}{1.5}
        \caption{Comparison of RMS errors in solutions for Example 3 obtained using the developed MCCNIM and MNIM schemes for different node sizes ($\Delta x$) using a smaller time step ($\Delta t = 0.0001$).}\label{tab:7}
		\resizebox{1\textwidth}{!}{		
		  \begin{tabular}{|c|c|c|c|c|c|c|}
		      \hline
			    $\Delta x$ & \multicolumn{2}{c|}{$Re = 10$ at $T = 0.2$} & \multicolumn{2}{c|}{$Re = 25$ at $T = 0.3$} & \multicolumn{2}{c|}{$Re = 50$ at $T = 0.4$} \\ \cline{2-7}
			     & MNIM ($\bar{u}^{x}$)  & MCCNIM ($\bar{u}^{xt}$) & MNIM ($\bar{u}^{x}$)  & MCCNIM ($\bar{u}^{xt}$) & MNIM ($\bar{u}^{x}$)  & MCCNIM ($\bar{u}^{xt}$) \\ 
			  \hline
			    0.125     & \(1.10 \times 10^{-2}\)  & \(9.85 \times 10^{-3}\)  & \(3.23 \times 10^{-2}\)  & \(2.89 \times 10^{-2}\)  & \(5.95 \times 10^{-2}\)  & \(5.32 \times 10^{-2}\)  \\
                0.0625    & \(3.10 \times 10^{-3}\)  & \(2.92 \times 10^{-3}\)  & \(9.67 \times 10^{-3}\)  & \(9.11 \times 10^{-3}\)  & \(2.46 \times 10^{-2}\)  & \(2.32 \times 10^{-2}\)  \\
                0.03125   & \(7.98 \times 10^{-4}\)  & \(7.74 \times 10^{-4}\)  & \(2.66 \times 10^{-3}\)  & \(2.58 \times 10^{-3}\)  & \(7.34 \times 10^{-3}\)  & \(7.12 \times 10^{-3}\)  \\
                0.015625  & \(2.01 \times 10^{-4}\)  & \(1.98 \times 10^{-4}\)  & \(6.83 \times 10^{-4}\)  & \(6.72 \times 10^{-4}\)  & \(2.00 \times 10^{-3}\)  & \(1.97 \times 10^{-3}\)  \\
                0.0078125 & \(5.03 \times 10^{-5}\)  & \(5.00 \times 10^{-5}\)  & \(1.72 \times 10^{-4}\)  & \(1.70 \times 10^{-4}\)  & \(5.13 \times 10^{-4}\)  & \(5.09 \times 10^{-4}\)  \\
                0.00390625& \(1.26 \times 10^{-5}\)  & \(1.25 \times 10^{-5}\)  & \(4.30 \times 10^{-5}\)  & \(4.29 \times 10^{-5}\)  & \(1.29 \times 10^{-4}\)  & \(1.29 \times 10^{-4}\)  \\
                0.001953125& \(3.15 \times 10^{-6}\) & \(3.14 \times 10^{-6}\)  & \(1.08 \times 10^{-5}\)  & \(1.07 \times 10^{-5}\)  & \(3.23 \times 10^{-5}\)  & \(3.22 \times 10^{-5}\)  \\
                0.000976563& \(8.32 \times 10^{-7}\) & \(7.81 \times 10^{-7}\)  & \(2.68 \times 10^{-6}\)  & \(2.68 \times 10^{-6}\)  & \(8.07 \times 10^{-6}\)  & \(8.07 \times 10^{-6}\)  \\
			  \hline
		  \end{tabular}}
	\end{center}
    }
\end{table}
\begin{table}[!t]
    \add{
	\begin{center}\renewcommand{\arraystretch}{1.5}
		\caption{RMS errors for different {time step ($\Delta t$)} in Example 3, utilizing a fine node size $\Delta x = 0.0001$.}\label{tab:8}
		\resizebox{0.75\textwidth}{!}
		{
		\begin{tabular}{|c|c|c|c|c|c|} \hline
			\multicolumn{2}{|c|}{ $Re = 10$ at $T = 0.2$} & \multicolumn{2}{c|}{$Re = 25$ at $T = 0.3$} & \multicolumn{2}{c|}{$Re = 50$ at $T = 0.4$} \\ \hline
			$\Delta t$ & RMS errors & $\Delta t$ & RMS errors & $\Delta t$ & RMS errors \\ 
			\hline
			0.2      & 2.49$\times 10^{-2}$  & 0.3      & 3.80$\times 10^{-2}$  & 0.4      & 7.94$\times 10^{-2}$  \\
			0.1      & 1.15$\times 10^{-2}$  & 0.15     & 1.85$\times 10^{-2}$  & 0.2      & 3.88$\times 10^{-2}$  \\
			0.05     & 1.75$\times 10^{-3}$  & 0.075    & 4.64$\times 10^{-3}$  & 0.1      & 1.17$\times 10^{-2}$  \\
			0.025    & 3.54$\times 10^{-4}$  & 0.0375   & 1.02$\times 10^{-3}$  & 0.05     & 2.37$\times 10^{-3}$  \\
			0.0125   & 8.05$\times 10^{-5}$  & 0.01875  & 2.39$\times 10^{-4}$  & 0.025    & 5.38$\times 10^{-4}$  \\
			0.00625  & 1.89$\times 10^{-5}$  & 0.009375 & 5.66$\times 10^{-5}$  & 0.0125   & 1.24$\times 10^{-4}$  \\
			0.003125 & 4.14$\times 10^{-6}$  & 0.0046875 & 1.23$\times 10^{-5}$ & 0.00625  & 2.56$\times 10^{-5}$  \\
			0.0015625 & 4.82$\times 10^{-7}$ & 0.00234375 & 1.61$\times 10^{-6}$ & 0.003125 & 4.65$\times 10^{-6}$  \\
			\hline
		\end{tabular}
		}
	\end{center}
    }
\end{table}
\add{Subsequently, MCCNIM scheme has been tested for higher Reynolds numbers ($Re=10^2 - 10^6$) and numerical solutions for Example 3 obtained at different time ($T=0-2$) using a fixed grid and time steps ($\Delta x=0.0125$, $\Delta t=0.001$) are presented in \fig\ref{fig:14}.  It is evident from the figure that the scheme remains stable and operates smoothly even at large Reynolds numbers, without introducing spurious oscillations near the discontinuous region ($x=0.5$). This highlights the robustness of the developed MCCNIM scheme in handling highly nonlinear problems characterized by shock-like structures.}

\add{The accuracy and convergence properties of the MCCNIM scheme are further analyzed through RMS error computations for various Reynolds numbers ($Re$), using different spatial grid resolutions (\tab\ref{tab:7}) and time step sizes (\tab\ref{tab:8}), in comparison with the MNIM scheme. By closely following the approach in Example 2, the methodology ensures consistency in the comparative analysis, and again demonstrate excellent convergence of MCCNIM scheme, as the RMS error decreases with refinement in grid ($\Delta x$) and time step ($\Delta t$). This trend holds regardless of Reynolds number ($Re$) values, consistent with the MNIM scheme, thereby confirming the accuracy and convergence of the MCCNIM scheme.	
Further, the RMS errors (\tabs\ref{tab:7} and \ref{tab:8}) are analyzed to establish the rate of convergence (ROC) to obtain the convergence order of the MCCNIM scheme, as discussed in Example 2.  The variations of RMS error with nodal spacing ($\Delta x$) and time step ($\Delta t$) are illustrated to be linear (slope of approxminately $\approx 2$)  for a range of Reynolds number ($Re$) in panels (a) and (b) of \fig\ref{fig:15}, respectively.
This demonstration further emphasizes that the developed MCCNIM scheme has quadratic convergence, i.e., second-order accurate in both space and time, O[($\Delta x$)$^2$, ($\Delta t$)$^2$], even for higher Reynolds numbers ($Re$) and for complex shock-like problem, as indicated in panels (a) and (b) of \fig\ref{fig:15} with spatial rate of convergence (SROC $\approx 2$) and temporal rate of convergence (TROC $\approx 2$). 
This consistency highlights the robustness of MCCNIM in handling nonlinear problems while preserving accuracy and computational efficiency.}
\begin{table}[!t]
    \add{
    \begin{center}\renewcommand{\arraystretch}{1.5}
		\caption{Performance evaluation of the proposed MCCNIM scheme for Example 3 at various Reynolds numbers ($Re$) using $T=2$ and $\Delta t=0.1$.}\label{tab:9}
		\resizebox{0.8\textwidth}{!}
		{
		\begin{tabular}{|c|c|c|c|c|c|c|}
        \hline
        \multirow{2}{*}{$\Delta x$} & \multicolumn{2}{c|}{Picard iterations (Nonlinear)} & \multicolumn{2}{c|}{Krylov iterations (Linear)} & \multicolumn{2}{c|}{CPU Runtime (seconds)} \\ \cline{2-7}
        & MNIM & MCCNIM & MNIM & MCCNIM & MNIM & MCCNIM \\
        \hline
        \multicolumn{7}{|c|}{Re = 50} \\
        \hline
        0.1      & 106  & 122  & 1112   & 1521   & 0.008846  & 0.015411  \\
        0.05     & 116  & 128  & 2471   & 2043   & 0.020811  & 0.022013  \\
        0.025    & 124  & 141  & 5178   & 3557   & 0.054423  & 0.045158  \\
        0.0125   & 130  & 146  & 10631  & 6160   & 0.182759  & 0.100295  \\
        0.00625  & 133  & 152  & 22004  & 11031  & 0.788949  & 0.302347  \\
        0.003125 & 137  & 154  & 47141  & 19046  & 5.640219  & 1.104069  \\
        0.0015625& 192  & 157  & 146416 & 32800  & 63.497402 & 5.275467  \\
        \hline
        \multicolumn{7}{|c|}{Re = $10^3$} \\
        \hline
        0.05     & 120  & 139  & 1878   & 2789   & 0.016389  & 0.030211  \\
        0.025    & 120  & 135  & 3029   & 2821   & 0.030432  & 0.034988  \\
        0.0125   & 136  & 155  & 5111   & 3809   & 0.070571  & 0.058586  \\
        0.00625  & 165  & 186  & 8989   & 5713   & 0.202839  & 0.113486  \\
        0.003125 & 195  & 214  & 18189  & 8995   & 0.983258  & 0.309982  \\
        0.0015625& 217  & 255  & 39350  & 16058  & 6.195129  & 1.263147  \\
        0.00078125& 236 & 331  & 92539  & 36922  & 60.230415 & 16.393041 \\
        \hline
        \multicolumn{7}{|c|}{Re = $10^6$} \\
        \hline
        0.05     & 118  & 134  & 2048   & 2938   & 0.01729   & 0.03269   \\
        0.025    & 122  & 135  & 4259   & 4562   & 0.042169  & 0.058966  \\
        0.0125   & 124  & 139  & 8444   & 5897   & 0.134861  & 0.093245  \\
        0.00625  & 128  & 143  & 16433  & 7527   & 0.528363  & 0.173682  \\
        0.003125 & 132  & 147  & 33263  & 11235  & 3.287563  & 0.523028  \\
        0.0015625& 134  & 150  & 67057  & 19069  & 21.581253 & 2.328984  \\
        0.00078125& 137 & 153  & 136197 & 34715  & 214.66407 & 28.528517 \\
        \hline
    \end{tabular}
		}
	\end{center}
    }
\end{table}

\add{Furthermore, Table \ref{tab:9} presents the computational perofmance, of the proposed MCCNIM, in terms of the Picard iterations (for non-linear solver), Krylov iterations (for linear solver), and CPU runtime for Example 3 for various Reynolds numbers ($Re = 50, 10^3, 10^6$), with $T=2$ and $\Delta t=0.1$. For comparison purpose, the results for MNIM scheme have also been obtained and included in \tab\ref{tab:9}. These results follow similar trends to Example 2, demonstrating that MCCNIM generally requires slightly more Picard iterations but significantly reduces Krylov iterations, particularly at finer spatial resolutions. This reduction translates into improved computational efficiency, with MCCNIM achieving lower CPU runtime for smaller $\Delta x$. The findings reaffirm superior performance of MCCNIM in handling nonlinear systems while maintaining stability and accuracy.}
\subsubsection{Verification of Energy Conservation in MCCNIM for One-Dimensional Burgers' Equation}
%
\add{Energy conservation is an essential characteristics in numerical simulations of nonlinear partial differential equations, such as the Burgers equation \citep{anguelov2008energy}. In the inviscid case, total kinetic energy remains constant unless shock formation introduces discontinuities. In contrast, viscous cases exhibit energy dissipation over time due to diffusive effects, as observed in the viscous Burgers' equation.  An accurate numerical scheme must capture this energy behavior to ensure physical reliability.  In particular, numerical dissipation should align with theoretical dissipation to model energy decay correctly. Several studies \citep{anguelov2008energy,murugan2020suppressing,evans2022partial,kundu2024fluid} have investigated energy conservation and dissipation in different numerical methods for Burgers' equation . This study also examines the energy behavior of the developed MCCNIM scheme for Example 3, comparing numerical dissipation with expected theoretical dissipation to assess accuracy.}
\begin{figure}[!b]%
	\centering
	\subfigure[$\Delta x = 0.05$, $\Delta t = 0.05$] {\includegraphics[width=0.48\linewidth]{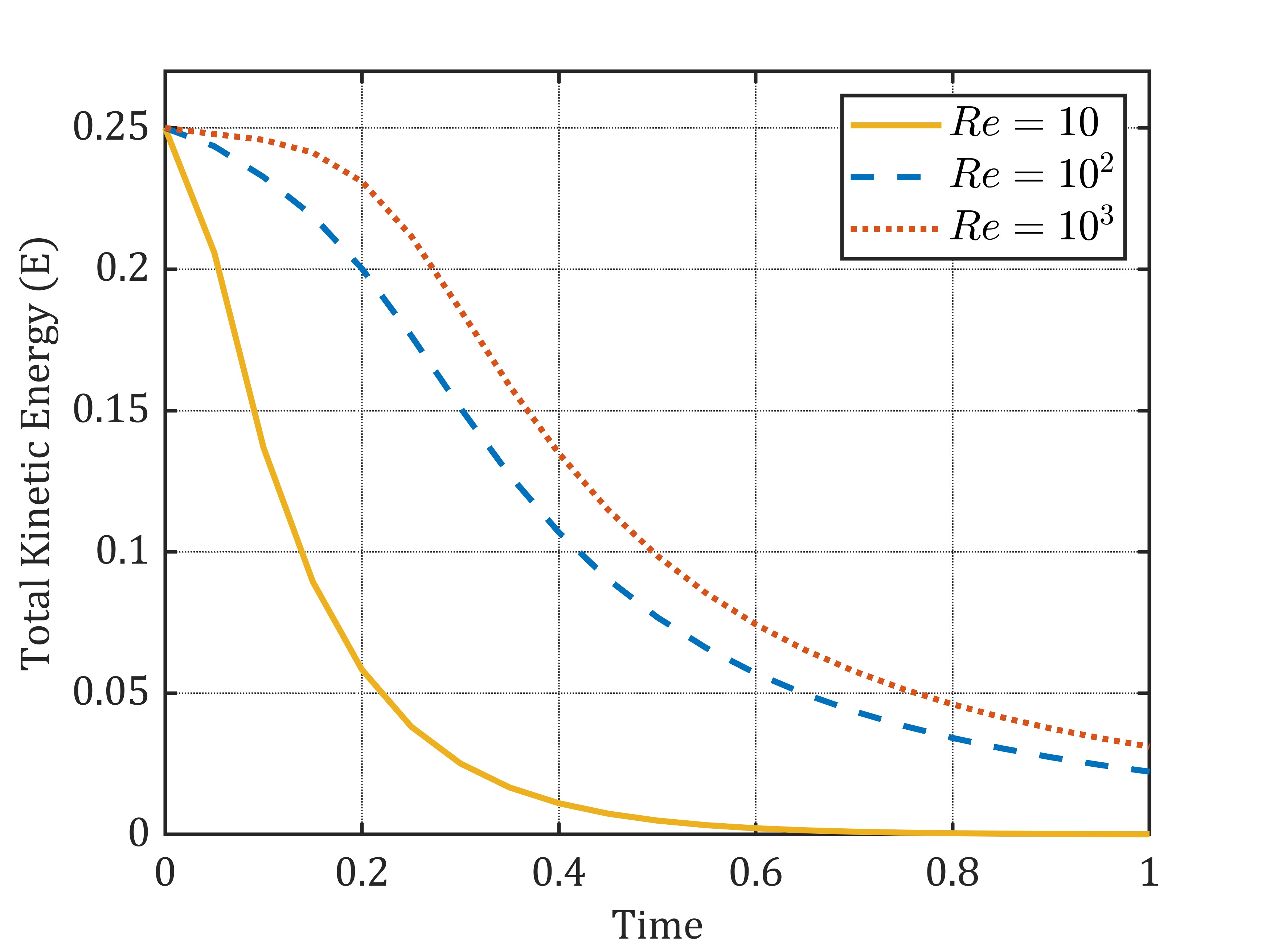}\label{fig17a}}
	\subfigure[$\Delta x = 0.01$, $\Delta t = 0.01$] {\includegraphics[width=0.48\linewidth]{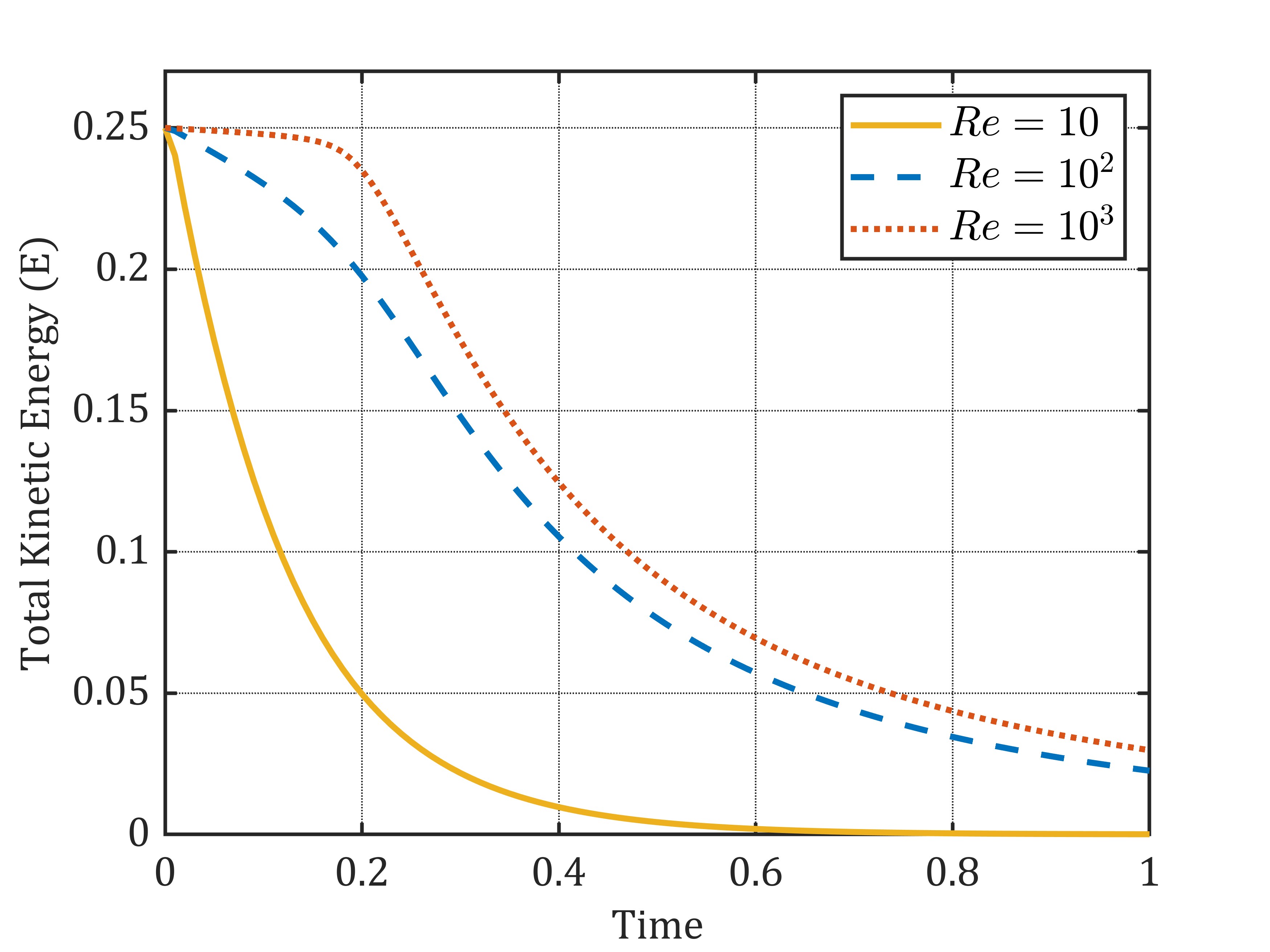}\label{fig17b}}
	\caption{\add{Total kinetic energy ($E$) at $T=1$ as  a function of time ($t$) and Reynolds number ($Re$)  for Example 3.}}%
	\label{fig:17new}%
\end{figure}

\add{The kinetic energy equation is obtained by multiplying the Burgers' equation (\eqn\ref{eq:001}) by the velocity field ($u$) as follows.}
\begin{equation}
    \add{u \frac{\partial u}{\partial t} + u^2 \frac{\partial u}{\partial x} = \frac{1}{Re} u \frac{\partial^2 u}{\partial x^2}.}
\end{equation}
\add{Using the chain rule,}
\begin{equation}
 \add{\frac{\partial E}{\partial t} + \frac{\partial}{\partial x} \left( \frac{1}{2} u^3 \right) = \frac{1}{Re} \left(u \frac{\partial^2 u}{\partial x^2}\right),\qquad E = \frac{1}{2} u^2}
 \label{eq:109}
\end{equation}
\add{The first term represents the rate of change of kinetic energy ($E$), while the second term corresponds to the convective flux, which does not contribute to energy dissipation.}
\begin{figure}[!b]%
	\centering
	\subfigure[$Re = 10$, $\Delta x = 0.02$, $\Delta t = 0.001$] {\includegraphics[width=0.48\linewidth]{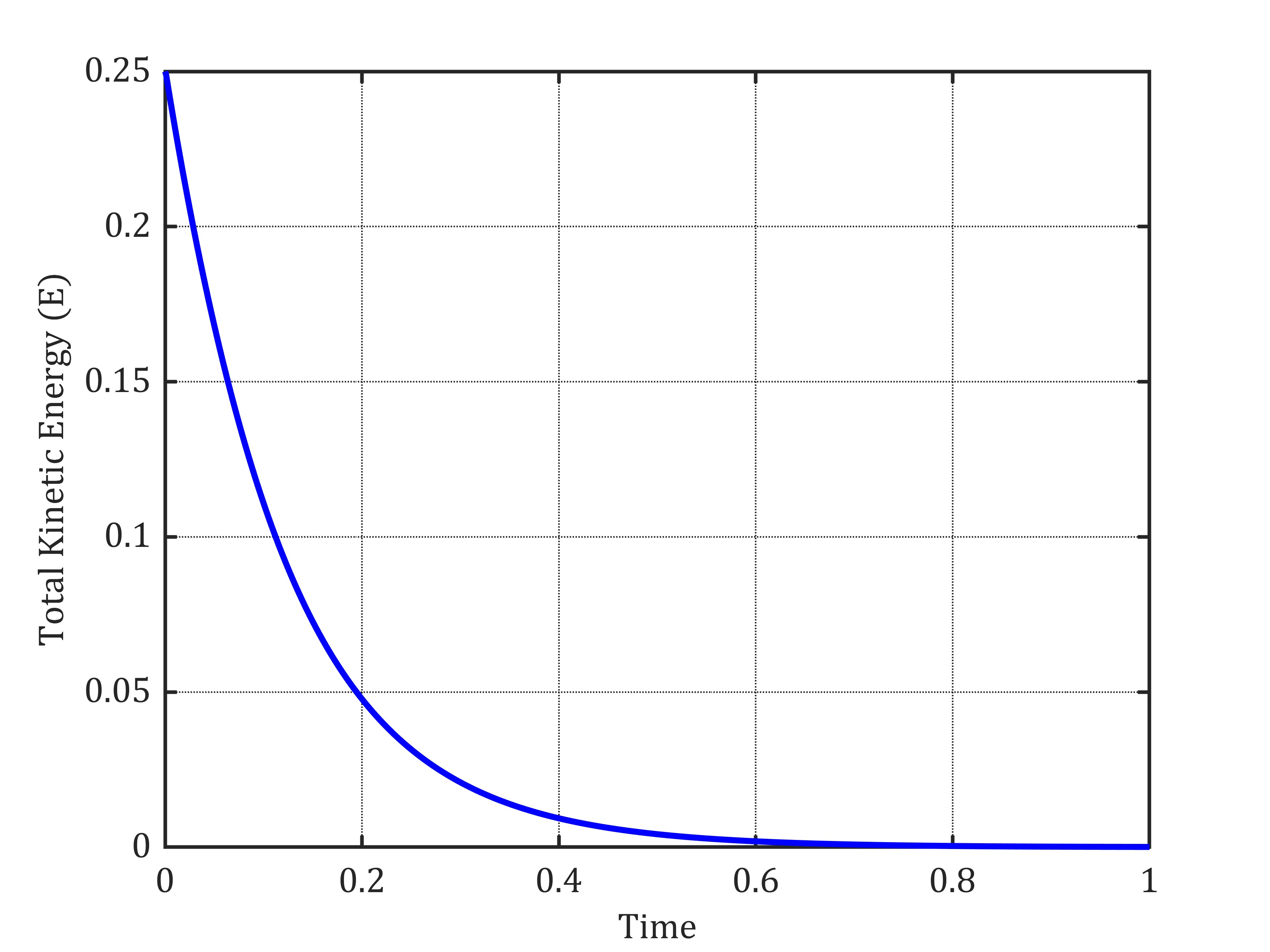}\label{fig18a}}
	\subfigure[$Re = 10$, $\Delta x = 0.02$, $\Delta t = 0.001$] {\includegraphics[width=0.48\linewidth]{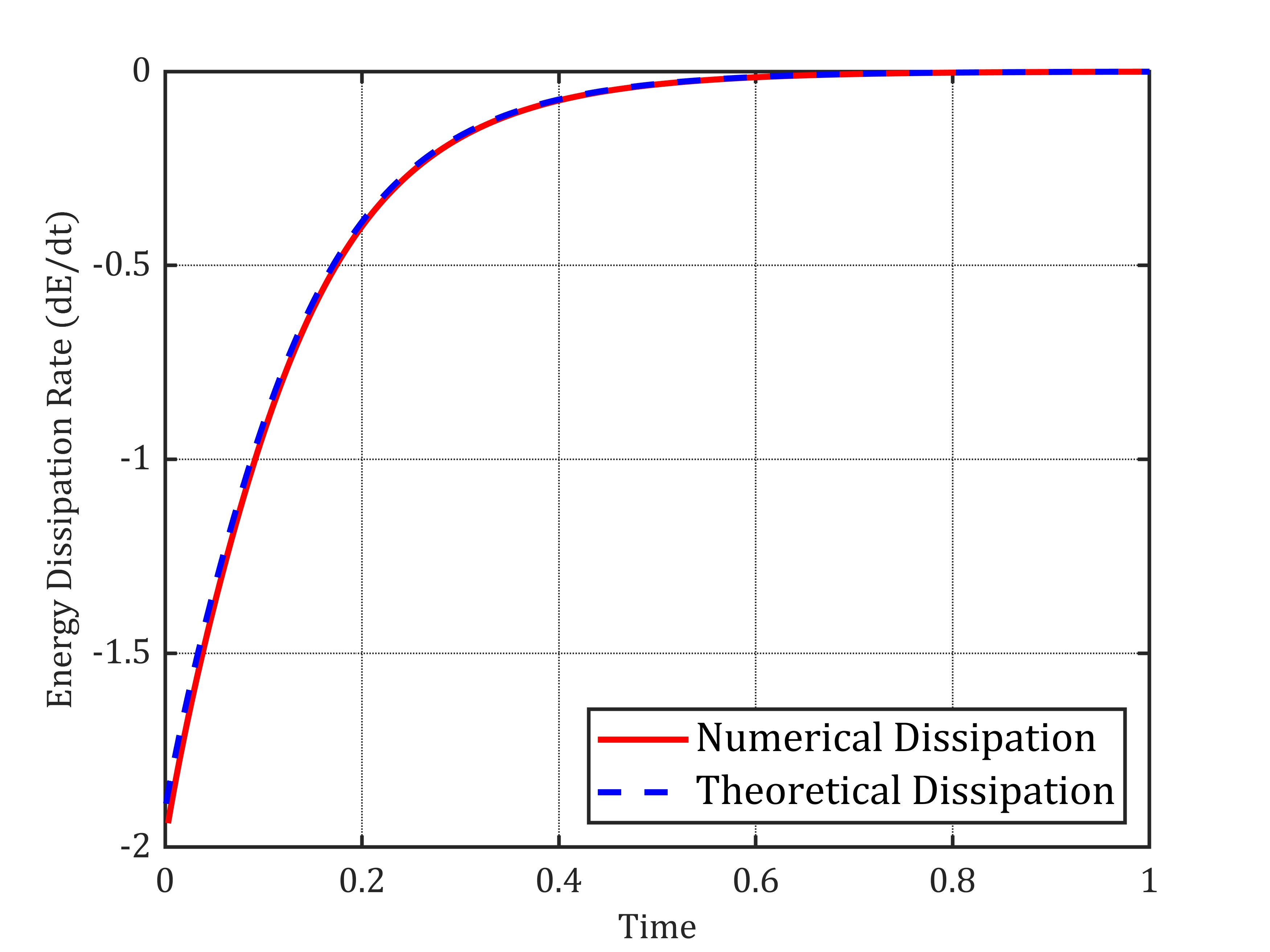}\label{fig18b}} \\
%
%
	\centering
	\subfigure[$Re = 10^2$, $\Delta x = 0.002$, $\Delta t = 0.001$] {\includegraphics[width=0.48\linewidth]{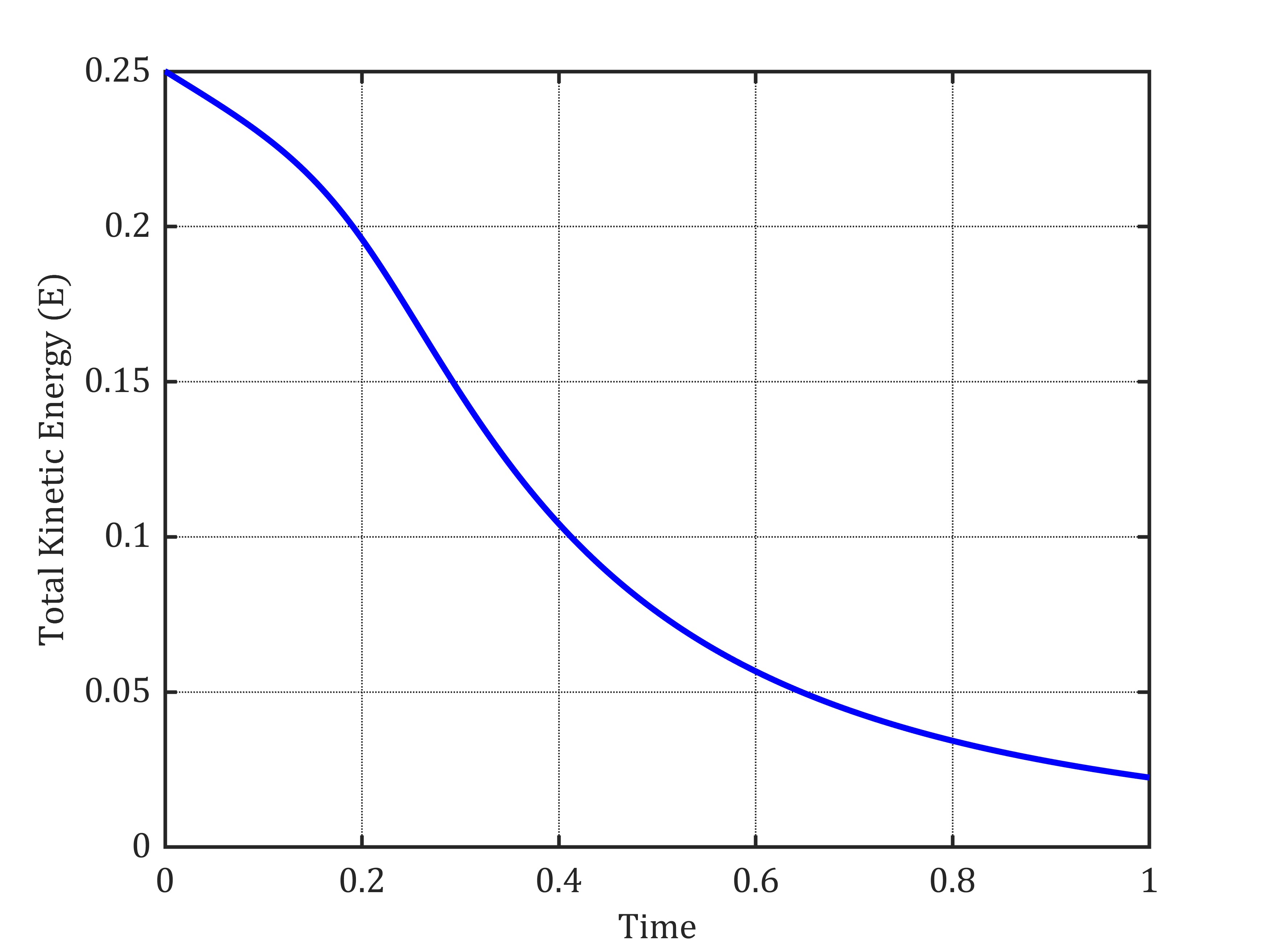}\label{fig19a}}
	\subfigure[$Re = 10^2$, $\Delta x = 0.002$, $\Delta t = 0.001$] {\includegraphics[width=0.48\linewidth]{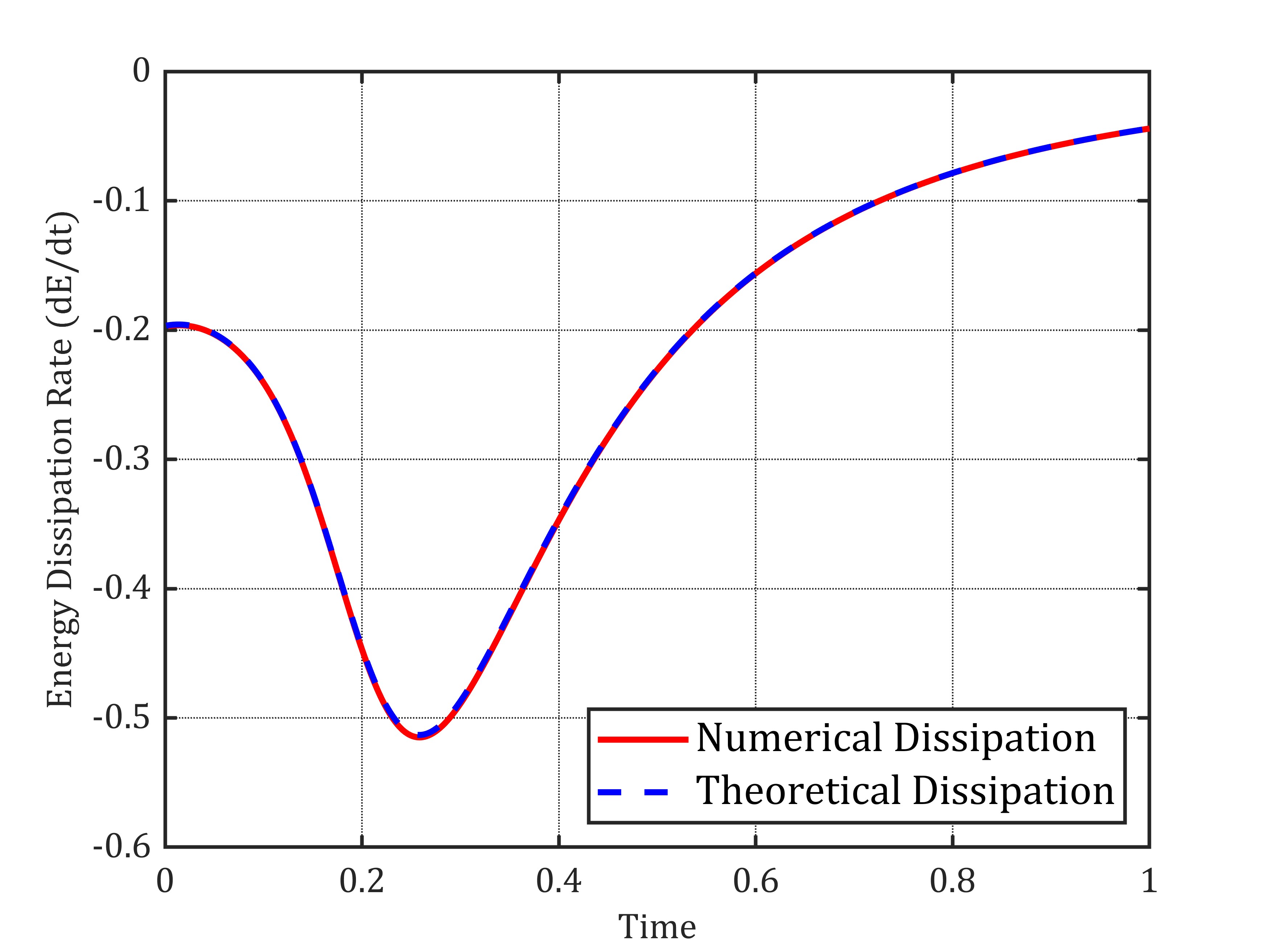}\label{fig19b}} \\
%
%
	\centering
	\subfigure[$Re = 10^3$, $\Delta x = 0.0002$, $\Delta t = 0.001$] {\includegraphics[width=0.48\linewidth]{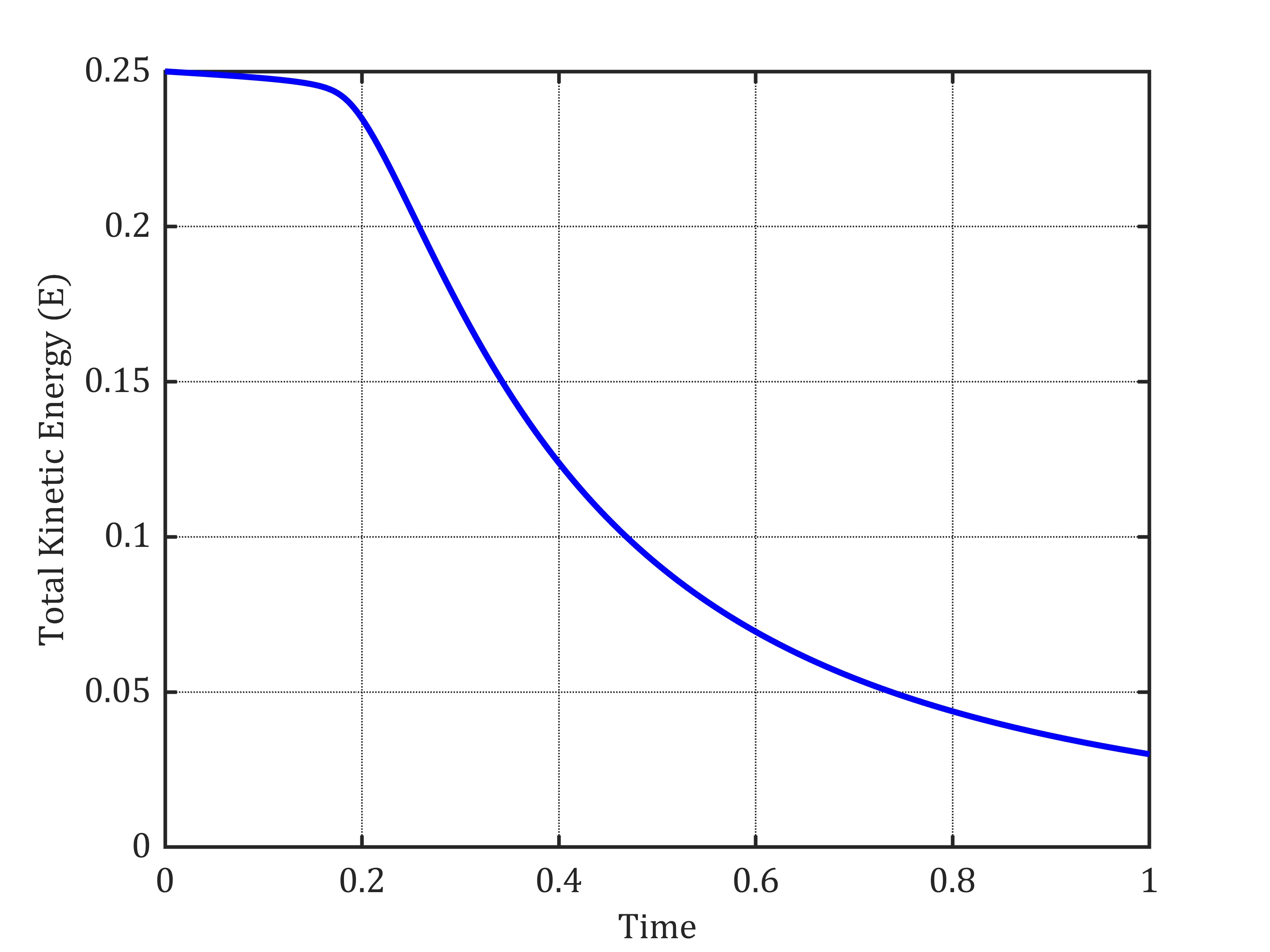}\label{fig20a}}
	\subfigure[$Re = 10^3$, $\Delta x = 0.0002$, $\Delta t = 0.001$] {\includegraphics[width=0.48\linewidth]{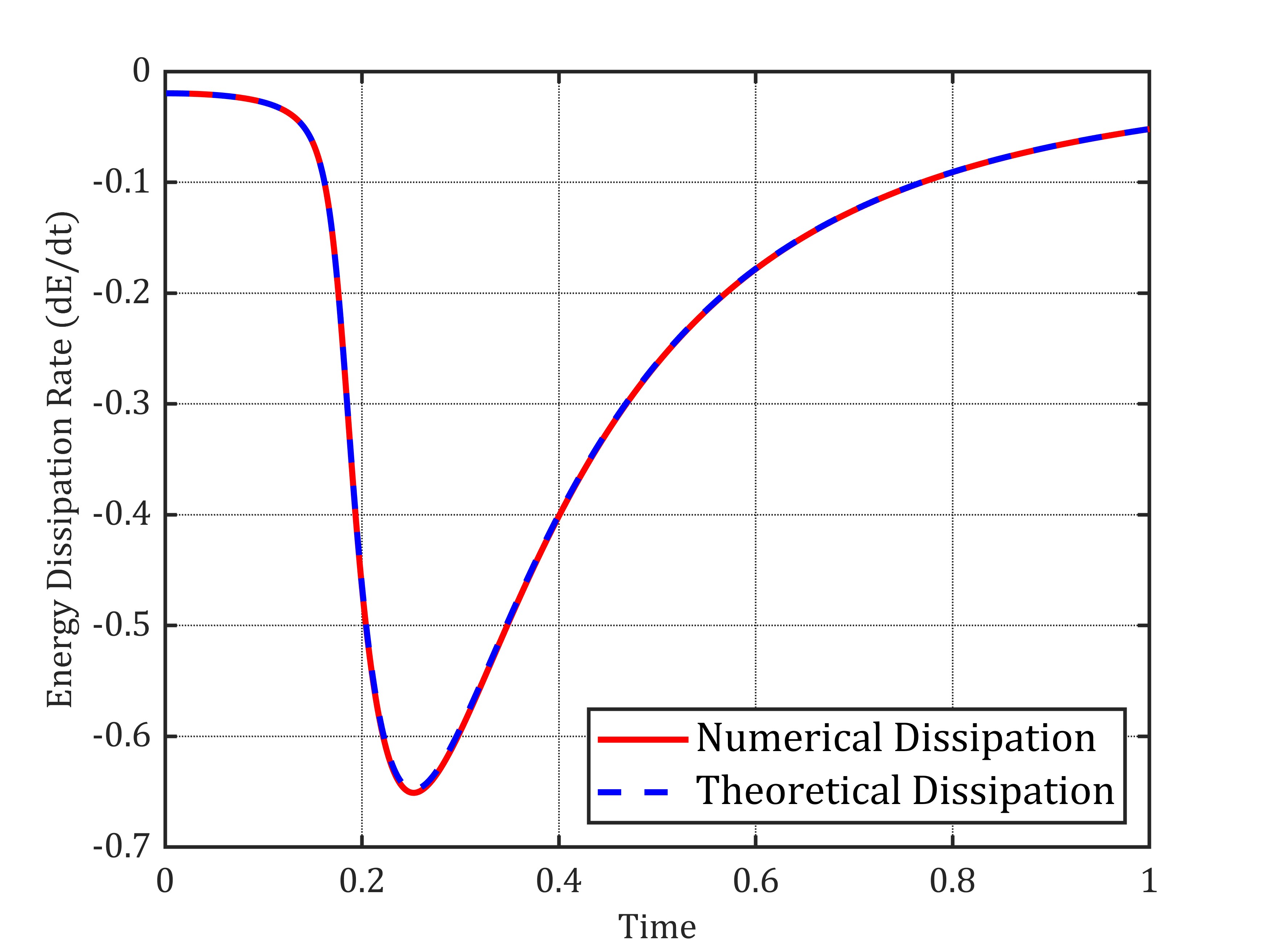}\label{fig20b}}
	\caption{\add{Energy conservation analysis as a function of time ($t$) and Reynolds number ($Re$)  for Example 3: (left) Evolution of total kinetic energy ($E$), and (right) Comparison of numerical and theoretical energy dissipation rate ($\mathrm{d}E/\mathrm{d}t$).}}%
	\label{fig:20new}%
\end{figure}

\add{To analyze energy dissipation further, upon integrating the kinetic energy equation (\eqn\ref{eq:109}) over the spatial domain $x \in [0, L]$, we get}
\begin{equation}
 \add{\int_0^L \frac{\partial E}{\partial t} \mathrm{d}x + \int_0^L \frac{\partial}{\partial x} \left( \frac{1}{2} u^3 \right) \mathrm{d}x = \frac{1}{Re} \int_0^L u \frac{\partial^2 u}{\partial x^2} \mathrm{d}x}
\end{equation}
After simplification, we get
\begin{equation}
	\add{\int_0^L \frac{\partial E}{\partial t} \, \mathrm{d}x + \left[ \frac{1}{2} u^3 \right]_0^L = \frac{1}{Re} \left( \left[ u \frac{\mathrm{d}u}{\mathrm{d}x} \right]_0^L - \int_0^L \left( \frac{\mathrm{d}u}{\mathrm{d}x} \right)^2 \, dx \right)}
\end{equation}
\add{On application of the Dirichlet boundary conditions, as given in \eqn\eqref{eq:096} for Example 3, we obtain the kinetic energy decay equation as follows.}
\begin{equation}
	\add{\int_0^L \frac{\partial E}{\partial t} \, \mathrm{d}x  = -\frac{1}{Re} \int_0^L \left( \frac{\mathrm{d}u}{\mathrm{d}x} \right)^2 \, \mathrm{d}x }
	\qquad\Rightarrow\qquad
		\add{\frac{\mathrm{d} E}{\mathrm{d} t} = -\frac{1}{Re}  \int_0^L \left( \frac{\mathrm{d}u}{\mathrm{d}x} \right)^2 \, \mathrm{d}x }
\end{equation}
\add{Since the right-hand side is always negative, the total kinetic energy ($E$) monotonically decreases due to viscous dissipation.}

\add{For a discrete numerical scheme, the numerical energy dissipation rate is approximated using finite difference as follows:}
\begin{equation}
\add{\frac{dE}{dt} \approx \frac{E_{\ell+1} - E_{\ell}}{\Delta t}}
\label{eq:113}
\end{equation}
\add{where ${\ell}$ represent the current time index, and the total kinetic energy is evaluated as follows.}
\begin{equation}
	\add{E_{\ell} = \frac{1}{2} \sum_{i} u_i^2 \Delta x}
\end{equation}
\add{On the other hand, the theoretical dissipation rate is computed from the discrete version as follows.}
\begin{equation}
\add{\frac{dE}{dt} = - \frac{1}{Re} \sum_i \left( \frac{u_{i+1} - u_{i-1}}{2\Delta x} \right)^2 \Delta x}
\label{eq:115}
\end{equation}
\add{The matching numerical and theoretical dissipation rates (\eqns\ref{eq:113} and \ref{eq:115}) confirm that the numerical scheme accurately captures energy conservation and dissipation.}

\add{To understand the energy  conservation and dissipation characteristics of the proposed MCCNIM scheme,  \figs\ref{fig:17new}  and \ref{fig:20new} present the instantaneous variation of total kinetic energy ($E$) and energy dissipation rate ($\mathrm{d}E/\mathrm{d}t$) for Example 3 for a wide range of Reynolds number ($Re$). As expected, \fig\ref{fig:17new} illustrates the energy decay patterns, which aligns with those reported in the literature \citep{anguelov2008energy} for the viscous Burgers' equation. This demonstrates that the proposed MCCNIM scheme accurately preserves the expected energy behavior, consistent with the reference \citep{anguelov2008energy}.}
\add{Furthermore,  \fig\ref{fig:20new} analyzes the evolution of total kinetic energy ($E$), and numerical and theoretical energy dissipation rate ($\mathrm{d}E/\mathrm{d}t$) as a function of Reynolds number ($Re$). It clearly demonstrates an excellent match between the numerical and theoretical energy dissipation rates. This strongly indicates that MCCNIM preserves the correct energy balance without introducing spurious energy growth or artificial dissipation.}
\subsubsection{Example 4: Two-Dimensional Propagating Shock Problem}
\label{sec:3.1.3}
The last example focuses on a shock wave that propagates at a 45-degree angle to the node orientation. This two-dimensional propagating shock problem \citep{25Wescott_2001} serves as an extension of the previously described one-dimensional case (Example 1), and the exact solution follows a similar pattern, which is represented as follows.
\begin{figure}[!hb]%
	\centering
	\subfigure[~]{\includegraphics[width=0.45\linewidth]{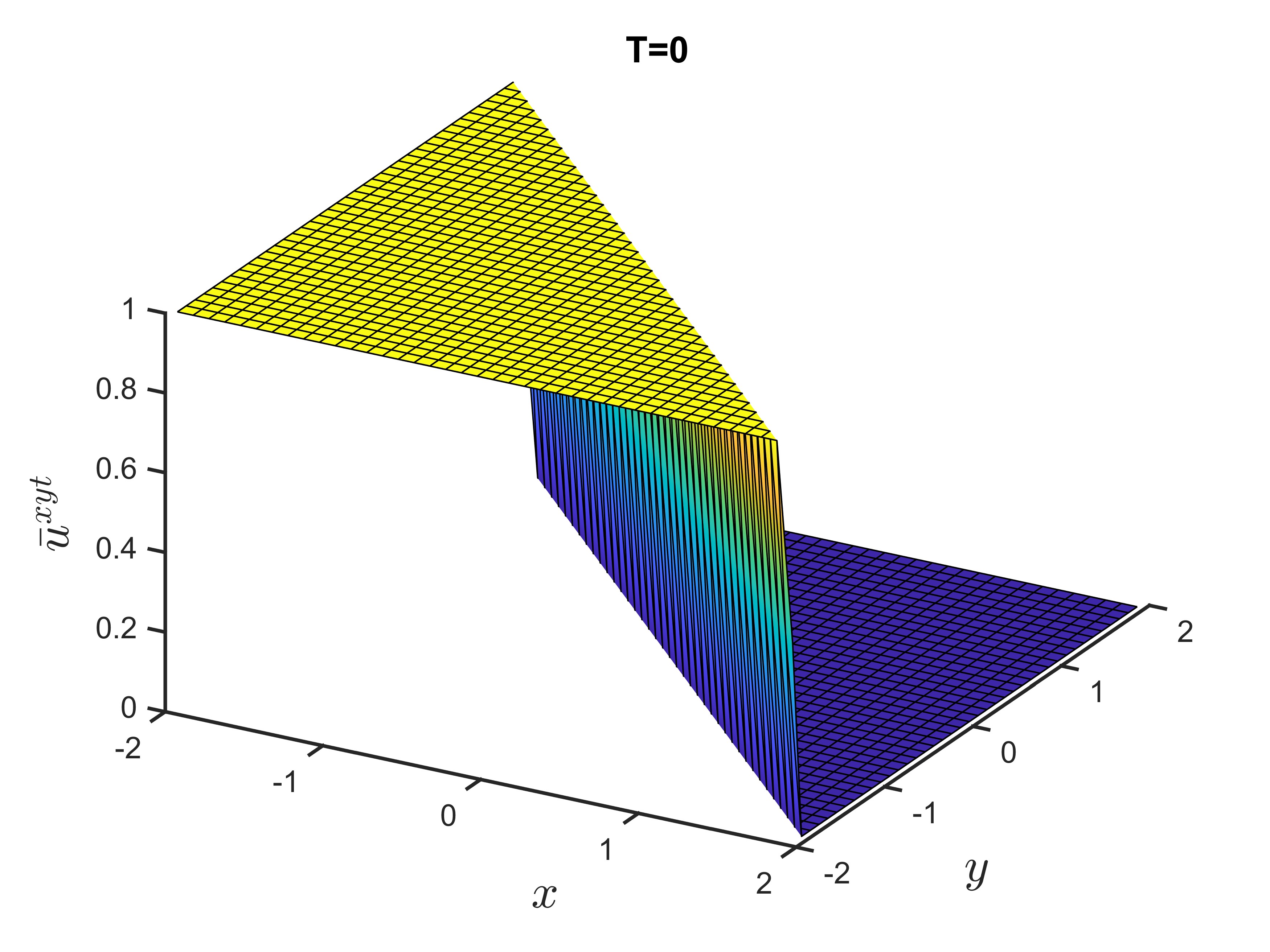}\label{fig10a}}
	\subfigure[~] {\includegraphics[width=0.45\linewidth]{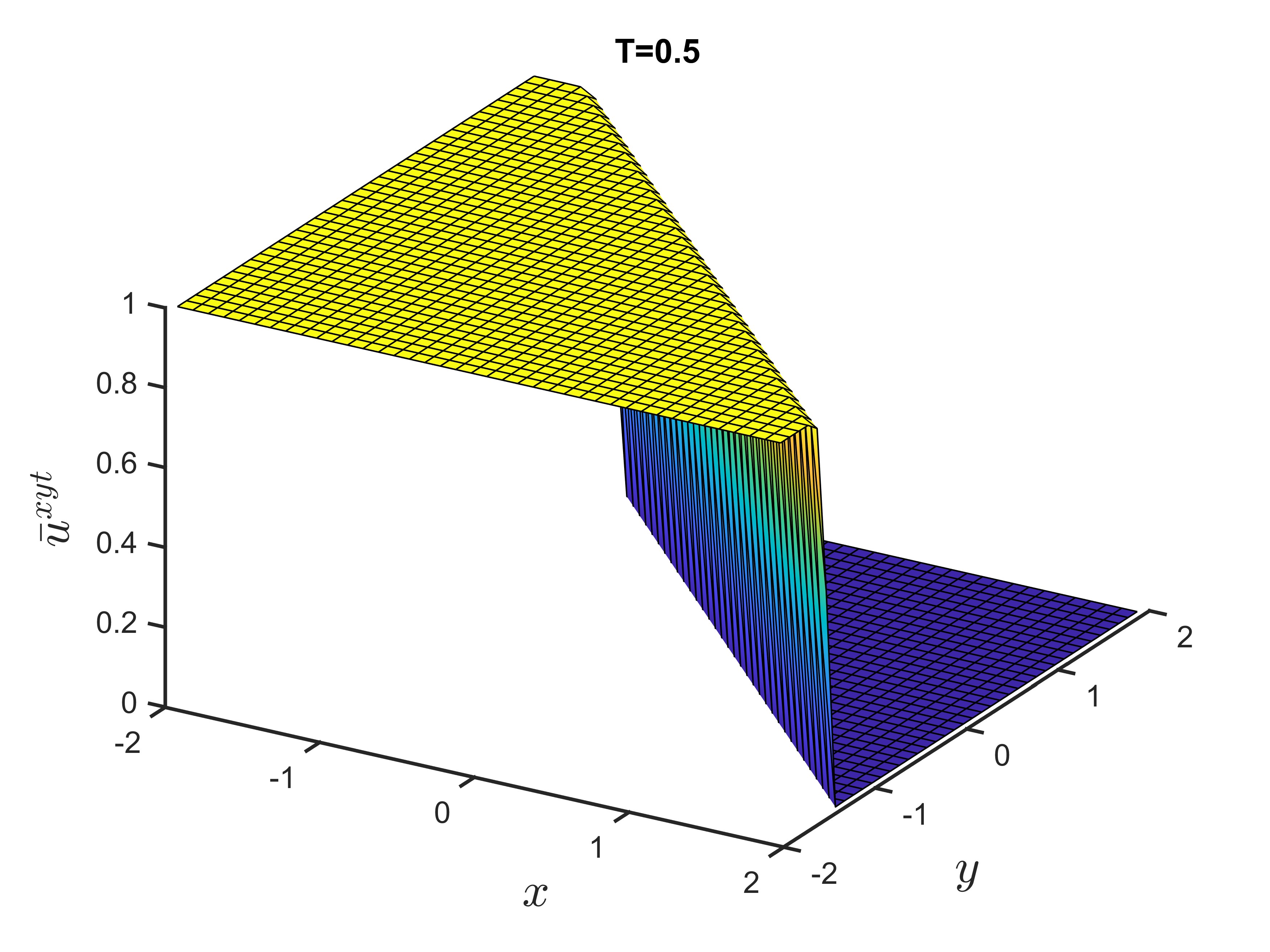}\label{fig10b}}\\
	\subfigure[~]{\includegraphics[width=0.45\linewidth]{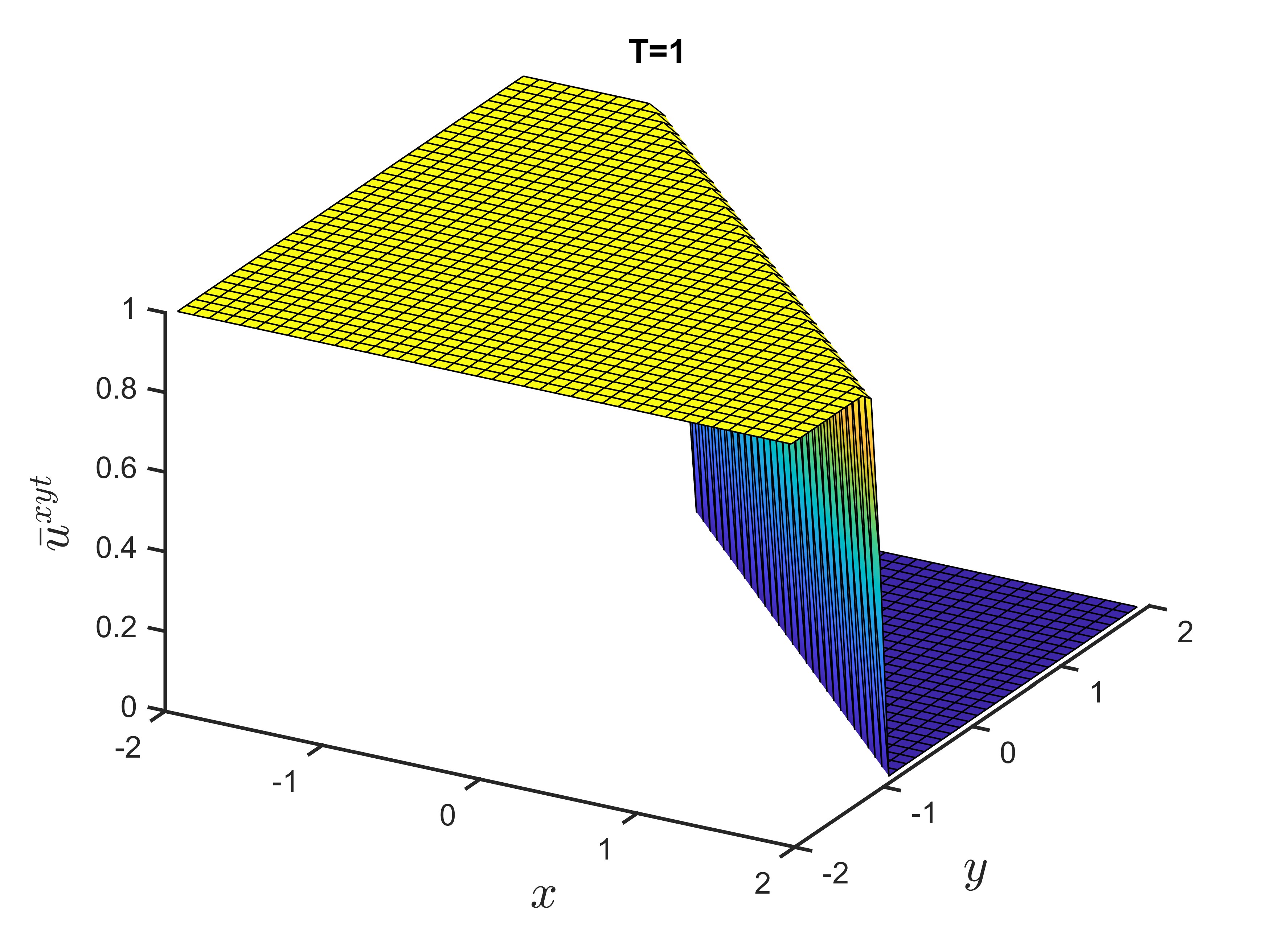}\label{fig10c}}
	\subfigure[~] {\includegraphics[width=0.45\linewidth]{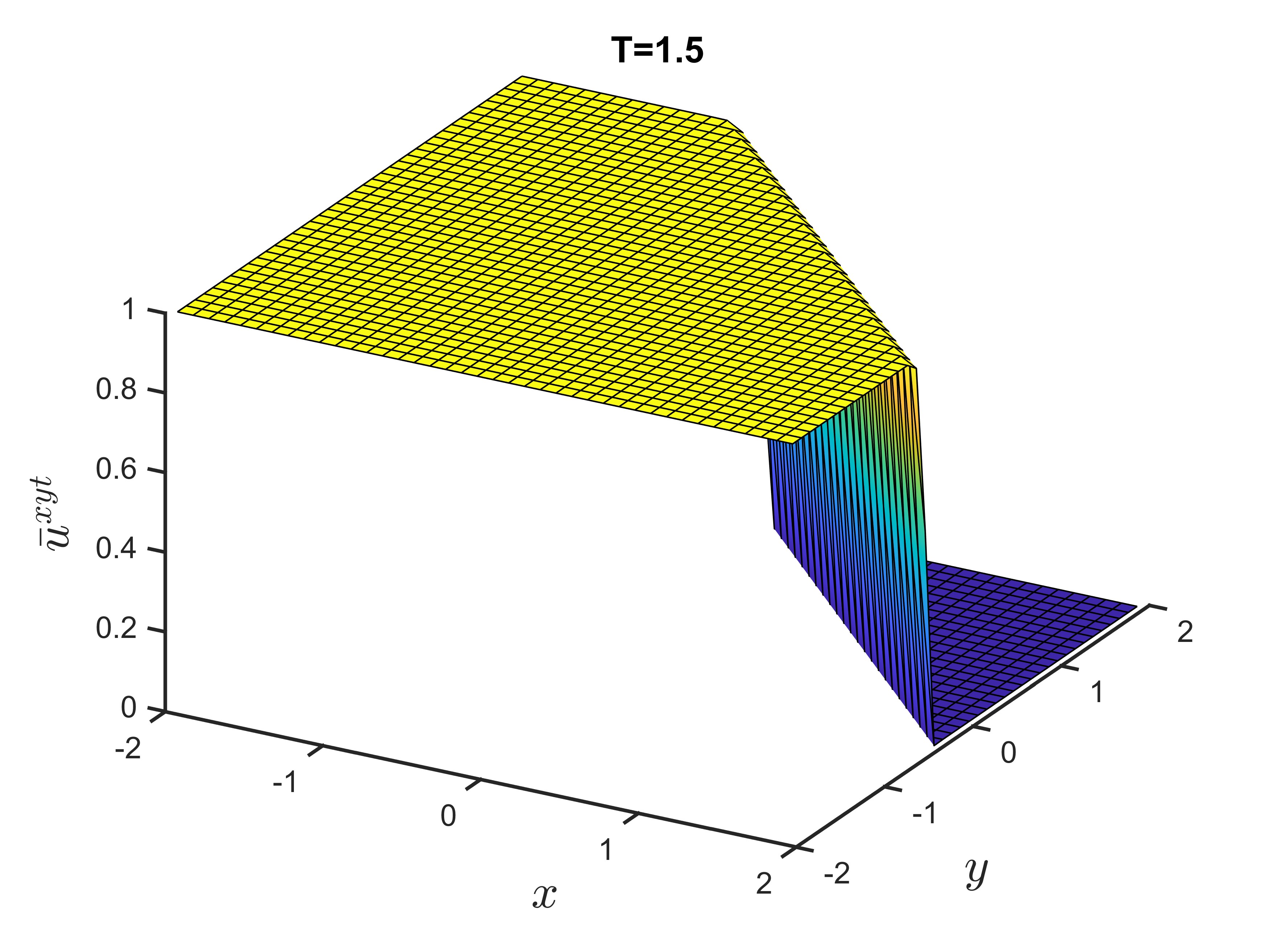}\label{fig10d}}\\
	\subfigure[~]{\includegraphics[width=0.45\linewidth]{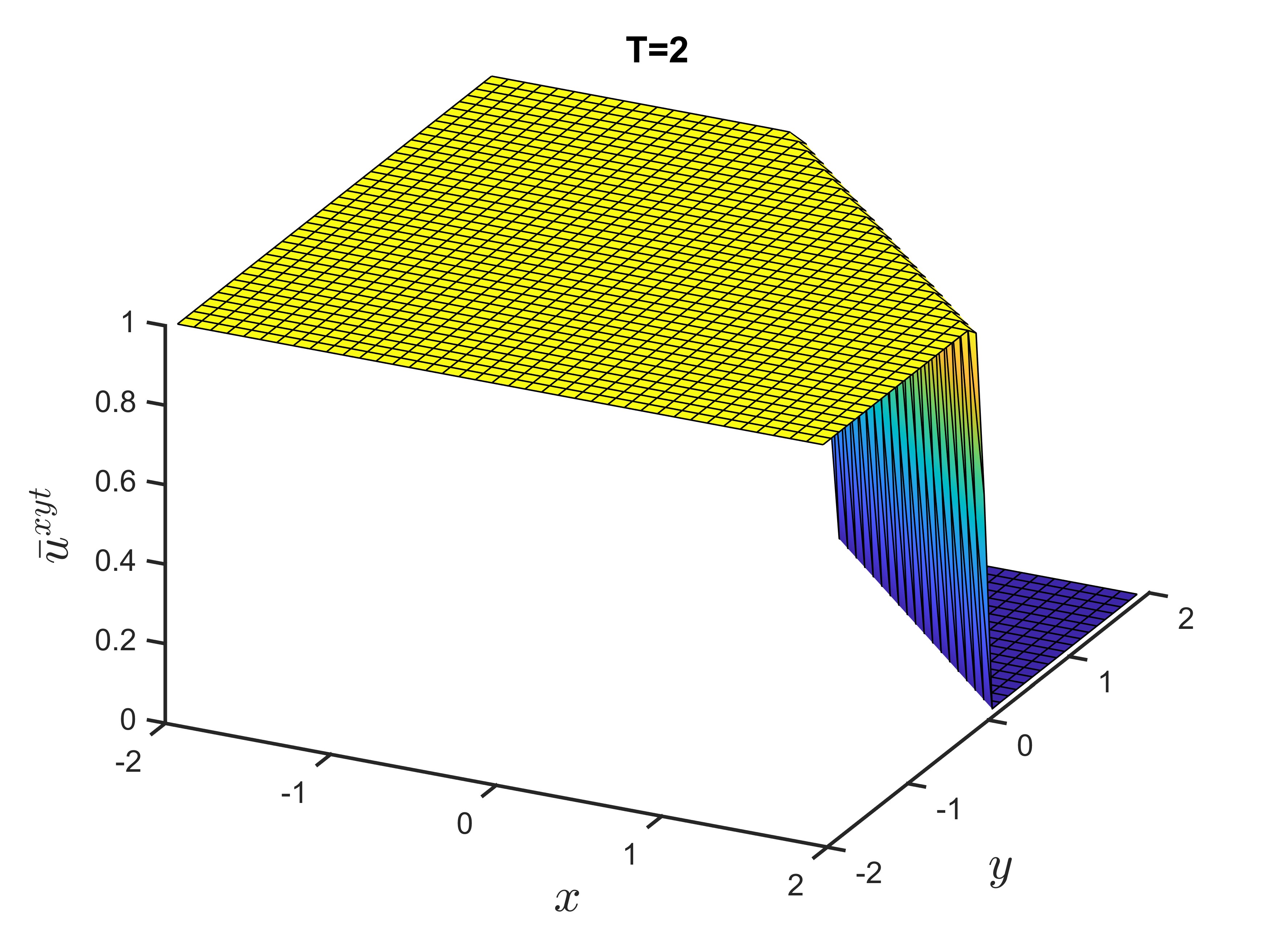}\label{fig10e}}
	\subfigure[~] {\includegraphics[width=0.45\linewidth]{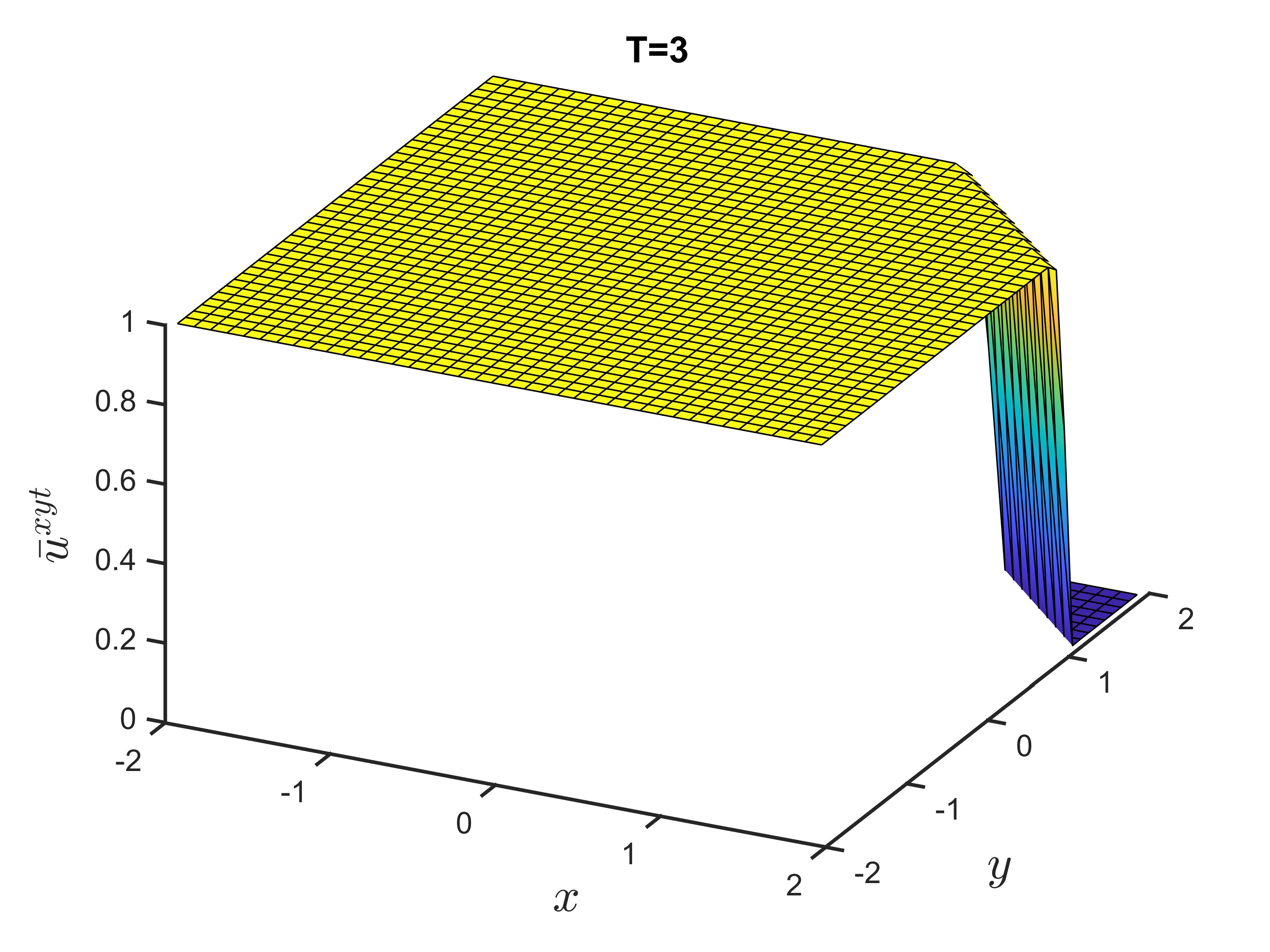}\label{fig10f}}\\
	\caption{Numerical solution for Example 4 using MCCNIM, illustrating the temporal evolution of the shock wave for $Re = 5000$, with a time step $\Delta t=0.01$ and spatial resolution $\Delta x=0.1$.}%
	\label{fig:10}%
\end{figure}
\begin{gather}
	\begin{split}
		u\left(x,y,t\right)=\frac{1}{2}\left[1-\tanh\left(\frac{xRe}{4}+\frac{yRe}{4}-\frac{tRe}{4}\right)\right] 
	\end{split}
	\label{eq:103}
\end{gather}
\begin{gather}
	\begin{split}
		v\left(x,y,t\right)=\frac{1}{2}\left[1-\tanh\left(\frac{xRe}{4}+\frac{yRe}{4}-\frac{tRe}{4}\right)\right] 
	\end{split}
	\label{eq:104}
\end{gather}
The problem domain is restricted to $[-2\ \ 2]\times[-2\ \ 2]$, with initial and boundary conditions determined from the exact solution (\eqns\ref{eq:103} and \ref{eq:104}).
\begin{figure}[!b]%
	\centering
	{\includegraphics[width=0.75\linewidth]{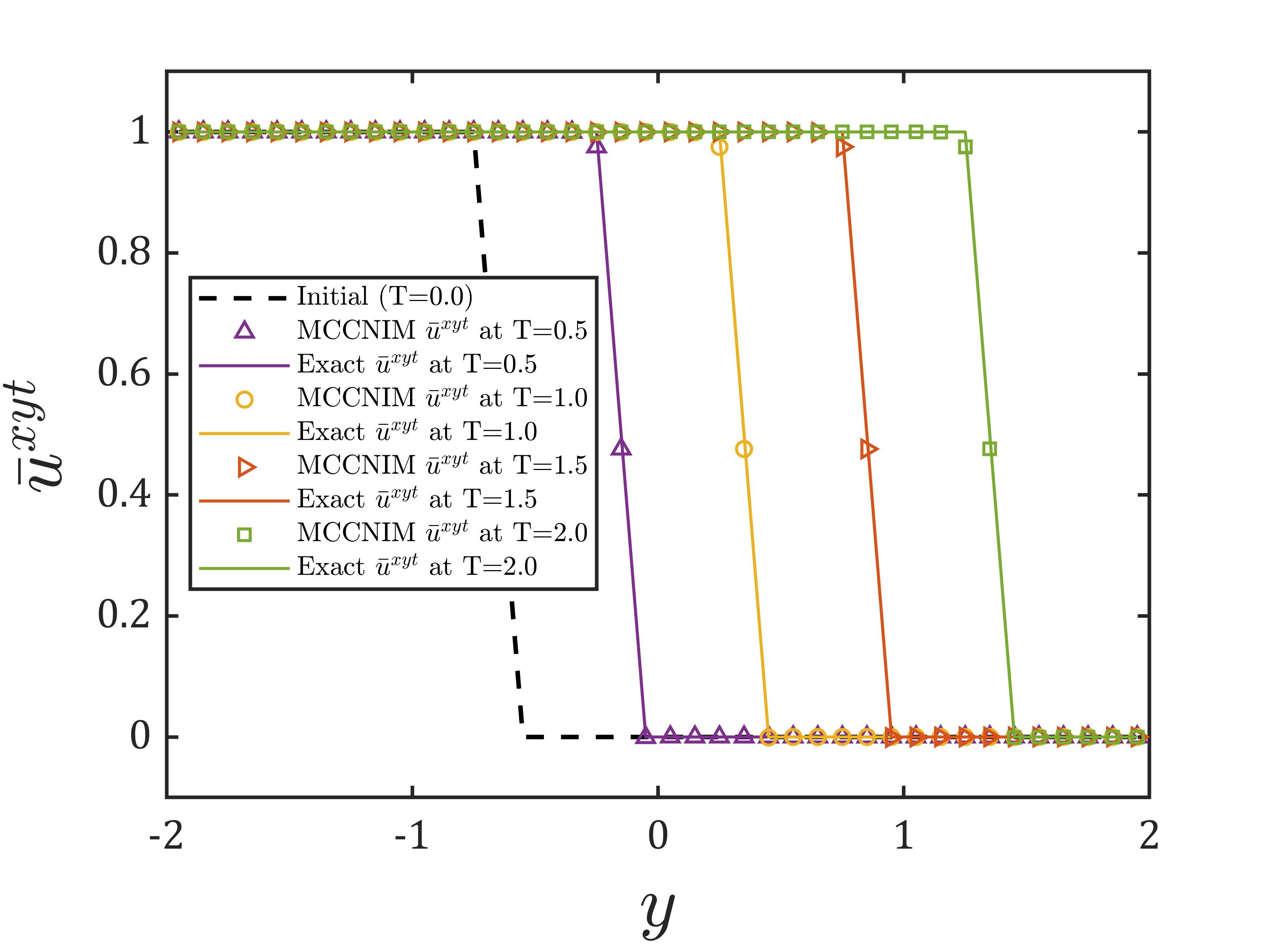}}
	\caption{Comparison of numerical and exact velocity profiles (${\bar{u}}^{xyt}$) for Example 4 at $x=0.7$, evaluated at various time points for $Re\ =\ 5000$, with a time step $\Delta t=0.01$ and spatial resolution $\Delta x=0.1$.}%
	\label{fig:11}%
\end{figure}
\begin{figure}[!b]%
	\centering
	\subfigure[~]{\includegraphics[width=0.48\linewidth]{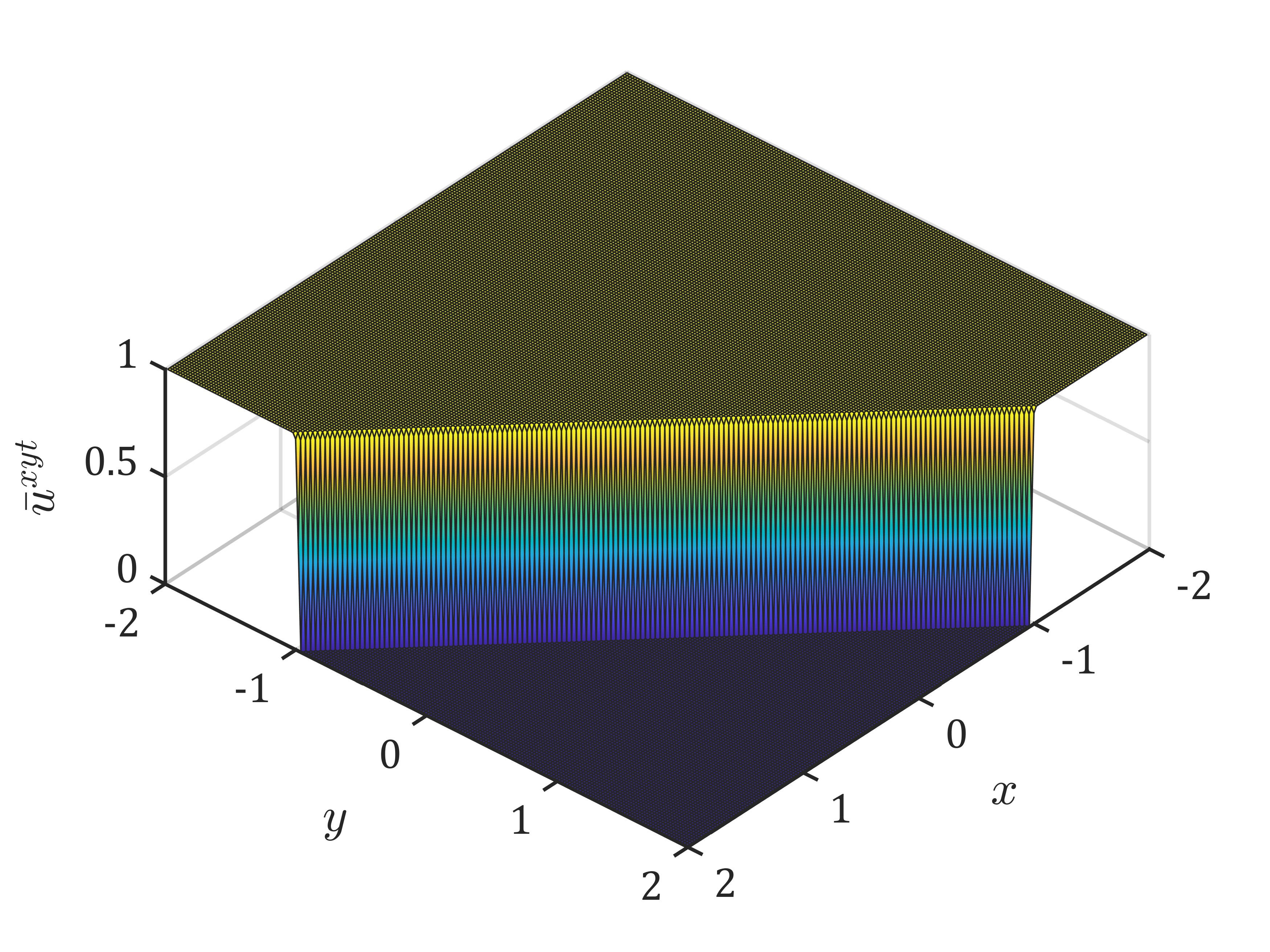}\label{fig12a}}
	\subfigure[~] {\includegraphics[width=0.48\linewidth]{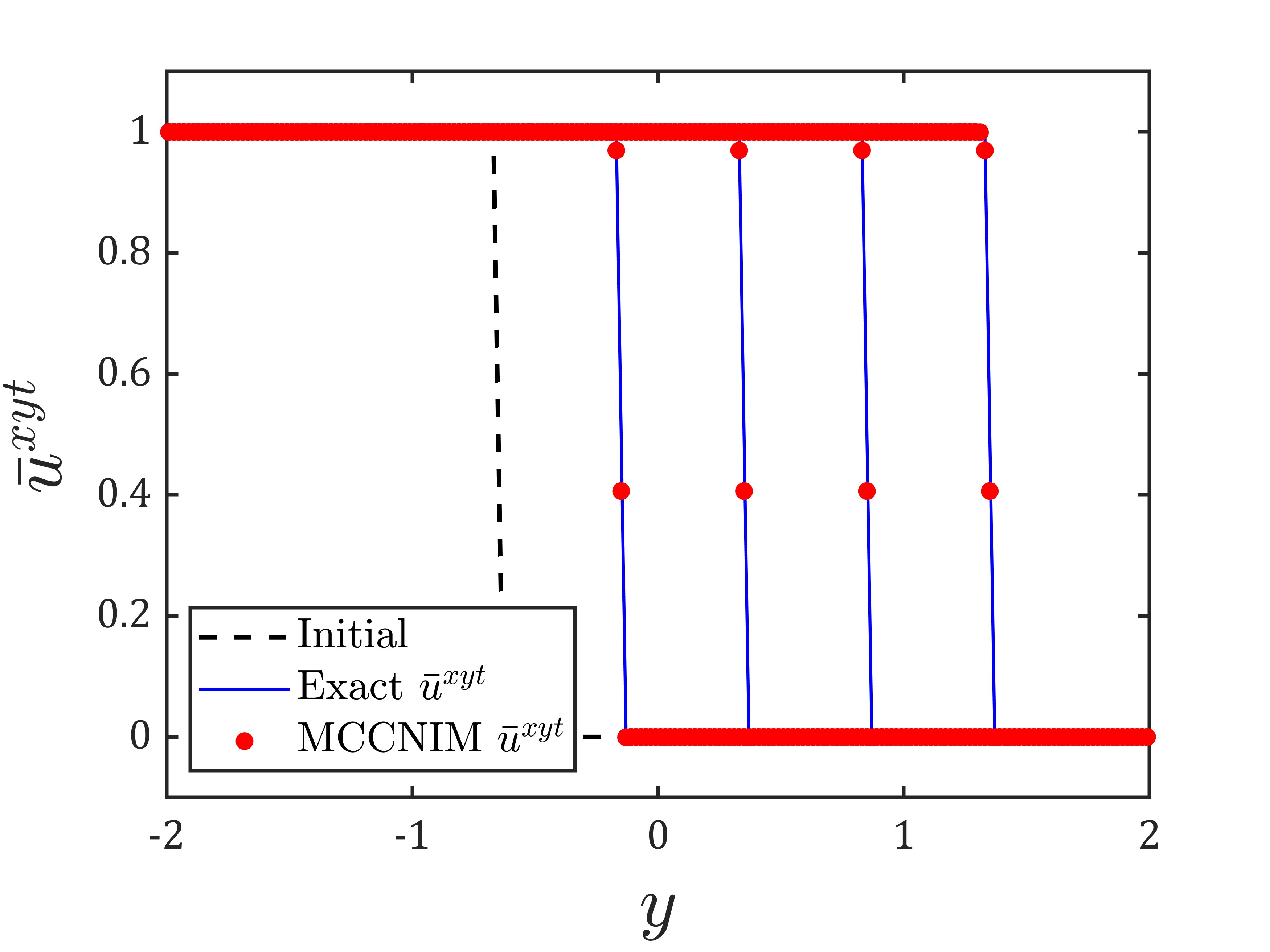}\label{fig12b}}
	\caption{Numerical solution obtained using MCCNIM for $Re=40,000$ (a) Surface plot illustrating the complete solution pattern at $T=1$ with the grid size $\Delta x=0.02$ and time step $\Delta t=0.005$. (b) Velocity profiles (${\bar{u}}^{xyt}$) at $x=0.7$, compared at four distinct time points ($T=0.5$, 1, 1.5, and 2) with analytical solution.}%
	\label{fig:12}%
\end{figure}

\fig\ref{fig:10} presents a 3D surface plot that visually captures the dynamic behavior of the propagating wave for the two-dimensional case at distinct time points ($T=0$ to $3$). These plots provide a clear depiction of the wave evolution over time. Notably, at a high Reynolds number ($Re\ =\ 5000$), the wavefront exhibits remarkable smoothness, even when computed on a very coarse grid ($\Delta x=0.1$). This observation highlights the robustness and accuracy of the developed MCCNIM scheme, validating its effectiveness for solving two-dimensional problems. The progression of the shock wave shown in \fig\ref{fig:10} is further analyzed at each time point using 2D line plots, where the value of $x=0.7$ is held constant, as illustrated in \fig\ref{fig:11}. These plots emphasize the comparison between the numerical and exact solutions. The numerical solution shows excellent agreement with the exact solution, even on coarse grids, despite the challenge of capturing a steep front propagating at an angle to the grid orientation. Notably, \fig\ref{fig:11} demonstrates the excellent performance of the MCCNIM scheme, achieving a high degree of alignment with the exact solution at a high Reynolds number ($Re\ =\ 5000$) using a relatively coarse grid resolution ($40\times 40$).
To further evaluate the efficacy of the developed scheme for a two-dimensional case at higher Reynolds numbers, the numerical simulations are performed for a very high Reynolds number ($Re=40,000$) and corresponding results are illustrated in \fig\ref{fig:12}. The numerical and the exact solutions agree closely (\fig\ref{fig:12}), underscoring the effectiveness of the proposed scheme for two-dimensional case.

It is notable that this two-dimensional problem serves as an extension of the earlier one-dimensional case, for which a detailed error analysis is reported in previous sections. Given the favorable performance observed in the one-dimensional case using the MCCNIM scheme, it is reasonable to anticipate its superior performance in the two-dimensional context as well, surpassing other nodal schemes. Additionally, due to the lack of error data for NIM in existing literature for the two-dimensional case, our comparison relies exclusively on evaluating the numerical solution against the exact solution. Notably, this comparison reveals a strong agreement between the numerical and exact solutions across a range of Reynolds numbers, from lower to higher values. Therefore, this work clearly demonstrates the reliability and efficacy of the proposed MCCNIM scheme in solving the one- and two-dimensional non-linear Burgers' equation with excellent accuracy and computational efficiency.
%
\section{Concluding remarks}
\label{sec:4}
\noindent 
\add{In this study, the MCCNIM scheme, previously limited to linear cases, has been successfully extended to solve non-linear multidimensional Burgers’ equations. Its effectiveness was first demonstrated through a comparative analysis of RMS errors with the traditional NIM for the one-dimensional Burgers’ equation, alongside validation using shock velocity and energy conservation properties. Further tests confirmed the second-order accuracy of the scheme in both spatial and temporal domains. A detailed performance evaluation, including comparisons of iteration counts (Picard and Krylov) and CPU runtime (in seconds) with the traditional NIM, showed that MCCNIM is either comparable or superior in terms of accuracy and computational efficiency. To further establish its robustness, the scheme is applied to two-dimensional coupled Burgers’ equations using classical benchmark problems with well-known analytical solutions. The results demonstrated the accuracy and stability of MCCNIM across both one- and two-dimensional cases.} Comparative analysis with traditional numerical schemes from existing literature highlights the ability of MCCNIM to efficiently and accurately address nonlinear problems. These findings establish a strong foundation for extending MCCNIM to more complex fluid dynamics problems, including applications to the Navier-Stokes equations.
%
\section*{Declaration of Competing Interest}
\noindent 
The authors declare that they have no known competing financial interests or personal relationships that could have appeared to influence the work reported in this paper.
%
\section*{Authors Contributions Statement}
%
\begin{table*}[!h]
\vspace{-1em} 
	\begin{center}\renewcommand{\arraystretch}{1.5}
	\begin{tabular}{|l|p{0.75\linewidth}|}
		\hline
		Author(s)  & Contribution(s) \\\hline
		1. Nadeem Ahmed	& 
		Conceptualization, Methodology, Validation, Formal analysis, Investigation, Data Curation, Visualization, Writing - Original Draft\\\hline
		2. Ram Prakash Bharti	& 
		Supervision, Funding, Conceptualization, Methodology, Software, Validation, Resources, Formal analysis, Writing - Review \& Editing\\\hline
		3. Suneet Singh	& 
		Conceptualization, Methodology, Formal analysis, Writing - Review \& Editing\\\hline
	\end{tabular}
\end{center}
\vspace{-2em} 
\end{table*}
%
\section*{Acknowledgements}
Nadeem Ahmed acknowledges the support received under the Institute Post-Doctoral Fellow (I-PDF) scheme of the Indian Institute of Technology Roorkee, Roorkee, India.
%
%
\begin{spacing}{1.5}
\fontsize{10}{10pt}\selectfont
\renewcommand{\nomgroup}[1]{%
	\ifthenelse{\equal{#1}{0}}{~}{
		\ifthenelse{\equal{#1}{A}}{\item[\textbf{Abbreviations}]}{%
			\ifthenelse{\equal{#1}{G}}{\item[\textbf{Greek Symbols}]}{%
				\ifthenelse{\equal{#1}{D}}{\item[\textbf{Dimensionless Groups / Constants}]}{%
					\ifthenelse{\equal{#1}{S}}{\item[\textbf{Subcripts / Superscripts}]}
					{}}}}}
}
%
\nomenclature[A]{CC}{cell-centered}
\nomenclature[A]{CCNIM}{cell-centered nodal integral method}
\nomenclature[A]{CN-4PU}{Crank-Nicolson 4-point upwind}
\nomenclature[A]{NIM}{nodal integral method}
\nomenclature[A]{MCCNIM}{modified cell-centered nodal integral method}
\nomenclature[A]{M$^2$NIM}{modified-modified nodal integral method}
\nomenclature[A]{CDE}{convection-diffusion equation}
\nomenclature[A]{DAE}{differential-algebraic equation}
\nomenclature[A]{ODE}{ordinary differential equation}
\nomenclature[A]{PDE}{partial differential equation}
\nomenclature[A]{TIP}{transverse integration process}
%
\nomenclature[0a]{$a$}{half of the node size in the $x$-direction}
\nomenclature[0b]{$b$}{half of the node size in the $y$-direction}
\nomenclature[0j]{$\bar{J}$}{surface-averaged flux}	
\nomenclature[0s]{$\bar{S}$}{pseudo-source term}	
\nomenclature[0u]{${\bar{u}}^{xy}$}{$x$ and $y$-space averaged velocity in the $t$-direction}
\nomenclature[0u]{${\bar{u}}^{xt}$}{$x$-space and time-averaged velocity in the $y$-direction}
\nomenclature[0u]{${\bar{u}}^{yt}$}{$y$-space and time-averaged velocity in the $x$-direction}
\nomenclature[0u]{${\bar{u}}^0$}{approximate convective velocity in the $x$-direction at current timestep}
\nomenclature[0v]{${\bar{v}}^{xy}$}{$x$ and $y$-space averaged velocity in the $t$-direction}
\nomenclature[0v]{${\bar{v}}^{xt}$}{$x$-space and time-averaged velocity in the $y$-direction}
\nomenclature[0v]{${\bar{v}}^{yt}$}{$y$-space and time-averaged velocity in the $x$-direction}
\nomenclature[0v]{${\bar{v}}^0$}{approximate convective velocity in the $y$-direction at current timestep}
%
%
\nomenclature[G]{$\Delta t$}{time step}
\nomenclature[G]{$\Delta x$}{node size in the $x$-direction}
\nomenclature[G]{$\Delta y$}{node size in the $y$-direction}
\nomenclature[G]{$\tau$}{half of the time step}
%
%
\nomenclature[D]{$Re$}{Reynolds number}
\nomenclature[D]{${Reu}_{i,j}$}{local Reynolds number}
%
\nomenclature[s]{$i$}{spatial index in the $x$-direction }
\nomenclature[s]{$j$}{spatial index in the $y$-direction}
\nomenclature[s]{$\ell$}{temporal index}
\nomenclature[s]{$t$}{transverse averaging in time}
\nomenclature[s]{$x$}{transverse averaging in the $x$-direction}
\nomenclature[s]{$y$}{transverse averaging in the $y$-direction}

%

{\printnomenclature[5em]}
\end{spacing}
%
%
%
%
\bibliography{references}
%
%
%
%
%
\newpage
\appendix
\section{MCCNIM coefficients}\label{app:A}
\renewcommand{\thesection}{\Alph{section}}
\subsection{MCCNIM coefficients for 1D Burgers’ equation}
The coefficients for the averaged flux equations (i.e., \eqn\eqref{eq:020} and \eqn\eqref{eq:021}) are as follows:
\begin{gather}
	\begin{split}
		A_{31}=\frac{2ae^{Reu_{i,\ell}}{Re}\add{({{\bar{u}}_{i,\ell}^0})^2}}{1-e^{Reu_{i,\ell}}+2ae^{Reu_{i,\ell}}Re{\bar{u}}_{i,\ell}^0} 
	\end{split}
\end{gather}
\begin{gather}
	\begin{split}
		A_{32}=\frac{-1+e^{Reu_{i,\ell}}-2ae^{Reu_{i,\ell}}Re{\bar{u}}_{i,\ell}^0+2a^2e^{Reu_{i,\ell}}{Re}^2\add{({{\bar{u}}_{i,\ell}^0})^2}}{Re{\bar{u}}_{i,\ell}^0\left(-1+e^{Reu_{i,\ell}}-2ae^{Reu_{i,\ell}}Re{\bar{u}}_{i,\ell}^0\right)} 
	\end{split}
\end{gather}
\begin{gather}
	\begin{split}
		A_{51}=\frac{2a{Re}\add{({{\bar{u}}_{i,\ell}^0})^2}}{1-e^{Reu_{i,\ell}}+2aRe{\bar{u}}_{i,\ell}^0} 
	\end{split}
\end{gather}
\begin{gather}
	\begin{split}
		A_{52}=\frac{1-e^{Reu_{i,\ell}}+2aRe{\bar{u}}_{i,\ell}^0+2a^2{Re}^2\add{({{\bar{u}}_{i,\ell}^0})^2}}{Re{\bar{u}}_{i,\ell}^0\left(1-e^{Reu_{i,\ell}}+2aRe{\bar{u}}_{i,\ell}^0\right)} 
	\end{split}
\end{gather}
The coefficients of the final set of algebraic equations (i.e., \eqn\eqref{eq:035} - \eqn\eqref{eq:037}) for the general interior node in one-dimensional Burgers’ case are given as follows:
\begin{gather}
	\begin{split}
	 F_{31}=\frac{A_{32}\left(A_{51,i+1}-{\bar{u}}_{i,\ell}^0\right)}{2a\left(A_{31}-A_{51,i+1}\right)}+\frac{A_{52}\left(A_{31,i-1}-{\bar{u}}_{i,\ell}^0\right)}{2a\left(A_{31,i-1}-A_{51}\right)} 
	\end{split}
\end{gather}
\begin{gather}
	\begin{split}
	F_{32}=\frac{A_{32,i-1}\left(-A_{51}+{\bar{u}}_{i,\ell}^0\right)}{2a\left(A_{31,i-1}-A_{51}\right)} 
	\end{split}
\end{gather}
\begin{gather}
	\begin{split}
		F_{33}=\frac{A_{52,i+1}\left(-A_{31}+{\bar{u}}_{i,\ell}^0\right)}{2a\left(A_{31}-A_{51,i+1}\right)} 
	\end{split}
\end{gather}
\begin{gather}
	\begin{split}
		F_{34}=\frac{A_{31}\left(A_{51,i+1}-{\bar{u}}_{i,\ell}^0\right)}{2a\left(A_{31,i-1}-A_{51}\right)}+\frac{A_{51}\left(A_{31,i-1}-{\bar{u}}_{i,\ell}^0\right)}{2a\left(A_{31,i-1}-A_{51}\right)} 
	\end{split}
\end{gather}
\begin{gather}
	\begin{split}
		F_{35}=\frac{A_{31,i-1}\left(-A_{51}+{\bar{u}}_{i,\ell}^0\right)}{2a\left(A_{31,i-1}-A_{51}\right)} 
	\end{split}
\end{gather}
\begin{gather}
	\begin{split}
		F_{36}=\frac{A_{51,i+1}\left(-A_{31}+{\bar{u}}_{i,\ell}^0\right)}{2a\left(A_{31}-A_{51,i+1}\right)} 
	\end{split}
\end{gather}
\begin{gather}
	\begin{split}
		F_{51}=\frac{2\tau F_{31}}{1-2\tau F_{34}};\quad F_{52}=\frac{2\tau F_{32}}{1-2\tau F_{34}};\quad F_{53}=\frac{2\tau F_{33}}{1-2\tau F_{34}};\quad F_{54}=\frac{\tau}{1-2\tau F_{34}};\\ 
		F_{55}=\frac{\tau}{1-2\tau F_{34}} \quad  F_{56}=\frac{1}{1-2\tau F_{34}};\quad F_{57}=\frac{2\tau F_{35}}{1-2\tau F_{34}};\quad F_{58}=\frac{2\tau F_{36}}{1-2\tau F_{34}} 
	\end{split}
\end{gather}
The coefficients of the final set of algebraic equations (i.e., \eqn\eqref{eq:045} and \eqn\eqref{eq:046}) for the right boundary (Dirichlet boundary) node in one-dimensional Burger’s case are given as follows:
\begin{gather}
	\begin{split}
		F_{31}^R=\frac{A_{52}\left(A_{31,i-1}+{\bar{u}}_{i,\ell}^0\right)}{2a\left(A_{31,i-1}-A_{51}\right)}-\frac{A_{32}}{2a} 
	\end{split}
\end{gather}
\begin{gather}
	\begin{split}
		F_{32}^R=\frac{A_{32,i-1}\left(-A_{51}+{\bar{u}}_{i,\ell}^0\right)}{2a\left(A_{31,i-1}-A_{51}\right)} 
	\end{split}
\end{gather}
\begin{gather}
	\begin{split}
		F_{33}^R=0 
	\end{split}
\end{gather}
\begin{gather}
	\begin{split}
		F_{34}^R=\frac{A_{51}\left(A_{31,i-1}-{\bar{u}}_{i,\ell}^0\right)}{2a\left(A_{31,i-1}-A_{51}\right)}-\frac{A_{31}}{2a} 
	\end{split}
\end{gather}
\begin{gather}
	\begin{split}
		F_{35}^R=\frac{A_{31,i-1}\left(-A_{51}+{\bar{u}}_{i,\ell}^0\right)}{2a\left(A_{31,i-1}-A_{51}\right)} 
	\end{split}
\end{gather}
\begin{gather}
	\begin{split}
		F_{36}^R=\frac{\left(A_{31}-{\bar{u}}_{i,\ell}^0\right)}{2a} 
	\end{split}
\end{gather}
The coefficients, $F_5^R$’s retain the same form as those defined for the general node ($F_5$’s ), with the key difference being that these coefficients now depend on $F_3^R$'s instead of $F_3$. For instance, $F_{51}^R$ is given by $F_{51}^R=\frac{2\tau F_{31}^R}{1-2\tau F_{34}^R}$. Similarly, expressions for the coefficients at the left boundary node can be derived following the same approach. Additionally, the coefficients corresponding to the Neumann boundary condition can be obtained in a comparable manner.
\subsection{MCCNIM coefficients for 2D Burgers’ equation}
The coefficients for the two-dimensional Burgers’ equation, as detailed in \sect\ref{sec:2.2}, are presented here for clarity. These include ($A$’s, $B$’s, $F$’s, $G$’s, $L$’s, $M$’s, $N$’s). Notably:
\begin{itemize}
	\item 	The definitions of $A$’s and $F$’s remain identical to those in the one-dimensional case.
	\item The initial $B$’s coefficients, which are analogous to $A$’s, are given as follows:
\end{itemize}
\begin{gather}
	\begin{split}
		B_{31}=\frac{2be^{Rev_{i,j,\ell}}{Re}\add{({{\bar{v}}_{i,j,\ell}^0})^2}}{1-e^{Rev_{i,j,\ell}}+2be^{Rev_{i,j,\ell}}Re{\bar{v}}_{i,j,\ell}^0} 
	\end{split}
\end{gather}
\begin{gather}
	\begin{split}
		B_{32}=\frac{-1+e^{Rev_{i,j,\ell}}-2be^{Rev_{i,j,\ell}}Re{\bar{v}}_{i,j,\ell}^0+2b^2e^{Rev_{i,j,\ell}}{Re}^2\add{({{\bar{v}}_{i,j,\ell}^0})^2}}{Re{\bar{v}}_{i,j,\ell}^0\left(-1+e^{Rev_{i,j,\ell}}-2be^{Rev_{i,j,\ell}}Re{\bar{v}}_{i,j,\ell}^0\right)} 
	\end{split}
\end{gather}
\begin{gather}
	\begin{split}
		B_{51}=\frac{2b{Re}\add{({{\bar{v}}_{i,j,\ell}^0})^2}}{1-e^{Rev_{i,j,\ell}}+2bRe{\bar{v}}_{i,j,\ell}^0} 
	\end{split}
\end{gather}
\begin{gather}
	\begin{split}
		B_{52}=\frac{1-e^{Rev_{i,j,\ell}}+2bRe{\bar{v}}_{i,j,\ell}^0+2b^2{Re}^2\add{({{\bar{v}}_{i,j,\ell}^0})^2}}{Re{\bar{v}}_{i,j,\ell}^0\left(1-e^{Rev_{i,j,\ell}}+2bRe{\bar{v}}_{i,j,\ell}^0\right)} 
	\end{split}
\end{gather}
The coefficients of the final set of algebraic equations (i.e., \eqn\eqref{eq:081} - \eqn\eqref{eq:084}) for the two-dimensional Burgers’ equation:
\begin{gather}
	\begin{split}
	 G_{31}=\frac{B_{32}\left(B_{51,j+1}-{\bar{v}}_{i,j,\ell}^0\right)}{2b\left(B_{31}-B_{51,j+1}\right)}+\frac{B_{52}\left(B_{31,\ell-1}-{\bar{v}}_{i,j,\ell}^0\right)}{2b\left(B_{31,\ell-1}-B_{51}\right)} 
	\end{split}
\end{gather}
\begin{gather}
	\begin{split}
	G_{32}=\frac{B_{32,\ell-1}\left(-B_{51}+{\bar{v}}_{i,j,\ell}^0\right)}{2b\left(B_{31,\ell-1}-B_{51}\right)} 
	\end{split}
\end{gather}
\begin{gather}
	\begin{split}
		G_{33}=\frac{B_{52,j+1}\left(-B_{31}+{\bar{v}}_{i,j,\ell}^0\right)}{2b\left(B_{31}-B_{51,j+1}\right)} 
	\end{split}
\end{gather}
\begin{gather}
	\begin{split}
		G_{34}=\frac{B_{31}\left(B_{51,j+1}-{\bar{v}}_{i,j,\ell}^0\right)}{2b\left(B_{31,\ell-1}-B_{51}\right)}+\frac{B_{51}\left(B_{31,\ell-1}-{\bar{v}}_{i,j,\ell}^0\right)}{2b\left(B_{31,\ell-1}-B_{51}\right)} 
	\end{split}
\end{gather}
\begin{gather}
	\begin{split}
		G_{35}=\frac{B_{31,\ell-1}\left(-B_{51}+{\bar{v}}_{i,j,\ell}^0\right)}{2b\left(B_{31,\ell-1}-B_{51}\right)} 
	\end{split}
\end{gather}
\begin{gather}
	\begin{split}
		G_{36}=\frac{B_{51,j+1}\left(-B_{31}+{\bar{v}}_{i,j,\ell}^0\right)}{2b\left(B_{31}-B_{51,j+1}\right)} 
	\end{split}
\end{gather}
\begin{gather}
	\begin{split}
	\ L_{51}=\frac{\tau}{1+2\tau F_{34}+2\tau G_{34}};\ \ L_{52}=\frac{\tau}{1+2\tau F_{34}+2\tau G_{34}};\\ L_{53}=\frac{1}{1+2\tau F_{34}+2\tau G_{34}};\ \ L_{54}=\frac{2\tau}{1+2\tau F_{34}+2\tau G_{34}} 	\end{split}
\end{gather}
\begin{gather}
	\begin{split}
		M_{51}=\frac{2\tau F_{31}}{1+2\tau F_{34}+2\tau G_{34}};\ \ M_{52}=\frac{2\tau F_{32}}{1+2\tau F_{34}+2\tau G_{34}};\ \ M_{53}=\frac{2\tau F_{33}}{1+2\tau F_{34}+2\tau G_{34}}; \\
		M_{54}=\frac{2\tau F_{35}}{1+2\tau F_{34}+2\tau G_{34}};\ \ M_{55}=\frac{2\tau F_{36}}{1+2\tau F_{34}+2\tau G_{34}}; 
	\end{split}
\end{gather}
\begin{gather}
	\begin{split}
	N_{51}=\frac{2\tau G_{31}}{1+2\tau F_{34}+2\tau G_{34}};\ \ N_{52}=\frac{2\tau G_{32}}{1+2\tau F_{34}+2\tau G_{34}};\ \ N_{53}=\frac{2\tau G_{33}}{1+2\tau F_{34}+2\tau G_{34}}; \\
	N_{54}=\frac{2\tau G_{34}}{1+2\tau F_{34}+2\tau G_{34}};\ \ N_{55}=\frac{2\tau G_{35}}{1+2\tau F_{34}+2\tau G_{34}} 
	\end{split}
\end{gather}
\subsection{MCCNIM coefficients for approximated (node-averaged) convective velocity in both 1D and 2D Burgers’ equations}
The coefficients for the equation node-averaged velocity $u_{i,\ell}^0$ given by \eqn\eqref{eq:040} in the one-dimensional Burgers' case are expressed as follows:
\begin{gather}
	\begin{split}
	F_{71}=\frac{A_{32,i-1}}{2\left(A_{31,i-1}-A_{51}\right)} 
	\end{split}
\end{gather}
\begin{gather}
	\begin{split}
		F_{72}=\frac{A_{32}}{2\left(A_{31}-A_{51,i+1}\right)}-\frac{A_{52}}{2\left(A_{31,i-1}-A_{51}\right)} 
	\end{split}
\end{gather}
\begin{gather}
	\begin{split}
		F_{73}=-\frac{A_{52,i+1}}{2\left(A_{31}-A_{51,i+1}\right)} 
	\end{split}
\end{gather}
\begin{gather}
	\begin{split}
		F_{74}=\frac{A_{31,i-1}}{2\left(A_{31,i-1}-A_{51}\right)} 
	\end{split}
\end{gather}
\begin{gather}
	\begin{split}
	F_{75}=-\frac{A_{51}}{2\left(A_{31,i-1}-A_{51}\right)}+\frac{A_{31}}{2\left(A_{31}-A_{51,i+1}\right)} 
	\end{split}
\end{gather}
\begin{gather}
	\begin{split}
		F_{76}=-\frac{A_{51,i+1}}{2\left(A_{31}-A_{51,i+1}\right)} 
	\end{split}
\end{gather}
Similarly, the equation of the node-averaged velocities $u_{i,j,\ell}^0$ and $v_{i,j,\ell}^0$ given by \eqn\eqref{eq:093} in the two-dimensional Burgers' case are expressed as follows:
\begin{gather}
	\begin{split}
	G_{71}=\frac{B_{32,\ell-1}}{2\left(B_{31,\ell-1}-B_{51}\right)} 
	\end{split}
\end{gather}
\begin{gather}
	\begin{split}
		G_{72}=\frac{B_{32}}{2\left(B_{31}-B_{51,j+1}\right)}-\frac{B_{52}}{2\left(B_{31,\ell-1}-B_{51}\right)} 
	\end{split}
\end{gather}
\begin{gather}
	\begin{split}
		G_{73}=-\frac{B_{52,j+1}}{2\left(B_{31}-B_{51,j+1}\right)} 
	\end{split}
\end{gather}
\begin{gather}
	\begin{split}
	G_{74}=\frac{B_{31,\ell-1}}{2\left(B_{31,\ell-1}-B_{51}\right)} 
	\end{split}
\end{gather}
\begin{gather}
	\begin{split}
		G_{75}=-\frac{B_{51}}{2\left(B_{31,\ell-1}-B_{51}\right)}+\frac{B_{31}}{2\left(B_{31}-B_{51,j+1}\right)} 
	\end{split}
\end{gather}
\begin{gather}
	\begin{split}
		G_{76}=-\frac{B_{51,j+1}}{2\left(B_{31}-B_{51,j+1}\right)} 
	\end{split}
\end{gather}
The coefficients of $F_7$’s in the two-dimensional Burger's equation will mirror those used in the one-dimensional Burger's equation.
%
%
%
%
%
%
%
%
%
\end{document}